\documentclass[10pt]{amsart}

\newcommand{\ntop}[2]{\genfrac{}{}{0pt}{1}{#1}{#2}}

\let\newpf\proof \let\proof\relax 
\newenvironment{pf}{\newpf[\proofname]}{\qed\endtrivlist}

\def\be{\begin{equation}}
\def\ee{\end{equation}}
\def\bm{\begin{pmatrix}}
\def\em{\end{pmatrix}}

\def\RDC{\mathrm{{RDC}}}
\def\P{\mathbb{{P}}}

\def\op{\overline \partial}

\def\ba{{\begin{align}}}
\def\ea{{\end{align}}}

\def\Ad{{\mathrm {Ad}}}
\def\pa{\partial}

\def\g{{\gamma}}

\def\SO{{\mathrm {SO}}}
\def\GL{{\mathrm {GL}}}
\def\so{{\mathrm {so}}}
\def\Sr{{\mathbb {S}}}

\def\Sl{\mathrm {sl}}

\def\SU{{\mathrm {SU}}}

\def\H{{\mathbb H}}

\def\0{{\mathbf 0}}

\def\cal{\mathcal}

\def\SL{{\mathrm {SL}}}

\def\PSL{\mathrm{{PSL}}}

\newtheorem{thm}{Theorem}[section]
\newtheorem{cor}[thm]{Corollary}

\newtheorem{lemma}[thm]{Lemma}

\newtheorem{prop}[thm]{Proposition}

\theoremstyle{remark}
\newtheorem{rem}{Remark}[section]

\newtheorem{problem}{Problem}[section]

\numberwithin{equation}{section}

\def \bn {\hfill \\ \smallskip\noindent}

\theoremstyle{definition}

\def\proof{\bn {\bf Proof.} }

\def\note#1
{\marginpar
{\tiny $\leftarrow$
\par
\hfuzz=20pt \hbadness=9000 \hyphenpenalty=-100 \exhyphenpenalty=-100
\pretolerance=-1 \tolerance=9999 \doublehyphendemerits=-100000
\finalhyphendemerits=-100000 \baselineskip=6pt
#1}\hfuzz=1pt}

\newcommand{\bignote}[1]{\begin{quote} \sf #1 \end{quote}}

\newcommand{\inter}{\operatorname{int}}
\renewcommand{\mod}{\operatorname{mod}}

\newcommand{\id}{\operatorname{id}}

\newcommand{\BB}{{\cal B}}

\newcommand{\MM}{{\cal M}}

\newcommand{\RR}{{\cal R}}

\newcommand{\C}{{\mathbb C}}
\newcommand{\D}{{\mathbb D}}

\newcommand{\Q}{{\mathbb Q}}
\newcommand{\R}{{\mathbb R}}
\newcommand{\T}{{\mathbb T}}

\newcommand{\U}{{\mathbb U}}

\newcommand{\Z}{{\mathbb Z}}

\def\B0{{\bold{0}}}

\catcode`\@=12

\def\Empty{}
\newcommand\oplabel[1]{
  \def\OpArg{#1} \ifx \OpArg\Empty {} \else
  	\label{#1}
  \fi}
		
\newcommand{\comm}[1]{}
\newcommand{\comment}[1]{}

\begin{document}


\title[Monotonic cocycles]{Monotonic cocycles}

\author{Artur Avila and Rapha\"el Krikorian}

\address{
CNRS UMR 7586, Institut de Math\'ematiques de Jussieu - Paris Rive Gauche,
B\^atiment Sophie Germain, Case 7012, 75205 Paris Cedex 13, France
\&
IMPA, Estrada Dona Castorina 110, 22460-320, Rio de Janeiro, Brazil
}

\email{artur@math.jussieu.fr}

\address{ Laboratoire de Probabilit\'es et Mod\`eles al\'eatoires\\
  Universit\'e Pierre et Marie Curie--Boite courrier 188\\
  75252--Paris Cedex 05\\
  France 
}

\email{raphael.krikorian@upmc.fr}

\begin{abstract}

We develop a ``local theory'' of multidimensional quasiperiodic
$\SL(2,\R)$ cocycles which are not homotopic to a constant.  It describes a
$C^1$-open neighborhood of cocycles of rotations and
applies irrespective of arithmetic
conditions on the frequency, being
much more robust than the local theory of
$\SL(2,\R)$ cocycles homotopic to a constant.  Our analysis is centered
around the notion of monotonicity with respect to some dynamical variable.
For such {\it monotonic cocycles},
we obtain a sharp rigidity result, minimality of the projective action,
typical nonuniform hyperbolicity, and
a surprising result of smoothness of the Lyapunov exponent (while
no better than H\"older can be obtained in the case of
cocycles homotopic to a constant, and only under arithmetic restrictions).
Our work is based on complexification ideas, extended ``\`a la Lyubich''
to the smooth setting (through the use of asymptotically holomorphic
extensions).  We also develop a counterpart of this theory centered around
the notion of monotonicity with respect to a parameter variable, which
applies to the analysis of $\SL(2,\R)$ cocycles over more general dynamical
systems and generalizes key aspects of Kotani Theory.  We conclude with a
more detailed discussion of one-dimensional monotonic cocycles,
for which results about rigidity and typical nonuniform hyperbolicity can be
globalized using a new result about convergence of renormalization.

\end{abstract}

\setcounter{tocdepth}{1}

\maketitle

\tableofcontents

\section{Introduction}

Let $f:X \to X$ be a homeomorphism of a compact metric space, preserving a
probability measure $\mu$.  Given a map $A \in C^0(X,\SL(2,\R))$ the
skew-product map on $X \times \R^2$
given by $(x,w) \mapsto (f(x),A(x) \cdot w)$, denoted by $(f,A)$, is called
an $\SL(2,\R)$ cocycle over $f$.

We will be particularly interested in {\it quasiperiodic cocycles} where
$X=\R^d/\Z^d$ and $f$ is a translation, $f(x)=x+\alpha$ for some $\alpha \in
\R^d$ and $\mu$ is Lebesgue measure.
Since in this case $f$ is a diffeomorphism of a manifold, it makes sense to
consider quasiperiodic cocycles with various degrees of smoothness.

Quasiperiodic cocycles $(f,A)$ have been primarily studied in the case where
$A$ is homotopic to a constant (in large part because this is the situation
arising in the consideration of quasiperiodic Schr\"odinger operators).  One
important aspect was the development of a {\it local theory}, starting with
the KAM based work of Dinaburg-Sinai \cite {DS}.  This local theory concerns
perturbations of the simplest cocycles homotopic to a constant, which are
just the constant ones.  The development of the
KAM approach involves, as usual, non-resonance assumptions which here are
coded in arithmetic conditions involving the frequency vector but also the
{\it fibered rotation number}, and a key achievement, due to Eliasson \cite
{E}, was the development of a local theory covering all cocycles with Diophantine
frequency vector.  Except for the one-dimensional case, where considerably
more can be achieved by a range of techniques
(both non-KAM as in \cite {BJ2}, \cite {AJ2}, and non-standard KAM,
\cite {AFK}), the work of Eliasson remains
basically the best description of the local theory in the case of cocycles
homotopic to a constant.

One of our goals here is to develop a local theory of cocycles non-homotopic
to a constant, covering (in the ergodic case)
perturbations of the simplest cocycles arising in
this case, which are the $\SO(2,\R)$-valued ones.  As it will turn out, the
theory we develop is considerably more robust then the usual one, in several
respects.  For instance, the frequency vector plays no role at all in our
considerations, and we are able to treat quite low regularity ($C^1$)
perturbations.  Moreover, several of our conclusions are in a sense also
much stronger, and even surprising from the point of view of the intuition
developed in the case of cocycles homotopic to a constant.  Specific
questions addressed in this paper concern the
regularity of the Lyapunov exponents,
rigidity arising from zero Lyapunov exponents, and minimality of the
associated projective action.

Our local theory centers around the crucial property of monotonicity
with respect to some phase (dynamical)
variable, a kind of twist condition that arises naturally (in
the ergodic case) for
$\SO(2,\R)$-valued cocycles not homotopic to a constant.
There is a counterpart, which actually precedes logically the analysis of
monotonic cocycles, which describes the consequences of monotonicity
with respect to parameter variables, and works for cocycles over general
dynamical systems.  Our results, abstract and
generalize key aspects of the theory of Schr\"odinger cocycles
(particularly Kotani Theory), in particular to the case of
non-analytic dependence on parameters.  Besides being fundamental to our
analysis of monotonic cocycles, we would also like to point out that
the larger flexibility afforded by this theory has been recently applied
back to address problems about the Schr\"odinger case \cite {A1}.

Our third focus in this paper concerns the analysis of
one-dimensional quasiperiodic cocycles non-homotopic to a constant from the
{\it global} point of view.  As in \cite {AK}, the basic plan is to reduce
global questions to local ones by renormalization.  To this end, we prove
a ``convergence of renormalization'' result that guarantees that, under a
natural (from the point of view of parameter analysis) hypothesis that
renormalizations eventually become monotonic.  As a consequence, we obtain a
global $L^2$-rigidity result, and conclude that typical cocycles
non-homotopic to a constant are nonuniformly hyperbolic.

We will next present more formally some key results of each of the three
topics mentioned above.

\subsection{Monotonic cocycles}

\comm{
The iterates of a cocycle $(f,A)$, we denote its iterates
We first recall the basic
definitions of the Lyapunov exponent
\be
L(f,A)=\lim_{n \to \infty} \frac {1} {n} \int \|A_n(x)\| dx,
\ee
and of conjugacy between cocycles $(f,A)$ and
}

Below we
will use the notation $(f^n,A_n)$ for the
$n$-iterate of the cocycle $(f,A)$: thus if $n \geq 1$, $A_n(x)=A(f^{n-1}(x)
) \cdots A(x)$.  Let us also recall basic definitions of the Lyapunov
exponent
\be
L(f,A)=\lim_{n \to \infty} \frac {1} {n} \int \ln \|A_n(x)\| dx,
\ee
and of a conjugacy between cocycles $(f,A)$ and $(f,A')$, which is given by
a map $B:X \to \SL(2,\R)$ satisfying
\be
A'(x)=B(f(x)) A(x) B(x)^{-1}.
\ee

Now, and for the remaining of this section, fix $d \geq 1$ and
let $X=\R^d/\Z^d$.  For $\alpha \in \R^d$ denote by
$f_\alpha:X \to X$ the map $f_\alpha(x)=x+\alpha$.  Since such dynamics are
regular, it makes sense to speak of regular cocycles and
regular conjugacies.

\comm{
Once $\alpha$ is fixed,
any $A \in C^0(X,\SL(2,\R))$ will determine a cocycle, and we will denote
the $n$-th iterate of $(f_\alpha,A)$ by $(f_\alpha^n,A_n)$: thus, if $n \geq
1$, $A_n(x)=A(x+(n-1) \alpha) \cdots A(x)$.
We recall basic
definitions of the Lyapunov exponent
\be
L(f_\alpha,A)=\lim_{n \to \infty} \frac {1} {n} \int \|A_n(x)\| dx,
\ee
and of conjugacy between cocycles $(f,A)$ and
}

We say that $A \in C^1(X,\SL(2,\R))$ is {\it monotonic} if there exists
some $w \in \R^d$ such that for every $x \in X$ and
$y \in \R^2 \setminus \{0\}$,
(any determination of) the argument of $A(x+t w) \cdot y$, $t \in \R$,
has positive derivative with respect to $t$.

This condition clearly determines a $C^1$-open
subset $\MM^1$ of $C^1(X,\SL(2,\R))$.
Notice that the monotonicity condition only makes reference to $A$ and is
thus independent of a frequency vector.  For this reason, it is in
particular not invariant by conjugacy.  Given a frequency vector $\alpha$,
it is thus natural to define a set $\MM^1_\alpha \subset C^1(X,\SL(2,\R))$
consisting of all $A$ for which there exists $n \geq 1$ such that
$(f_\alpha,A)^n$ is $C^1$-conjugated to some
monotonic cocycle.  We call cocycles $(f_\alpha,A)$ with $A \in
\MM^1_\alpha$ {\it premonotonic}.
It is also clear
that cocycles homotopic to a constant can never be premonotonic.\footnote
{It is somewhat delicate (in the ergodic case)
to construct examples of cocycles not homotopic to
a constant which are not premonotonic.  A non-negligible (positive measure
on parameters) set of such examples can be obtained by forcing a certain
behavior of the ``critical points'' that arise in the
approach of Lai-Sang Young \cite {Y1} (of Benedicks-Carleson \cite {BC}
flavor).  We will come back to this issue elsewhere.}
On the other hand, if $(f_\alpha,A)$ is $C^0$ conjugate to a
{\it cocycle of rotations}
(that is, to an $\SO(2,\R)$-valued one), and $f_\alpha$ is ergodic,
then $(f_\alpha,A)$ is premonotonic as long as
$A$ is $C^1$ and not homotopic to a constant.

Obviously any cocycle of rotations, or merely measurably
conjugate to such, must have a zero Lyapunov exponent.
Our first result shows that for regular (pre-) monotonic cocycles, one can go the
other way around:

\begin{thm} \label {l=0cr}

Let $(f_\alpha,A)$ be $C^r$, $r=\infty,\omega$ and premonotonic.  If
$L(f_\alpha,A)=0$ then $(f_\alpha,A)$ is $C^r$ conjugate to an
$\SO(2,\R)$-valued cocycle.

\end{thm}

By the previous discussion, Theorem \ref {l=0cr}
contains in it a rigidity result.

\begin{cor} \label {c0cr}

Let $f_\alpha$ be ergodic.  If
$(f_\alpha,A)$ be $C^r$, $r=\infty,\omega$ is non-homotopic to a constant
and $C^0$ conjugate to rotations then it is $C^r$-conjugate to rotations.

\end{cor}

Next we look at the dependence of the Lyapunov exponent.  We recall first
that in the case of cocycles homotopic to a constant, simple examples show
that one should never expect more (as far as the modulus of continuity is
concerned) than $1/2$-H\"older regularity, and Eliasson's
local theory does provide precisely this regularity
\cite {sana}, but only under arithmetic conditions on the
frequency vector.  In fact, it is easy to see
that in general no statement about the modulus of continuity can be made.
Even continuity is only known in the analytic case (a deep global
result of Bourgain \cite {B}).
Thus the following results were very surprising to the authors:

\begin{thm} \label {lyapanal}

Let $(f_\alpha,A_s)$ be a one-parameter analytic
family of analytic premonotonic cocycles.  Then $s \mapsto L(f_\alpha,A_s)$
is analytic.

\end{thm}

\begin{thm} \label {lyapsmooth}

Let $(f_{\alpha(s)},A_s)$ be a one-parameter $C^\infty$
family of $C^\infty$ premonotonic cocycles.
Then $s \mapsto L(f_{\alpha(s)},A_s)$ is $C^\infty$.

\end{thm}

\begin{rem}

Those results are all the more striking in view of the recent discovery by
Wang-You \cite {WY} of smooth
one-frequency cocycles at which the Lyapunov exponent is discontinuous.

\end{rem}

Since the Lyapunov exponent is a regular function, its zero set can be
expected to be some kind of variety, hence Lyapunov exponents should be rare
unless the equation $L=0$ is very degenerate.
In fact, since the Lyapunov exponent
can not become negative, $DL$ must be zero whenever $L=0$.  We are however
able to show non-degeneracy of $D^2 L$, which implies:

\begin{thm} \label {premonotonicpositive}

For fixed $\alpha$, a typical $C^r$, $r=\omega,\infty$,
premonotonic cocycle over $f_\alpha$ has positive Lyapunov
exponent.\footnote {The set of premonotonic cocycles with zero Lyapunov
exponent has infinite codimension when the number of frequencies $d$ is at
least $2$, and finite codimension when $d=1$.  We will come back to this
issue elsewhere.}

\end{thm}

Though we do obtain several other results, particularly addressing less
regular situations, we would like to conclude our discussion, at this
introduction, with a
result of different flavor.  Given a cocycle $(f,A)$, we may define its
projective action, which is just the projectivized skew-product on
$X \times \P\R^2$.  The topological dynamics of the projective action is a
very interesting subject in itself (see for instance \cite {Bj1}, \cite {Bj2}, \cite
{BjJo}, \cite {J2}).
Here we prove:

\begin{thm} \label {preminimal}

If $f_\alpha$ is ergodic and $(f_\alpha,A)$ is premonotonic with $A \in
C^{1+\epsilon}$, then the projective action is minimal.

\end{thm}

Let us point out that dynamical notions of monotonicity also make sense for
dynamical systems which are not strict translations, such as the
skew-shifts, and some of our results can be carried to a larger generality,
see Remark \ref {skewshift}.

\subsection{Monotonicity in the parameter space}

We return now to the consideration of more general dynamics $f:X \to X$. 
This time we will be interested in parametrized families of cocycles
$(f,A_\theta)$, and we will often assume some base regularity of this
dependence (just with respect to $\theta$, since nothing beyond continuity
can be made sense with respect to the dynamical variable under our
hypothesis).  Moreover, we will require the dependence of
$A_\theta$ on $\theta$ displays monotonicity: assuming that
$\theta \mapsto A_\theta$ is $C^1$, this means that for every $x \in X$,
$y \in \R^2 \setminus \{0\}$,
(any determination of) the argument of $\theta \mapsto
A_\theta(x) \cdot y$ has positive derivative.

Let us consider two key examples where monotonicity arises.

Recall that in the dynamical approach to ergodic Schr\"odinger operators,
the basic object considered is a one-parameter family of cocycles, depending
on a parameter $E$, of the form $A^{(E)}(x)=\bm E-v(x) & -1 \\ 1 & 0
\em$.  Though the
family $(f,A^{(E)})$ is not monotonic in $E$, its second
iterate is (which in fact is just as good for our purposes).  This family
has of course been considered intensively and much is known about it:
Kotani Theory (\cite {kotani}, \cite {S}, and more dynamically \cite {CJ}),
in particular, gives much dynamical information about the set
of parameters where the Lyapunov exponent vanishes.

Another family type that has been considered is of the
form $\theta \mapsto R_\theta A$, where $A$ is arbitrary and $R_\theta \in
\SO(2,\R)$ is the rotation of angle $2 \pi \theta$.  This family displays
obvious monotonicity.  Some global aspects of this family were first
exploited in \cite {H} to yield examples of cocycles with positive Lyapunov
exponents: the average Lyapunov exponent, with respect to $\theta$,
is zero if and only if $A$ is a
cocycle of rotations.  In fact, a later refinement \cite {AB} shows that
\be \label {abformula}
\int_{\R/\Z} L(f,R_\theta A) d\theta=\int_X \ln \frac {\|A\|+\|A\|^{-1}} {2}
d\mu(x),
\ee
so the average Lyapunov exponent depends on $A$ through a very simple
formula.  The dynamical aspects of
Kotani Theory have also been extended to such families.

Both examples we mentioned have in common, besides (some) monotonicity,
a very nice global behavior of the {\it holomorphic}
dependence on $\theta$ when $\theta$ is complex.  Here we will
show that a theory can be constructed without taking into account global
aspects: in fact even analyticity can be bypassed.  But complexification is
still fundamental, and what allows us to consider the smooth case is the use
of {\it asymptotically holomorphic} extensions, a tool first used in
dynamics by Lyubich \cite {Ly1}, in the context of unimodal maps.

In doing this, our key motivation has been the understanding of monotonic
cocycles.  In fact, if $(f_\alpha,A)$ is monotonic in the dynamical sense,
then one can construct a monotonic family $(f_\alpha,A_\theta)$ by setting
$A_\theta=A(x+\theta w)$, for some $w \in \R^d$.  When we change the
parameter, we are
not really changing the dynamics (merely the coordinates), but we do get
something non-trivial out of it, by applying the parameter results we will
obtain.  For this purpose, we will obtain analogous of
Theorems \ref {l=0cr} and \ref {lyapanal}. 
Those parametrized versions correspond respectively to
a non-global version of (\ref {abformula}) and
to a well known result of Kotani Theory (see \cite {CJ}).
Instead of presenting formal versions of those here,
we prefer to mention a different
Kotani-type application.  Let us say that $(f,A)$ is $L^2$-conjugate to
rotations if it is measurably conjugate and the conjugacy $B$
satisfies $\int_X \|B(x)\|^2 d\mu(x)<\infty$.

\begin{thm} \label {l2ae}

Let $A_\theta \in C^0(X,\SL(2,\R))$, be a one-parameter family which is
monotonic and $C^{2+\epsilon}$ in $\theta$.  Then for almost every
$\theta$, if $L(f,A_\theta)=0$ then $(f,A_\theta)$ is $L^2$-conjugate to
rotations.

\end{thm}

As we will see in the next section, $L^2$-conjugacy to rotations is a
fundamental hypothesis in renormalization theory of one-dimensional
quasiperiodic cocycles, so in a sense this
result enlarges the set of families along which parameter exclusion
arguments can be made before applying renormalization.  The ability to
analyze in this way arbitrary monotonic deformations turns out to be
relevant even if one is ultimately interested in the
Schr\"odinger case, see \cite {A1}.

\subsection{One-dimensional quasiperiodic cocycles non-homotopic to a
constant}

We consider again the quasiperiodic case, but now restrict attention to the
one-dimensional case.  In this section $X=\R/\Z$ and $f_\alpha$ will always
denote an ergodic translation (thus $\alpha \in \R \setminus \Q$).

Renormalization is a classical tool in the analysis of
diffeomorphisms of the circle, where it can be used to reduce global
questions to local ones \cite {KS}.
Application of renormalization ideas to the case of quasiperiodic cocycles
has also proved to be fruitful though in this case the renormalization
operator does not always lead to the local situation due to the
possible presence of positive Lyapunov exponents.  It should thus be
basically considered as a
tool to explore cocycles with zero Lyapunov exponent; see \cite{K1} for the case of $SU(2)$-valued cocycle.  However, zero
Lyapunov exponents are not a sufficient condition to achieve the
global-local reduction.\footnote {Particularly, the analysis of the spectrum
of the critical almost Mathieu operator, where the Lyapunov exponent is
still zero, does not reduce to the local situation.}

In \cite {AK}, it is shown that the existence of an $L^2$-conjugacy
to rotations is enough to guarantee {\it precompactness} of the
renormalization operator, which is used to extract limits which are
cocycles of rotations (up to constant conjugacy).  Using Kotani Theory and
the local (KAM) description of cocycles homotopic to a constant, this yields
a dicothomy (under suitable regularity requirements and arithmetic
conditions) for typical energies: the associated Schr\"odinger cocycle has
either a positive Lyapunov exponent, or it is conjugate to a
cocycle of rotations.\footnote {Since this result indeed assumes at least
a Diophantine condition on $\alpha$, and the cocycle is homotopic to a
constant, the conclusion is in fact equivalent to the existence of a
conjugacy to a {\it constant} cocycle of rotations.}

Though cocycles of rotations which are not homotopic
to a constant are premonotonic {\it if the basis is ergodic},
renormalization affects the base dynamics and in particular may lead to
non-ergodic limits.

In order to obtain more precise results, we prove here
that the limits of renormalization are in fact of a very special kind,
namely of the form $x \mapsto R_{\theta+|\deg| x}$, where $\deg$ is the
topological degree.  In particular, if $|\deg| \neq 0$, one does reach
monotonicity.  We conclude the following global rigidity result:

\begin{thm} \label {globall2}

Let $(f_\alpha,A)$, $\alpha \in \R \setminus \Q$,
be $C^r$, $r=\omega,\infty$, and
non-homotopic to a constant.  If $(f_\alpha,A)$
is $L^2$-conjugate to rotations, then
it is $C^r$-conjugate to rotations.

\end{thm}

Combined with Theorems \ref {l2ae} and \ref {premonotonicpositive} we conclude:

\begin{thm}

For fixed $\alpha \in \R \setminus \Q$,
a typical $C^r$, $r=\omega,\infty$,
cocycle over $f_\alpha$ which is not homotopic to a constant
has positive Lyapunov exponent.

\end{thm}

\begin{rem}

The first result about the {\it existence} of positive Lyapunov exponents
for cocycles non-homotopic to a constant was obtained in \cite {H}, for
cocycles of a very specific form (with respect to their global
holomorphic extensions, say
\be \label {lam}
x \mapsto \bm \lambda&0\\0&\lambda^{-1} \em R_x,
\ee
$\lambda>1$ (irrespective of the frequency).  Later
on, Lai-Sang Young \cite {Y1}
showed that for more general cocycles, but under a largeness
condition, for instance
perturbations of (\ref {lam}) with $\lambda>>1$,
positive Lyapunov exponents have large
probability with respect to the choice of the frequency vector.  The method
used by Young, quite different from ours,
of Benedicks-Carleson type \cite {BC}, is based on an inductive scheme
which loses control of a positive measure set of parameters, and needs some
initial condition to start (creating the need for a largeness assumption).

\end{rem}

We would like to point out that our more precise results about
convergence of renormalization have been recently applied \cite {AFK}
to the case of cocycles homotopic to a constant: together with
developments in the local theory, it yields the basic \cite {AK}
dichotomy without arithmetic conditions.

\subsection{Structure of the paper}

In section \ref {sec:2}, we analyze the consequences of monotonicity with respect to
parameters.  We start with some aspects of the dynamics of certain
$\SL(2,\C)$ cocycles, whose action by M\"obius transformations
preserve the upper half plane when going forward, but not necessarily
backwards.  We also discuss the crucial notion of
{\it variation of the fibered rotation number}, which is necessary even to
formulate several key results.  After a few simple applications of the
complexification idea in the analytic case, we describe an asymptotically
holomorphic framework that allows us to exploit the basic monotonicity
phenomenon, and we carry out the basic computations of Kotani Theory
in this setting.

Basically, monotonicity is used to guarantee that when the
parameter turns complex, the dynamics becomes {\it uniformly hyperbolic},
and everything depends nicely on parameters. 
The focus is thus to recover some information when the imaginary
part approaches zero.  In the analytic case, this is done by appealing to
theorems of complex analysis (such as Fatou's Theorem on existence of
non-tangential limits).  In the asymptotically holomorphic setting, we would
like to show that the discrepancy from holomorphicity corresponds to a
regular correction, say, of the non-tangential limit.  There are competing
factors though: while the cocycle behaves ``more holomorphically'' near real
parameters, it also behaves ``less uniformly hyperbolic''.  A key point is
thus to give an estimate of the resulting regularity
(in practice giving up some derivatives in the process).  After this, we are
in good shape to collect results such as Theorem \ref {l2ae}.

In section \ref {sec:3}
we then move on to the analysis of monotonic cocycles.  We
prove easily Theorems \ref {l=0cr} and \ref {lyapanal}
using the results of section \ref {sec:2}, and proceed to look at
Theorem \ref {lyapsmooth}, whose proof does not really fall in the same context:
since the dynamics changes we have to revisit the estimates
regarding the
competition between asymptotic holomorphicity and uniform hyperbolicity,
incorporating this additional parameter.  Using a somewhat different
approach, we next discuss results in low
regularity, such as continuity of the Lyapunov exponent in the Lipschitz
category.  We come back to a more regular situation in the proof of
Theorem \ref {preminimal}, where we apply
the Schwarz Reflection Principle, or rather, use that it can not be
applied.  We conclude with an estimate on the second derivative of the
Lyapunov exponent, which implies Theorem \ref {premonotonicpositive}.

In section \ref {sec4} we specify further to the one-dimensional case, but now with a
global focus.  We introduce formally the renormalization operator and
explain how ``convergence of renormalization'' combined with the local
theory indeed implies Theorem \ref {globall2}.
We conclude with a proof of convergence of renormalization.

We include two appendices.  The first discusses an estimate about the
conformal barycenter which is used when taking limits of $L^2$-conjugacies
to rotations, which is used to avoid unnecessary parameter exclusions.  The
second gives a proof of transitivity of the projective action for
quasiperiodic cocycles non-homotopic to a constant.

{\bf Acknowlegements:} R.K. would like to thank the hospitality of IMPA. 
This research was partially conducted during the period A.A. was a
Clay Research Fellow.  We are grateful to Zhenghe Zhang for detailed
comments about a preliminary version of this paper.

\comm{
\comm{
This result also contains a rigidity theorem: if $(f,A)$ were merely
$C^0$-conjugate to rotations then
}

The theory developed

\begin{thm}

Let $f:X \to X$ be an ergodic translation of the torus $X=\R^d$,
and let $A_s \in C^\omega(X,\SL(2,\R))$ be an analytic family.

Let $d \geq 1$ and let $f:\R^d/\Z^d \to \R^d/\Z^d$ be a translation.  A
quasiperiodic $\SL(2,\R)$ cocycle
}

\comm{
A {\it $d$-dimensional quasiperiodic $\SL(2,\R)$ cocycle} is a
pair $(\alpha,A)$ where $\alpha \in \R^d$ and $A \in
C^0(\R^d/\Z^d,\SL(2,\R))$ (we call $d$ the {\it number of frequencies}).
A cocycle should be viewed as a
{\it skew-product}:
\begin{align} \label {skpro}
(\alpha,A):\R^d/\Z^d \times \R^2 &\to \R/\Z \times \R^2\\
\nonumber
(x,w) &\mapsto (x+\alpha,A(x) \cdot w).
\end{align}
We call $\alpha$ the {\it frequency} of the cocycle.  We say that $\alpha$
is irrational if $x \mapsto x+\alpha$ acts minimally on $\R^d/\Z^d$.

There is a fairly developed theory of cocycles homotopic to a constant,
especially in the one-dimensional case. 
This is partially motivated by the theory of Schr\"odinger operators, which
leads to the study of a certain family of such cocycles. 

Our aim in this paper is to develop the theory in the case of cocycles not
homotopic to the identity.  It turns out that such a theory is quite rich in
its own right.  Most of our work will center on the construction of a
``local theory'' which is completely different, much more robust, and
sometimes much more complete than that of Schr\"odinger cocycles.
Indeed, the results we obtain are quite surprising: for
instance, while in the local theory of Schr\"odinger cocycles,
KAM schemes play a determinant role, arithmetic properties of the frequency
turn out to be irrelevant for all applications considered here.

While this local theory is the main novelty of this work, our original
motivation was to obtain global results from local results
via renormalization in the one-dimensional case.
We will start by discussing those.
}

\comm{
\subsection{Global results: typical non-uniform hyperbolicity}

The Lyapunov exponent of $(\alpha,A)$ is defined as
\be
L(\alpha,A)=\lim \frac {1} {n} \int_{\R^d/\Z^d} \ln \|A_n(x)\| dx \geq 0,
\ee
where $A_n(x)=\prod_{j=n-1}^0 A(x+j\alpha)=A(x+(n-1)\alpha) \cdots A(x)$
(we will keep the dependence on $\alpha$ implicit).

A key property of a cocycle (and other dynamical systems)
is whether it has a positive Lyapunov exponent:
this is a very good starting point to a description of the dynamics.

In the case of cocycles homotopic to a constant, there is a severe
obstruction for a smooth cocycle $(\alpha,A)$ to have a
positive Lyapunov exponent, namely, to be
conjugate to a constant elliptic matrix $A_* \in \SL(2,\R)$:
\be
A(x)=B(x+\alpha) A_* B(x)^{-1},
\ee
where $B:\R^d/\Z^d \to \SL(2,\R)$ is smooth.
We say that this obstruction is severe because it
is reflected on a positive measure set in parametrized families (this is a
consequence of KAM theory, see \cite {eliasson} for the most sophisticated
results).

In \cite {AK} the following result was proved.  Let
\be
R_\theta=
\begin{pmatrix}
\cos 2\pi\theta & -\sin 2\pi\theta\\
\sin 2\pi\theta & \cos 2\pi\theta
\end{pmatrix}.
\ee

\begin{thm}[\cite {AK}]

Let $\alpha \in \RDC$ (a subset of $\R$ of full Lebesgue measure).  For every
$A \in C^\omega(\R/\Z,\SL(2,\R))$ which is homotopic to the identity,
and for almost every $\theta \in \R/\Z$, eiher $L(\alpha,R_\theta A)>0$ or
$(\alpha,R_\theta A)$ is $C^\omega$-conjugate to a constant elliptic matrix.

\end{thm}

In other words, {\it a typical one-dimensional
quasiperiodic cocycle which is homotopic
to a constant is either non-uniformly hyperbolic or conjugate to a
constant.}  In this result, typical corresponded both to a full measure
set of frequencies and a ``full measure'' set of cocycles.



In the case we are interested here, of quasiperiodic
cocycles non-homotopic to the identity,
the obstruction discussed above vanishes: if $(\alpha,A)$ is not
homotopic to a constant then it can not be conjugate to a cocycle
homotopic to a constant, and in particular it can not be reducible
(conjugate to a constant matrix).

The main global result of this paper is the following: {\it a typical
one-dimensional
quasiperiodic cocycle which is not homotopic to the identity has a
positive Lyapunov exponent.}  However, we are able to show this result
without exclusion of frequencies:

\begin{thm} \label {globalresult}

Let $\alpha \in \R \setminus \Q$.
Let $Z^r_\alpha \subset C^r(\R/\Z,\SL(2,\R))$, $r=\omega,\infty$,
be the set of all $A$
non-homotopic to a constant such
that $L(\alpha,A)=0$.  Then $Z^r_\alpha$ can be written as a disjoint union
$Z^r_\alpha=O^r_\alpha \cup NO^r_\alpha$, where
\begin{enumerate}
\item $O^r_\alpha$ has positive codimension in
$C^r({\R}/{\Z},\SL(2,{\R})$\footnote{By this we mean
that, locally, $O^r_\alpha$ is contained in the zero set
of a $C^r$ real-valued function for which $0$ is a regular value.},
\item $NO_\alpha$ is closed and for every
$A \in C^r(\R/\Z,\SL(2,\R))$, the set of $\theta \in \R$ such that
$R_\theta A \in NO_\alpha$ has zero Lebesgue measure.
\end{enumerate}

\end{thm}

Notice that item (1) in the above description can not be removed:
if $A \in C^\omega(\R/\Z,\SO(2,\R))$ then $L(\alpha,R_\theta A)=0$
for every $\theta \in \R/\Z$.  This event is ``essentially'' (modulo
conjugation) what item (1) covers.

\begin{rem}

The conclusion of item (2) above still holds if instead of the family
$\theta \mapsto R_\theta A$ one
considers any $C^r$ perturbation of it.

\end{rem}

\bignote{more on this later. Also, what about commentary below}

\comm{
The two following theorems should clarify in
which sense we understand ``full measure'' in the space of cocycles.

Let $\RR^r_\alpha \subset C^r(\R/\Z,\SL(2,\R))$ be the set of $A$ such that
$(\alpha,A)$ is conjugate to a cocycle of rotations.

Our first result stablishes that $\RR^r$ has ``full measure'' in the
complement of non-uniformly hyperbolic cocycles.

\begin{thm}

Let $A \in C^r(\R/\Z,\SL(2,\R))$, $r=\omega,\infty$
be non-homotopic to the identity and
let $\alpha \in \R$.  Then for almost every $\theta \in \R$, either
$L(\alpha,A)>0$ or $A \in \RR^r_\alpha$.

\end{thm}

Notice that the second possibility in the Theorem happens for a positive
measure set of values of $\theta$ (indeed all $\theta$)
whenever $A$ is itself a cocycle of
rotations.  However, it turns out that $\RR^r_\alpha$ is quite small
(this differs radically from the case of cocycles homotopic to the
identity): it is a ``subvariety'' with positive codimension.

\begin{thm}

Let $\alpha \in \R$.  Then $\RR^r_\alpha$
is the zero set of a $C^r$ function (defined in a neighborhood of
$\RR^r_\alpha$) with non-vanishing second derivative.

\end{thm}

\begin{rem}

Those results also hold in finite differentiability (with loss of
derivatives).

\end{rem}
}

\subsection{Local theory: cocycles with a twist}

In order to prove the theorem above we develop a ``local'' theory of
cocycles non-homotopic to a constant.  This local theory does extend to
arbitrary dimension.

In the case of cocycles homotopic
to the identity, ``local'' stands for ``cocycles close to a constant''. 
In our case, there are no constant cocycles, so instead we discuss cocycles
which are close to the simplest examples of cocycles not homotopic to a
constant, which are of the form $(\alpha,A)$ with $A(x)=R_{\phi(x)}$, where
$\phi$ is a non-constant affine map of $\R^d/\Z^d$.
The key reason to study such cocycles (and for us to call this setting
``local'') comes from renormalization and will be clear later.

Let us note that, previously, a local theory had been
developed (in the one-dimensional case)
using a KAM scheme \cite {K}, thus it was perturbative
(it depends on arithmetic properties of $\alpha$).  Here we will follow a
different path, based on complexification ideas.  One of the key advantages
of this approach is that it is non-perturbative.  Indeed the local setting
is quite robust: it is $C^1$-open (notice that a satisfactory theory is
impossible in the $C^0$ topology by \cite {bo}).

Our local setting consists of the class
of ``monotonic cocycles'' (cocycles satisfying a twist condition).  More
precisely, a $C^1$ cocyle $(\alpha,A)$ is said to be monotonic if its
projectivized skew-product action $\R^d/\Z^d \times \mathbb{P}^1 \to \R/\Z
\times \R\P^1$ takes the
horizontal foliation $\{\R/\Z \times \{z\}\}_{z \in \mathbb{P}^1}$ to a
foliation transverse to the horizontal foliation.

Notice that whether
$(\alpha,A)$ is monotonic only depends on $A$, so we may define the set of
monotonic functions $M^r \subset
C^r(\R/\Z,\SL(2,\R))$.  Clearly monotonic
cocycles can not be homotopic to the identity.

\comm{
More generally, one can consider premonotonic cocycles $(\alpha,A)$, that
is, cocycles which admit a monotonic iterate.  Properties of monotonic
cocycles pass quite easily to premonotonic ones.
If $(\alpha,A) \in (\R \setminus \Q) \times C^1(\R/\Z,\SL(2,\R))$
is continuously conjugate to a cocycle of rotations then it is
premonotonic.  Whether $(\alpha,A)$ is premonotonic depends on both $\alpha$
and $A$.  We denote by $P^r_\alpha \subset C^r(\R/\Z,\SL(2,\R))$ the set of
$A$ such that $(\alpha,A)$ is premonotonic.
}

In the context of cocycles homotopic to the identity, one has to work a lot
to obtain regularity properties of the Lyapunov exponent, even just
continuity\footnote {See \cite {GS} (which also addresses H\"older
regularity, and \cite {BJ} for continuity in the global analytic setting.
There are also perturbative local results in the analytic and
smooth setting, see \cite {eliasson} and \cite {AK2}.}.
Perhaps the most surprising result of our analysis is the following:

\begin{thm}\label{theo:1.3}

Let $r=\omega,\infty$.  The Lyapunov exponent of $(\alpha,A)$ is a
$C^r$ function 
\bignote{It seems that we prove only G\^ateaux differentiability
and not Fr\'echet}
of $A \in M^r$.

\end{thm}

This result is quite easy to prove (at least in the analytic case),
once we have chosen the right framework.  Notice that in the case of
cocycles homotopic to a constant, one only expects smoothness of the
Lyapunov exponent in the uniformly hyperbolic regime, while such cocycles
simply do not exist in our context.\footnote{The Lyapunov exponent may not
be better than $1/2$-H\"older,
even for the Almost Mathieu Operator with Diophantine frequencies
(which falls in the local setting of Eliasson \cite {eliasson}).
Indeed, in this setting, the Lyapunov exponent vanishes in the spectrum but
near the endpoints of the gaps one has square root behavior.}

It is not much more difficult to vary the frequency.

\begin{thm}\label{theo:1.4}

The Lyapunov exponent of $(\alpha,A)$
is a $C^\infty$ function of $(\alpha,A) \in \R^d \times M^\infty$.

\end{thm}

We have the following result which characterizes monotonic
cocycles with a zero Lyapunov exponent:

\begin{thm}\label{theo:1.5}

Let $(\alpha,A) \in \R^d \times M^r$, $r=\omega,\infty$.  If
$L(\alpha,A)=0$ then $L(\alpha,A)$ is $C^r$ conjugate to a cocycle of
rotations.

\end{thm}

This theorem can be seen as a rigidity result: a priori, a zero Lyapunov
exponent is related to measurable conjugation to some standard models, such
as cocycles of rotations \cite {thieullen}.

As an example of a result
going in a quite different direction, we get:

\begin{thm}\label{theo:1.6}

Let $\alpha \in \R^d$ be irrational and let $A \in M^\infty$.  Then
the projective skew-product action of $(\alpha,A)$ is minimal.

\end{thm}

We believe many other results on monotonic cocycles
are accessible by the techniques we use here. 

\subsubsection{Premonotonic cocycles}

The class of monotonic cocycles is not dynamically natural: it is not
invariant by real-analytic conjugacy.  Thus it is natural to define the
class of premonotonic cocycles $P^r \subset \R^d \times C^r(\R^d/\Z^d,
\SL(2,\R))$ as the set of cocycles $(\alpha,A)$ which are $C^1$
conjugate to a monotonic cocycle.  We let
$P^r=\cup_{\alpha \in \R^d} \{\alpha\} \times P^r_\alpha$
(whether $(\alpha,A)$ is premonotonic depends both on
$\alpha$ and on $A$).

Premonotonic cocycles form a $C^1$ open set, and all properties of monotonic
cocycles that we discussed transfer automatically to this larger setting.
It is not difficult to see that a cocycle $(\alpha,A) \in
\R^d \times C^1(\R^d/\Z^d,\SO(2,\R))$ is always premonotonic provided that
$\alpha$ is
irrational and $A$ is not homotopic to a constant.

In the statement of Theorem \ref{globalresult}
we make reference to a partition
$Z_\alpha=O_\alpha \cup NO_\alpha$ of cocycles with zero Lyapunov
exponent.  We may give a precise definition of those
sets now: $O_\alpha=Z_\alpha \cap P^\omega_\alpha$ is the set of
premonotonic $A$ such that $L(\alpha,A)=0$.

It is easy to see that the first derivative of $A \mapsto
L(\alpha,A)$ vanishes in $O_\alpha$, and it is possible to show
(using the rigidity theorem) that its second derivative does not vanish.
It immediately follows:

\begin{cor}\label{cor:1.7}

The set $O_\alpha$ has positive codimension in $P_\alpha$.

\end{cor}

This gives item (1) of Theorem \ref{globalresult}.

It could be hoped that all cocycles with irrational frequencies
(non-homotopic to the identity) are premonotonic.  We will show however that
this is not the case: there are many (positive measure set in an open set of
parametrized families) cocycles which are not premonotonic.  The examples we
will construct build on results of Young \cite {Y1}.

At this point, it is not clear if   premonotonic cocycles form
the largest class of cocycles that behave dynamically as monotonic cocycles.
For instance, it is not clear if the class of premonotonic cocycles is fully
invariant under iteration or renormalization.  However, this class is
more than enough for several purposes.  Some (possibly larger) classes of
cocycles which are invariant under several reasonable operations (including
renormalization) and share the properties of monotonic cocycles
will be discussed in this work.  None of them includes the examples
we mentioned above.

\subsection{Reduction from global to local}

\comm{
Let us remark that there is a difference between the global setting and the
local setting:

\begin{thm}

There exists an open subset of $\R \times C^\omega(\R/\Z,\SL(2,\R))$ of
cocycles which are not homotopic to the identity and not premonotonic.

\end{thm}

The examples of cocycles in the previous theorem have positive Lyapunov
exponent, and we do not know if there exist examples with zero Lyapunov
exponent and irrational frequency.  In any case, this result shows that one
should work more in order to understand the global case.
}

We now return to the one-dimensional case.
In order to prove Theorem \ref {globalresult}, we will proceed by reduction
to the local theory via a renormalization scheme.  Such a scheme was shown
in \cite {AK} to cover typical cocycles with zero Lyapunov exponent.
The argument is slightly more complicated here
because we want to take care of all frequencies.  Thus, instead
of considering limits of renormalization we are led to use cancelation
arguments to show convergence (of an appropriate sequence of
renormalizations) to ``standard models''.
(Cancellation schemes were developed in \cite
{K}, though we will present an alternative argument closer in spirit to
\cite {AK}.)  In particular, the renormalization scheme
leads us to monotonic cocycles.  The following consequence implies item (2)
of Theorem \ref {globalresult}.

\begin{thm}\label{theo:1.8}

Let $(\alpha,A) \in (\R \setminus \Q) \times C^r(\R/\Z,\SL(2,\R))$,
$r=\omega,\infty$, be
non-homotopic to the identity.  For almost every
$\theta \in \R/\Z$, either $L(\alpha,R_\theta A)>0$ or $(\alpha,A)$ is
$C^r$-conjugate to a cocycle of rotations (and automatically premonotonic).

\end{thm}

All results discussed here have analogues in finite differentiability, which
involve loss of derivatives.  Although we did not want to consider those
distractions in the introduction, we will state and prove the results also
in the case of finite differentiability, and we will give an estimate on the
loss of derivatives.  Many of the results are much stronger than we state in
the introduction, for instance, in Theorem \ref {globalresult}
one is not forced to consider families of the type $\theta \mapsto
R_\theta A$, for instance, any perturbation of such a family will do.

\subsection{Proof of Theorem \ref{globalresult}}
Item (2) of Theorem \ref{globalresult} follows directly from  Theorem \ref{theo:1.8}, while Item (1) is the content of Corollary \ref{cor:1.7}.

\subsection{Proof of Theorem \ref{theo:1.8}}
Let $A$ be $C^r$, $r=\infty,\omega$ non homotopic to the identity and $\alpha$ be irrational.  For Lebesgue almost every $\theta$ either $L(R_\theta A)>0$ or $(\alpha,R_\theta A)$ is $L^2$-conjugated to a cocycle with values in $SO(2,\R)$. This is a classical  fact from Kotani's theory which  in our paper is a consequence of Lemma \ref{key computa} and formula (\ref{2.26}) (in the proof of Theorem \ref{simple1}) or directly from Theorem \ref{L^2}. Consequently, one can apply the renormalization procedure to $(\alpha,R_\theta A)$ for $\theta$ in this set of full measure and, applying  Corollary \ref{cor:E.2} of the Appendix,  this provides us with a $C^r$ monotonic renormalization $(\tilde\alpha,A_{monot})$ of $(\alpha,R_\theta A)$ ($\tilde\alpha$ being the image of $\alpha$ under some iterate of the Gauss map).  Of course, we still have $L(\tilde\alpha,A_{monot})=0$. From Theoem \ref{theo:1.5} we get that $(\tilde\alpha, A_{monot}))$ is $C^r$-conjugated to a cocycle of rotations. But renormalization (and its inverse procedure)  clearly preserves  the property of being $SO(2,\R)$-valued. The cocycle $(\alpha,R_\theta A)$ is then conjugated to a cocycle of rotations and hence premonotonic. 

\subsection{Proof of Corollary \ref{cor:1.7}}
This follows from the analyticity of $A\mapsto L(A)$ on the set of $C^0$ non homotopic to the identity $A\in C^0(\R/\Z,\SL(2,\R))$, the comment after Lemma \ref{lemma:4.1} and the fact that every $C^r$ cocycle $(\alpha,A)$ in $MZ^r_\alpha$ is $C^r$ conjugated to a cocycle of rotations (Theorem \ref{theo:1.5}): if $\alpha$ is diophantine one can reduce to the case where $A=R_{nx}$ (since $Z(\cdot)\mapsto B(\cdot+\alpha)Z(\cdot)B(\cdot)^{-1}$ is clearly analytic) and  the computation after   Lemma \ref{lemma:4.1} shows that $D^2L(\alpha,A)\ne 0$; since $L(\alpha,A)$ is smooth with respect to $\alpha$ the result is true in the general case).

\subsection{Structure of the paper}

The key feature of monotonic cocycles is that they are amenable to
complexification techniques.  The typical complexification procedure, which
is quite developed for Schr\"odinger cocycles, is to perturb the cocycle by
a complex parameter (the energy in the Schr\"odinger case).
As far as we know, such complexification techniques were
restricted so far to very special kinds of perturbation, where
the dependence on the complex parameter is a holomorphic function with
very specific global properties.  In this work, we show that the
construction is much more robust, and the sole feature that enables it is
monotonicity of the perturbation parameter.  Our first step is to
generalize such parameter techniques to the case of general analytic
dependence.  Then, using a technique introduced for dynamical systems
by Lyubich \cite {Ly1} (in the context of unimodal maps), we extend those
results to the smooth setting. 
After obtaining several parameter style results, we obtain the
properties of monotonic cocycles by observing that the phase variable can be
considered as a parameter.

By these techniques, we will be able to prove  Theorem  \ref{theo:1.3}, 
Theorem \ref{theo:1.4} and Theorem \ref{theo:1.5} in section \ref{sec:3} using the results of section \ref{sec:2}.

\comm{
\section{$L^2$ estimates}

In this section we shall describe the complexification theory of Kotani. 
The objects and ideas introduced here will be the basis of our local theory.

\subsection{Lifts}

Let $A:\R/\Z \times \SL(2,\R)$ be continuous.  Let $\phi_A:\R/\Z
\times S^1 \to S^1$ be given by $\phi_A(x,v)=\frac {A(x) \cdot v}
{\|A(x) \cdot v\|}$.

A lift of $A$ is a continuous function
$\zeta:\R \times \R \to \R$ such that $\zeta(x,y)=y'$ if and only if
$\phi_A(x,(\cos 2 \pi y,\sin 2 \pi y))=(\cos 2 \pi y',\sin 2 \pi y')$.
Any lift satisfies $\zeta(x,y+1)=\zeta(x,y)+1$.
Two lifts differ by a constant integer.  Any lift satisfies
$\zeta(x+1,y)-\zeta(x,y)=d$ where $d$ is the degree of $A$.  Let
$\psi_\zeta:\R/\Z \times \R/\Z \to \R$ be given by
$\psi_\zeta(x,y)=\zeta(x,y)-y-dx$.

\subsection{Fibered rotation number}

Given $\alpha \in \R$ and a continuous $A:\R/\Z \to
\SL(2,\R)$, define $\phi_{\alpha,A}:\R \times S^1 \to \R \times S^1$ by
$\phi_{\alpha,A}(x,v)=(x+\alpha,\phi_A(x,v))$.  Given a lift $\zeta$ of $A$,
define $\zeta_{\alpha,\zeta}:\R \times \R/\Z \to \R$ by
$\zeta_{\alpha,\zeta}(x,y)=(x+\alpha,\zeta(x,y))$.  Define
\be
\psi^{(n)}_{\alpha,\zeta}(x,y)=\sum_{k=0}^{n-1} \psi_\zeta \circ
(\zeta_{\alpha,\zeta})^k(x,y).
\ee

Notice that
$|\psi^{(n)}_{\alpha,\zeta}(x,y)-\psi^{(n)}_{\alpha,\zeta}(x,y')|<1$.  In
particular, the {\it fibered rotation number}
\be
\rho_{\alpha,\zeta}=\int_{\R/\Z \times S^1} \psi_\zeta d\mu
\ee
is independent of the choice of a probability measure $\mu$ invariant by
$\phi_{\alpha,A}$.

A different choice of a lift $\zeta$ changes $\rho$ by a constant integer. 
This allows us to define $\rho_{\alpha,A} \in \R/\Z$.

\subsection{Families}

Let $A_\theta \in C^0(\R/\Z \to \SL(2,\R)$, $\theta \in \R$
be continuous in $\theta$.

A coherent lift of $A_\theta$ is a family of lifts $\zeta_\theta$ which
depends continuously on $\theta$.

A coherent determination of the fibered rotation number is a function.

measurable family of cocycles
such that for every $x \in \R/\Z$, $\theta \to A_\theta(x)$ is continuous.

A coherent lift of $A_\theta$ is a family of lifts $\zeta_\theta$ such that
$\theta \mapsto \zeta_\theta(x,y)$ is continuous for every $(x,y) \in \R
\times \R$.

A coherent lift is determined by the family $A_\theta$ and the value of
$\zeta_{\theta_0}$ for some $\theta_0$.

\section{Renormalization}

\section{Monotonic cocycles}

Let $A \in C^0(\R/\Z,\SL(2,\R))$.  Let $\phi_A:\R/\Z \times S^1 \to \R/\Z
\times S^1$ be defined by $\phi_A(x,v)=\phi_A(x,\frac {A \cdot v} {\|A \cdot
v\|}$.  There exists a continuous function
$\psi_A:\R \times \R \to \R \times \R$ such that
$\Pi \circ \psi_A=\phi_A \circ \Pi$ where $\Pi:\R \times \R \to \R/\Z \times
S^1$ is given by $\Pi(x,y)=(x,(\cos 2 \pi y,\sin 2 \pi y))$.  Two choices of
$\psi_A$ differ by a constant integer in the second coordinate.

We say that $A \in C^0(\R/\Z,\SL(2,\R))$ is $\epsilon$-monotonic if
for every $y \in \R$, and for every
$x \neq x' \in \R$ we have
\be
\left |\frac {\psi(x,y)-\psi(x',y)} {x-x'} \right | \geq \epsilon.
\ee
We say that $A$ is monotonic if it is $\epsilon$-monotonic for some
$\epsilon>0$.

Let $\MM^r(\epsilon) \subset C^r(\R/\Z,\SL(2,\R))$ be the set of
$\epsilon$-monotonic cocycles.  Let $\MM^r=\cup_{\epsilon>0}$
\MM^r(\epsilon)$ and $\MM(\epsilon)=\MM^0(\epsilon)$, $\MM=\MM^0$.

If $n \neq 0$ then $A(\theta)=R_{n \theta}$ is monotonic.

The set of monotonic cocycles is open in the Lipschitz topology

If $A \in \MM$ then $A$ is not homotopic to the identity.

If $A \in \MM^r(\epsilon)$ then $R_\theta A \in \MM^r(\epsilon)$.

\subsection{Continuous monotonic cocycles}

\begin{thm}

Let $A \in \MM^0$ and let $\alpha \in \R$.
Any coherent
determination $\rho$ of the fibered rotation number of $\theta \mapsto
(\alpha,R_\theta A)$
is a Lipschitz function of $\theta$.  More precisely, if $A$ is
$\epsilon$-monotonic and has degree $n$ then $\rho$ is
$|n|\epsilon^{-1}$-Lipschitz.

\end{thm}

\begin{pf}

It suffices to prove the second claim.  To fix ideas, let us assume $n>0$.

Let $\zeta_\theta$ be a coherent lift of $\theta \mapsto A_\theta$ and let
$\rho$ be the associated coherent determination of the
fibered rotation number.
Since $\theta \mapsto \rho(\theta)$ is a non-decreasing function of
$\theta$, it is enough to show that for every
$\theta_0>0$ and for every $\theta>0$ we have
$\rho(\theta_0+\theta) \leq \rho(\theta_0)+n \epsilon^{-1} \theta$.

Let $B_\theta(x)=R_{\theta_0} A(x+\epsilon^{-1}(\theta-\theta_0))$,
Let $\zeta'_\theta$ be a coherent lift of $\theta \mapsto B_\theta$
satisfying $\zeta'_{\theta_0}=\zeta_0$ and let $\rho'$ be the associated
coherent determinations of the fibered rotation number.  Then
$\rho'(\theta_0+\theta)=\rho(\theta_0)+n \epsilon^{-1} \theta$.  From the
definition of $\epsilon$-monotonicity, we see that for $\theta>\theta_0$ we
have $\zeta_\theta \leq \zeta'_\theta$.  But this implies that
$\rho(\theta) \leq \rho'(\theta)$ for $\theta>\theta_0$ as required.
\end{pf}

\begin{cor} \label {l2monotonic}

Let $A \in C^0(\R/\Z,\SL(2,\R))$ be monotonic and let $\alpha \in \R$. 
Then either $L(\alpha,A)>0$ or $(\alpha,A)$ is $L^2$-conjugated to a
cocycle of rotations.

\end{cor}

\begin{cor}

Let $A \in C^0(\R/\Z,\SL(2,\R))$ be monotonic and let $\alpha \in \R$.  Then
$\theta \to L(\alpha,R_\theta A)$ is continuous.  More precisely, it
has a derivative in BMO.

\end{cor}

\begin{cor}

For every $\epsilon>0$, the Lyapunov exponent is a continuous function of
$(\alpha,A) \in \R \times \MM(\epsilon)$.  In particular, $L(\alpha,A)$ is
continuous in $\R \times \MM^{\rm Lip}$.

\end{cor}

\subsection{Analytic monotonic cocycles}

Let $\MM^\omega(\epsilon,\delta)$ be the set of $A$ such that $A$ has an
extension to $\C/\Z \cap \{|\Im(z)|<\delta\}$ and if $0<|t|<\delta$ then
                                                       
\begin{enumerate}

\item If $t$ and $n$ have the same sign then
$Q A(x+it) Q^{-1}:\overline \C \to \overline \C$ takes $\D$ into
$\D_{1-\epsilon}$,

\item If $t$ and $n$ have different sign then
$Q A(x+it)^{-1} Q^{-1}:\overline \C \to \overline \C$ takes $\D$ into
$\D_{1-\epsilon}$,

\end{enumerate}

\begin{lemma} \label {delta}

Let $A \in \MM^\omega(\epsilon)$ have degree $n$.  There exists
$\delta>0$ (only depending on $\epsilon$ and on the bounds on the
holomorphic extension of $A$) such that $A$ belongs to
$\MM^\omega(\epsilon,\delta)$.

\end{lemma}

\begin{pf}

This follows from the Cauchy-Riemann equations.
\end{pf}

\begin{thm}

Let $(\alpha,A) \in \R \times
\MM^\omega(\epsilon,\delta)$ and let $0<|t|<\delta$.
Then $L(\alpha,A)=L(\alpha,x \mapsto A(x+it))-|t n|$ where $n$ is the degree
of $A$.

\end{thm}

\begin{pf}

To fix ideas, let $n>0$.  Then there exists a holomorphic function
$\tau:\C/\Z \cap \{0<\Im(z)<\delta\} \to \D$
such that $Q A(z) Q^{-1} \cdot \tau(z)=\tau(z+\alpha)$.  In particular, the
function $L(t) \equiv L(\alpha,x \mapsto A(x+i t))$ is an affine function of
$0<t<\delta$.  Since $L(-t)=L(t)$ and using subharmonicity, we see that
$L(t)-L(0)$ is a multiple of $|t|$.  Let $\rho(\sigma)$ be such that
$\sigma+i t \mapsto L(t)+2 \pi \rho(\sigma)$ is a holomorphic function of
$\sigma+i t \in \C/\Z \cap \{0<\Im(z)<\delta\}$.
We must show that $\rho(x)-n x$ is a constant.  As we saw before,
$\rho(\sigma)$ may be taken as a coherent determination of the fibered
rotation number of $\theta \mapsto
(\alpha,x \mapsto A(\sigma+\theta))$, so the result follows.
\end{pf}

\begin{cor}

The Lyapunov exponent is a smooth function of $(\alpha,A) \in \R \times
\MM^\omega$, and is analytic in the second coordinate.

\end{cor}

\begin{lemma} \label {conid}

Let $A \in \MM^{\rm Lip}$.  Then there is no measurable $B:\R/\Z \to
\SL(2,\R)$ such that $A(x)=B(x+\alpha) B(x)^{-1}$.

\end{lemma}

\begin{pf}

If $\alpha \in \Q$ then this is an easy exercise using monotonicity.  Let us
assume that $\alpha \in \R \setminus \Q$.

If there exists $B:\R/\Z \to \SL(2,\R)$ such that $A(x)=B(x+\alpha)
B(x)^{-1}$ then $L(\alpha,A)=0$.  By Corollary \ref {l2monotonic}, there
exists $C \in L^2(\R/\Z,\SL(2,\R))$ such that $C(x+\alpha) A(x) C(x)^{-1}
\in \SO(2,\R)$.  It follows that $D=CB$ satisfies $D(x+\alpha)D(x)^{-1} \in
\SO(2,\R)$, so $\|D(x)\|$ is almost everywhere constant.  Thus $B \in
L^2(\R/\Z,\SL(2,\R))$.  Using \cite {AK}, we see that for almost every $x_*
\in \R$, if $n_k \to \infty$ satisfies $\|n_k \alpha\|<\frac
{1} {n_k}$, then the sequence $B(x_*)^{-1} A_n(x_*+\frac {x} {n_k})
B(x_*) \in C^{\rm {Lip}}(\R,\SL(2,\R))$ converges to $\id$.  But this
implies that all sufficiently deep renormalizations of $A$ have degree $0$,
so $A$ has degree $0$ contradicting the monotonicity of $A$.
\end{pf}

\begin{rem}

We do not know if the degree is an invariant under measurable conjugacy
under reasonable smoothness assumptions.

\end{rem}

\begin{lemma} \label {sec}

Let $(\alpha,A) \in (\R \setminus \Q) \times C^0(\R/\Z,\SL(2,\R))$.
If there exist more than one function $\tau:\R/\Z \to \overline \D$
satisfying the invariance
relation $Q A(x) Q^{-1} \cdot \tau(z)=\tau(x+\alpha)$ then either

\begin{enumerate}

\item There are exactly two such functions, and both are real,

\item We have $A(x)=B(x+\alpha) B(x)^{-1}$ for some measurable
function $B:\R/\Z \to \SL(2,\R)$.

\end{enumerate}

\end{lemma}

\begin{lemma} \label {vab}

Let $\tau:\C/\Z \cap \{0<|\Im(z)|<\delta\} \to \D$ be such that
$\lim_{t \to 0} \tau(\sigma+i t)$ exists for almost every $\sigma$.
Then $\tau$ admits a holomorphic extension
$\tau:\C/\Z \cap \{0 \leq |\Im(z)|<\delta\} \to \D$.

\end{lemma}

\begin{pf}

Consider the Fourier series of $\tau$.
\end{pf}

\begin{thm}

Let $(\alpha,A) \in \R \times \MM^\omega$.  If $L(\alpha,A)=0$ then
$(\alpha,A)$ is $C^\omega$-conjugated to a cocycle of rotations.

\end{thm}

\begin{pf}

By Lemma \ref {delta}, we have $A \in
\MM^\omega(\epsilon,\delta)$ for some $\epsilon>0$, $\delta>0$.
Then we have a holomorphic function
$\tau:\{\sigma+i t,\, 0<|t|<\delta\} \to \D$
satisfying the invariance relation
$Q A(z) Q^{-1} \tau(z)=\tau(x+\alpha)$.
Define measurable functions
$\tau_\pm:\R \to \overline \D$ by
$\tau_\pm(\sigma)=\lim_{t \to 0\pm} \tau(\sigma+i t)$ for almost every
$\sigma \in \R$.  Notice that both $\tau_\pm$ satisfy the invariance
relation.

Let us first deal with the more interesting case $\alpha \in \R \setminus
\Q$.  If $\tau_+ \neq \tau_-$ then, using
Lemma~\ref {conid}, we must be in the first case of Lemma~\ref {sec}.
But this is impossible because of Corollary~\ref {l2monotonic}.
Thus $\tau_+=\tau_-$ and so it follows from removability, see Lemma~\ref
{vob}, $\tau$ admits a holomorphic extension to
$\tau:\{\sigma+it,\, 0 \leq |t| <\delta\} \to \D$
satisfying the invariance relation.  The result follows.

Let now $\alpha=\frac {p} {q} \in \Q$.
Then $L(\alpha,A)=0$ means that $A_q(x)$ is not hyperbolic, $x \in \R/\Z$.
It follows from monotonicity that $A_q(x)$ is elliptic for all $x$ but
finitely many.  We conclude easily that $\tau_\pm$ agree in $\R$ and the
result follows.
\end{pf}

\subsection{Smooth monotonic cocycles}

In the previous discussion we have thoroughly used the holomorphic extension
of $A$.  In order to extend those results to the smooth case, we will
consider asymptotically conformal extensions of smooth cocycles.  This idea
was first applied to dynamical systems by Lyubich \cite {teichfib}, in the
context of unimodal maps.

The main results of this section are the following:

\begin{thm}

Let $(\alpha,A) \in \R \times \MM^\infty$.  If $L(\alpha,A)=0$ then
$(\alpha,A)$ is $C^\infty$-conjugated to a cocycle of rotations.

\end{thm}

\begin{thm}

The Lyapunov exponent is a smooth function of $(\alpha,A) \in \R \times
\MM^\infty$.

\end{thm}
}

\comm{
\section{An estimate on the second derivative of the Lyapunov exponent for a
cocycle of rotations}

Since the Lyapunov exponent $L$
takes non-negative values, we must have $DL=0$ whenever $L=0$.
Here we are going to show that if $L=0$ then $D^2 L \neq 0$.  This implies
that $\{L=0\}$ is a subvariety of positive codimension.

\begin{lemma}

Let $B \in C^0(\R/\Z,\SL(2,\R))$.  Let $A_\theta(x)=R_{n x}
B(x-\theta)$.  Then
\be
\int_0^1 L(\alpha,A_\theta(x)) d\theta=\int^1_0 \ln \left
(\frac {\|B(x)\|+\|B(x)\|^{-1}} {2} \right ) dx.
\ee

\end{lemma}

\begin{pf}

Let $C_\theta(x)=R_{n \theta} R_{n x} B(x)$.
Notice that $A_\theta(x+\theta)=C_\theta(x)$.  In
particular, $L(\alpha,A_\theta)=L(\alpha,C_\theta)$.  The result follows by
\cite {AB}.
\end{pf}

Let $s \in C^0(\R/\Z,\Sl(2,\R))$, that is,
\be
s(x)=\begin{pmatrix}a(x)&b(x)+c(x)\\b(x)-c(x)&-a(x)\end{pmatrix},
\ee
where $a,b,c:\R/\Z \to \R$ are continuous functions.
Let $A_{\theta,t}=R_{n x} e^{t s(x-\theta)}$.  Then the previous lemma
implies that
\be
\lim_{t \to 0} \frac {1} {2 t^2} \int_0^1 L(\alpha,A_{\theta,t})
d\theta=\int_0^1 a^2(x)+b^2(x) dx.
\ee
In particular, the limit above is zero if and only if $a$ and $b$ vanish
identically, that is, if and only if $s$ takes values in $\so(2,\R)$.

\section{Codimension}

In this section we are going to show that $\{L=0\}$ has codimension $4n$.

\begin{thm}

Let $(\alpha,A) \in C^r(\R/\Z,\SL(2,\R))$ be conjugated to a cocycle of
rotations.  Then $D^L$ has rank $4n$.

\end{thm}

\begin{pf}

The result is obvious if $\alpha \in \Q$.

It is enough to consider the case where $A(x)=R_{n x}$.  If $\alpha$ is
close to $0$ then $D^2 L$ has rank $4n$.  The result follows by
renormalization.$
\end{pf}
}
}

\section{Monotonicity in parameter space}\label{sec:2}

In this section, $f:X \to X$ is a fixed homeomorphism of a compact metric
space $X$, preserving a fixed probability measure $\mu$, assumed to have
full support.  Given $A \in C^0(X,\SL(2,\C))$, we use the dynamics
$f:X \to X$ to define
the iterated matrix products $A_n(x)$, $n \in \Z$, by
\be \label {matrixproducts}
A_n(x)=A(f^{n-1}(x)) \cdot A_{n-1}(x),
\quad A_{-n}(x)=A_n(f^{-n}(x))^{-1}, \quad n \geq -1, \quad A_0(x)=\id.
\ee
The Lyapunov exponent is defined by
\be
L(A)=\lim_{n \to \infty} \frac {1} {n} \int_X \ln \|A_n(x)\| d\mu(x).
\ee

As discussed in the introduction, our aim here is to develop a ``smooth''
version of Kotani Theory, centered around the concept of monotonic
dependence with respect to the parameter.
Let $I \subset \R$ be an interval.  We say that a continuous function
$f:I \to \R$ is $\epsilon$-monotonic if for every $x \neq x'$ we have
\be
\frac {|f(x')-f(x)|} {|x'-x|} \geq \epsilon.
\ee
This definition naturally extends to functions defined on (or taking values
on) $\R/\Z$ (by considering lifts) and on the unit circle $\Sr^1 \subset
\R^2 \equiv \C$ (by considering the identification with $\R/\Z$ given by
$x \mapsto e^{2 \pi i x})$.  Naturally, we may distinguish between two types
of $\epsilon$-monotonicity, increasing or decreasing.

We say that a continuous one-parameter family of matrices
$A_\theta(\cdot) \in \SL(2,\R)$ is $\epsilon$-monotonic if,
for every $w \in \R^2 \equiv \C$, the function $\theta \mapsto
\frac {A_\theta \cdot w} {\|A_\theta \cdot w\|}$ is
$\epsilon$-monotonic.

We will be interested in one-parameter families of $\SL(2,\R)$ cocycles
displaying monotonicity with respect to the parameter variable.
Thus, a continuous one-parameter family $A_\theta \in C^0(X,\SL(2,\R))$ is
said to be $\epsilon$-monotonic increasing (respectively, decreasing)
if for every $x \in X$,
the family $\theta \mapsto A_\theta(x)$ is $\epsilon$-monotonic increasing
(respectively, decreasing).

For several results, we will need to assume
further regularity with respect to the parameter.  Let us say that the
family $A_\theta$ is $C^r$ in $\theta \in J$ if $\theta \mapsto
A_\theta(x)$ belongs to some fixed compact subset of $C^r(J,\SL(2,\R))$
(compact open topology) for each $x$.

Among the results we will obtain in this section, we highlight the following
ones:

\begin{thm} \label {rigidity}

Let $A_\theta \in C^0(X,\SL(2,\R))$ be monotonic and
$C^{r+1+\epsilon}$, $1 \leq r<\infty$, $C^\infty$, or $C^\omega$
in $\theta$.  If $L(A_\theta)=0$
for every $\theta$ in some open interval $J$
then there exists $B_\theta \in C^0(X,\SL(2,\R))$, $\theta \in J$
depending $C^r$, $C^\infty$ or $C^\omega$
on $\theta$ and conjugating $A_\theta$ to a cocycle of rotations.

\end{thm}

\begin{thm} \label {dependence}

Let $A_{\theta,s} \in C^0(\R/\Z,\SL(2,\R))$, $\theta \in \R/\Z$, $s$ a
one-dimensional real parameter, be monotonic in $\theta$ and
$C^{2r+1+\epsilon}$, $1 \leq r<\infty$, $C^\infty$ or $C^\omega$,
in $(\theta,s)$.
Then
\be
s \mapsto \int_{\R/\Z} L(A_{\theta,s}) d\theta
\ee
is $C^r$, $C^\infty$, or $C^\omega$.

\end{thm}

This theorem implies for instance that $A \mapsto
\int_{\R/\Z} L(R_\theta A) d\theta$
is an analytic function of $A \in C^0(\R/\Z,\SL(2,\R))$.  Indeed, it can be
shown (see \cite {AB}) that
\be
\int_{\R/\Z} L(R_\theta A) d\theta=\int_X \ln \frac
{\|A(x)\|+\|A(x)\|^{-1}} {2} d\mu(x).
\ee
This generalization beyond families with specific form such as $R_\theta A$
will be crucial when we start to mix phase and parameter in the analysis of
quasiperiodic cocycles displaying monotonicity with respect to some {\it
phase} variable.

We will also obtain several other results whose statements depend on the
concept of ``variation of the fibered rotation number'', which we will first
need to introduce.



\subsection{Complexification}

\def\mA{\mathring {{A}}}
\def\mB{\mathring {{B}}}

Much of the information we will get from matrices in $\SL(2,\C)$ will come
from their action on the Riemann Sphere $\overline \C$ through M\"obius
transformations:
\be
\bm a&b\\c&d\em \cdot z=\frac {az+b} {cz+d},
\ee

For $A \in \SL(2,\C)$, let
$\mA=Q A Q^{-1}$ where
\be
Q=\frac {-1} {1+i} \begin{pmatrix}1&-i\\1&i\end{pmatrix}.\label{2.4}
\ee
The map $A\mapsto \mA$ maps bijectively $SL(2,\R)$ to $SU(1,1)$,
the real Lie group of matrices
$\bm u&\bar v\\  v& \bar u\\\em$, $u,v\in\C$ such that
$|u|^2-|v|^2=1$.\footnote{
If one identifies the complex
one-dimensional projective space $\mathbb{CP}^1$ with $\overline
\C$, by associating to the line through $(0,0) \neq (x,y) \in \C^2$
the complex number $\frac {x-iy} {x+iy}$, the action of $A \in \SL(2,\C)$ on
$\mathbb{CP}^1$ is given precisely by $z \mapsto \mA \cdot z$.  In this
identification, the real one-dimensional projective space $\mathbb{RP}^1$
corresponds to the unit circle $\partial \D$, and $\SL(2,\R)$ matrices
preserve the unit disk $\D$.}





\subsection{Variation of the fibered rotation number}\label{sec:variationofthe}

In the analysis of Schr\"odinger cocycles, the notion of fibered rotation
number plays a fundamental role (often through the analysis of its close
cousin, the integrated density of states).  In our setting, it turns out
that it is not always possible to define ``the''
fibered rotation number of a cocycle.  However, we will be able to define
the notion of {\it variation of the fibered rotation number} along a path.

Define $\Upsilon$ as the space of $\SL(2,\C)$ matrices $A$ such that $\mA \cdot
\D \subset \D$ (or equivalently $A \cdot \H\subset\H$).

Let $A \in \Upsilon$.  Define $\tau_A:\overline \D \to \C \setminus \{0\}$ by
\be
\mA  \begin{pmatrix}z\\1\end{pmatrix}=\tau_A(z)
\cdot \begin{pmatrix}\mA \cdot z\\1\end{pmatrix}.
\ee
Since $\overline \D$ is simply connected there exists a map
$\hat\tau_A:\overline \D\to \C$ such that $e^{2\pi i\hat\tau_A(z)}=\tau_A(z)$;
any other lift is obtained by the addition of an integer.
If we denote by $\hat\Upsilon$ the universal cover of $\Upsilon$ considered
as a topological semi-group with unity $\hat \id$, there
exists a unique continuous map
$\hat \tau:\hat\Upsilon\times \overline \D\to \C$ such that $\hat \tau(\hat \id,z)=0$
and $e^{2\pi i\hat\tau(\hat A,z)}=\tau_A(z)$.  This map satisfies
\be
\hat\tau(\hat A_2 \hat A_1,z)=\hat\tau(\hat A_2,\mA_1\cdot z)+
\hat\tau(\hat A_1,z).\label{taulift}
\ee
We note that for any $\hat A\in\hat\Upsilon$ and any $z,z'\in \overline \D$
$$\Im\hat\tau(\hat A,z)=-\frac{1}{2\pi}|\ln\tau_A(z)|$$
$$|\Re \hat\tau(\hat A,z)-\Re\hat\tau(\hat A,z')|<1/2.$$
The equality is trivial, while the inequality follows from the fact that for
any $A$, $\tau_A(\overline \D)$ is contained in an open half plane
(it is enough to observe that if $\mA=\begin{pmatrix} u& \overline v\\v & \overline u \end{pmatrix}$ then
$\tau_A(z)=v z+\overline u$, so $\tau_A(\overline \D)$ does not intersect the line
through $i \overline{u}$, since $|u|^2-|v|^2=1$).

\comm{
It is enough to prove that for any $A\in\Upsilon$, $z,z'\in\D$ one has $|\Re\hat\tau_A(z)-\Re\hat\tau_A(z')|<1$. We observe that if $A\in\Upsilon$ there exist   $U\in SU(1,1)$ and $D_w=\bm w&0\\0&w^{-1}\em$ where $w\in\bar \D\setminus\{0\}$ such that $\mA=UD_wU^{-1}$ or changing notations $A=A_2A_1$, $A_1,A_2\in\Upsilon$ with $\mA_2=UD_z$, $\mA_1=U^{-1}$. A simple computation shows that 
$\tau_{A_2}(z)=wvz+w^{-1}\bar u$ if $U=\bm u&\bar v\\v&\bar u\em $ ($|u|^2-|v|^2=1$).
We claim that $|\Re\hat\tau_{A_2}(z)-\Re\hat\tau_{A_2}(z')|<\frac{1}{2}$. Indeed, if this were not the case there would exist two points $z_*,z'_*$ in the line segment $[z,z']$ ($\{z_*,z'_*\}\ne \{z,z'\}$) such that $\Re\tau_{A_2}(z_*)-\Re\tau_{A_2}(z'_*)=1/2$. In other words there would exist $t>0$ such that 
$$wvz_*+w^{-1}\bar u=-t(wvz_*'+w^{-1}\bar u).$$
This would imply $|v||w||z_*+tz_*'|=|u||w^{-1}|(1+t)$. But since $|u|>0$, $|w|\leq 1$ and  $|v|\leq |u|$ we would have $1+t\leq |z_*+tz_*'|$ which is in contradiction with the fact that $z_*,z'_*$ are in $\D$. To conclude, we observe that the same computation holds with $w=1$ and $U^{-1}$ in place of $U$ and we make use of (\ref{taulift}).
}

Now if $\gamma:[0,1]\to\Upsilon$ is continuous, and
$\hat \gamma:[0,1]\to \hat\Upsilon$ is a continuous lift, we define
$\delta_\gamma\hat\tau(z_0,z_1)=\hat\tau(\hat\gamma(1),z_1)-\hat\tau(\hat\gamma(0),z_0)$;
notice that this is independent of the choice of the lift $\hat\gamma$.


Let us note a nice composition rule: given $\gamma$ and $\gamma'$,
let $\gamma' \gamma(t)=\gamma'(t) \gamma(t)$.
Then
\be \label {formulacomp}
\delta_{\gamma' \gamma} \hat \tau(z_0,z_1)=
\delta_{\gamma'} \hat \tau(\mathring{\gamma}(0)z_0,\mathring{\gamma}(1)z_1)+
\delta_\gamma \hat \tau(z_0,z_1).
\ee


Consider now a continuous path $\gamma \in C^0([0,1],C^0(X,\Upsilon))=
C^0([0,1] \times X,\Upsilon)$.
Define $\delta_\gamma \xi:X \times \overline \D \times \overline \D \to \C$ by
$\delta_\gamma \xi(x,z_0,z_1)=
\delta_{\gamma_x} \hat \tau(z_0,z_1)$, where $\gamma_x(t)=\gamma(t,x)$.

Using the dynamics $f:X \to X$, we define
paths $\gamma_n \in C^0([0,1],C^0(X,\Upsilon))$ by putting
$\gamma_n(t,x)=\gamma(t,f^{n-1}(x)) \cdots \gamma(t,x)$.
Define $\delta_\gamma \xi_n \in C^0(X \times \overline \D \times \overline \D,\C)$ by
$\delta_\gamma \xi_n=\frac {1} {n} \delta_{\gamma_n} \xi$.
We have
an expression for $\delta_\gamma \xi_n$ as a Birkhoff average (for the
dynamical system $(x,z_0,z_1) \mapsto (f(x),\mathring {\gamma}
(0,x) \cdot z_0,\mathring {\gamma}(1,x) \cdot z_1)$ acting on
$X \times \overline \D \times \overline \D$):
\be
\delta_\gamma \xi_n(x,z_0,z_1)=\frac {1} {n} \sum_{k=0}^{n-1}
\delta_\gamma \xi(f^k(x),\mathring {\gamma}_k(0,x) \cdot z_0,
\mathring {\gamma}_k(1,x) \cdot z_1).
\ee
This is obtained by the composition formula (\ref {formulacomp}).

We claim that $\lim_{n \to \infty} \delta_\gamma \xi_n(x,z_0,z_1)$
exists for $\mu$-almost every $x$, and is independent of
$z_0,z_1 \in \D$. We will show this by proving convergence for the real and imaginary parts.

Since for every $z_0,z_1 \in \overline \D$ we have
$|\Re(\delta_\gamma\xi_n(x,z_0,z_1)-\delta_\gamma\xi_n(x,z_0',z_1'))|<\frac
{1} {n}$,
it follows  from Birkhoff Ergodic Theorem that
\be \label {fiberedrotation}
\lim_{n \to \infty} \Re \delta_\gamma \xi_n(x,z_0,z_1)
\ee
exists and is independent of
$z_0,z_1 \in \overline \D$ for $\mu$-almost every $x$.
We call the $\mu$-average of (\ref {fiberedrotation}) the
{\it variation of the fibered rotation number}
along $\gamma$ and we denote it by $\delta_\gamma\rho$.  It is obviously invariant by homotopy, and it is
a continuous function of $\gamma \in C^0([0,1],C^0(X,\Upsilon))$.

\begin{rem}

If $\mu$ is ergodic, (\ref {fiberedrotation}) is
$\mu$-almost everywhere constant.
If $f$ is uniquely ergodic, (\ref {fiberedrotation}) exists for every $x$
and is constant.  If the iterates of $f$ are uniformly equicontinuous
(for instance, if $f$ is a translation of the torus, but
perhaps non-ergodic),
(\ref {fiberedrotation}) exists for every $x$, and is a continuous function
of $x \in X$.

\end{rem}

The convergence of the imaginary part  of $ \delta_\gamma \xi_n(x,z_0,z_1)$ is somewhat more delicate,
and it is certainly less robust.  For
$A \in C^0(X,\SL(2,\C))$ the Lyapunov exponent is defined by 
\be
L(A,x)=\lim_{n \to \infty} \frac {1} {n}
\ln \left \|A_n(x)  \right \|
\ee
(which exists $\mu$-almost everywhere by subadditivity). When  $A \in C^0(X,\Upsilon)$,  a simple application of the Oseledets Theorem shows that for any
$z\in \D$ 
\be
L(A,x)=\lim_{n \to \infty} \frac {1} {n}
\ln \left \|\mA_n(x) \cdot\begin{pmatrix} z\\1\end{pmatrix} \right \|
\ee
since $\begin{pmatrix} z\\1 \end{pmatrix}$ can not belong to the {\it stable
direction}  (this is due to the fact that in that case $A_k(x) \in \Upsilon$
for every $k \geq 0$).

\begin{lemma}

For $\mu$-almost every $x \in X$ and for every $z_0,z_1 \in \D$,
\be
\lim_{n \to \infty} \Im \delta_\gamma \xi_n(x,z_0,z_1)=\frac {1} {2 \pi}
(L(\gamma(0),x)-L(\gamma(1),x)).
\ee

\end{lemma}

\begin{pf}

Notice that if $A \in C^0(X,\SL(2,\C))$
\be
 \ln \left \|\prod_{k=n-1}^0 \mA(f^k(x)) \cdot
\begin{pmatrix}z\\1\end{pmatrix} \right \|=\ln|\tau_{A_n(x)}(z)|+
\ln \left \|\bm\prod_{k=n-1}^0 \mA(f^k(x)) \cdot z\\1\em
\right \|,
\ee
and the second term in the right hand side is bounded.
\end{pf}

We let $\delta_\gamma \zeta$ be the $\mu$-average of $\lim \delta_\gamma
\xi_n(x,0,0)$.  As we have seen, its real part is $\delta_\gamma \rho$ and its
imaginary part is $\frac {1} {2 \pi} (L(\gamma(0))-L(\gamma(1))$, where
$L(A)$ is defined as in the introdution:
\be
L(A)=\lim_{n \to \infty} \frac {1} {n} \int \ln \|A_n\| d\mu.
\ee

\begin{rem} \label {concate}

Notice that $\delta_\gamma \zeta$ behaves well under concatenation: if
$\gamma$, $\gamma'$ and $\tilde \gamma$ are such that $\gamma(1)=\gamma'(0)$
and $\tilde \gamma$ is homotopic to the concatenation of $\gamma$ and
$\gamma'$ then
$\delta_{\tilde \gamma} \zeta=\delta_{\gamma'} \zeta+\delta_\gamma \zeta$.

\end{rem}

\subsubsection{Invariant section}\label{sec:invariantsection}

Assume that for every $x \in X$, we have $A(x) \cdot \overline \D
\subset \D$.
In this
case, by the Schwarz Lemma, $A(\cdot)$ uniformly contracts the Poincar\'e metric of the disk and  there exists $m \in
C^0(X,\D)$ satisfying
\be
m(f(x))=\mA(x) \cdot m(x),\label{2.12}
\ee
and we have for every $(x,z) \in X \times \D$,
\be
\lim_{n\to\infty} \mA_n(f^{-n}(x)) \cdot z=m(x).
\ee
Thus
\be
\delta_\gamma \zeta=\int \delta_\gamma
\xi(x,m_{\gamma(0)}(x),m_{\gamma(1)}(x)) d\mu(x),\label{eq:2.20}
\ee
where $m_{\gamma(t)}$ stands for the invariant section corresponding to
$A=\gamma(t)$.

Notice that there is another formula for the Lyapunov exponent in terms of
$m$.  Let $q(x)$ be the norm of the derivative of the holomorphic function
$z \mapsto \mA(x) \cdot z$ at $z=m(x)$, with respect to some
conformal Riemannian metric $\| \cdot \|_z$ on $\D$.
In this case
\be\label{2.17}
L=\frac {1} {2} \int -\ln q(x) d\mu(x).
\ee

The most convenient metric to consider is the Poincar\'e metric
$(1-|z|^2)^{-1}|dz|$, since it enables us to apply the Schwarz Lemma.
Notice that whenever 
\be \mA  \begin{pmatrix}z\\ 1\end{pmatrix}=\tau\begin{pmatrix}\tilde z\\ 1\end{pmatrix}
\ee
one has 
\be\label{2.16} \biggl|\frac{d\tilde z}{dz}\biggr|=\frac{1}{|\tau|^2},\qquad q=\biggl|\frac{d\tilde z}{dz}\biggr|\frac{1-|z|^2}{1-|\tilde
z|^2},
\ee
and (\ref {2.17}) can be also obtained as a consequence
of~(\ref{2.12}) (which implies that $L=\int \ln
|\tau_{A(x)}(m(x))|d\mu(x)$).

The Schwarz Lemma yields, for instance, the following estimate of
the Lyapunov exponent.  If $\mA(x) \cdot \D \subset \D_{e^{-\epsilon}}$
for every $x \in X$ then
the composition $m \mapsto \tilde m=\mA(x) \cdot m \mapsto e^\epsilon \tilde m$
sends $\D$ to itself
and hence weakly contracts Poincar\'e metric:
\be  \frac{|d (e^\epsilon \tilde m)|}{1-|e^\epsilon\tilde m|^2} \leq  \frac{|dm|}{1-|m|^2}
\ee
from which one gets 
\be \frac{|d\tilde m|}{1-|\tilde m|^2}\leq e^{-\epsilon}\frac{1-e^{2\epsilon}|\tilde m|^2}{1-|\tilde m|^2} \frac{|dm|}{1-|m|^2}.
\ee
Whence
\be\label{2.20}
q(x)^{-1} \geq e^{\epsilon}\frac {1-|m(f(x))|^2} {1-e^{2\epsilon}|m(f(x))|^2}\geq e^\epsilon
\ee
so that $L \geq \frac {\epsilon}{2}$.

\subsubsection{Fibered rotation function}

Let us now consider a continuous family $A_\theta \in
C^0(X,\Upsilon)$, where $\theta$ belongs to some connected
Banach manifold $M$. Since $\delta_\gamma\zeta$ depends only on the homotopy class of
$\gamma:[0,1]\to C^0(X,\Upsilon)$, and it behaves well under concatenation,
see Remark \ref {concate}, we can define a map $\zeta:\tilde M\to \C$
($\tilde M$ being the universal cover of $M$) such that for every
path $\tilde \gamma:[0,1] \to \tilde M$, we have $\delta_\gamma
\zeta=\zeta(\tilde \gamma(1))-\zeta(\tilde \gamma(0))$, where
$\gamma(t)=A_{\pi(\tilde \gamma(t))}$ and $\pi:\tilde M\to M$ is the
canonical projection.
Moreover, we can take $\zeta$ so that
for $\pi(\tilde \theta)=\theta$ we have $-2 \pi \Im
\zeta(\tilde \theta)=L(A_\theta)$.
We shall then  denote $\rho(\tilde \theta)=\Re \zeta(\tilde \theta)$
or with an abuse of notation
$\rho(\theta)$ or $\rho_{A_\theta}$.  Though $\rho$ is only defined up to a
real constant,
it makes sense to speak
about the derivative of $\rho$, and if $M=\R$ or $M=\R/\Z$, we may ask
whether $\rho$ is monotonic or not.




\subsection{Simple applications} \label {simple2}

We now turn to one-parameter continuous families
$\theta \mapsto A_\theta (\cdot)\in C^0(X,\SL(2,\R))$.
To keep definite and to avoid superfluous notations, we will consider
only the case when the parameter space is $\R/\Z$.
To fix ideas, we will always assume in the proofs below
that $A_\theta(\cdot)$ is monotonic decreasing.  Note that
due to our identification of $\P\R^2$ with $\partial \D$,
the projectivization map
$\Sr^1 \to \P\R^2$ is orientation reversing.  Thus if $\theta \mapsto
A_\theta(\cdot)$ is
monotonic decreasing then for each $x \in X$,
$z \in \partial \D$, $\theta \mapsto \mA_\theta(x) \cdot z$ is monotonic
increasing.


\comm{
For the results below, the key assumption on the family
$A_\theta(\cdot)$ will be monotonicity in
$\theta$.  To fix ideas, we will always assume in the proofs below
that $A_\theta(\cdot)$ is monotonic increasing.  One obvious consequence of
monotonicity (following directly from the definitions) is the following.

\begin{lemma}

Let $A_\theta(\cdot) \in C^0(X,\SL(2,\R))$,
be a one-parameter family
monotonic in $\theta$.  Then the fibered rotation number
is non-decreasing as a function of $\theta$.

\end{lemma}

Before dwelving into more complicated matters, let us illustrate
the relation of monotonicity and complexification with two simple
applications.
}

\comm{
Let us say that a family
$\theta \mapsto A_\theta (\cdot)\in C^0(X,\SL(2,\R))$
is $C^r$ in $\theta$ if there exists a compact subset $K
\subset C^r(X,\SL(2,\R))$ such that for every $x \in \R/\Z$,
$\theta \mapsto A_\theta(x)$ belongs to $K$.
}

Let us show how to deduce the analytic case of Theorem \ref {dependence}.
Let $\Omega_\delta=\{z \in \C/\Z,\, |\Im(z)|<\delta\}$,
$\Omega^\pm_\delta=\{z \in \C/\Z,\, 0<\pm \Im(z)<\delta\}$, and let $\tilde
\Omega_\delta$, $\tilde \Omega^\pm_\delta$, be their universal covers.

If $\theta\mapsto A_\theta$ is assumed to be monotonic and analytic,
the Cauchy Riemann equations imply that
there exists $\delta>0$ such that
the analytic extension of $A_\theta$ to $\theta
\in \Omega^+_\delta$ satisfies
\be\label{contr}
\mA_\theta \cdot \D \subset \D_{e^{-2\epsilon \Im(\theta)+
\kappa(\Im(\theta))}}, \quad 0<\Im \theta \leq \delta,
\ee
where $\kappa(t)<2 \epsilon t$ for $0<t \leq \delta$
and $\lim_{t \to 0} \kappa(t)/t=0$.  As we saw in section \ref{sec:invariantsection}, this implies that $L(A_\theta) \geq \epsilon
\Im \theta-\frac {1} {2}
\kappa(\Im \theta)$.

Indeed, for fixed $z\in \partial {\D}$, $x\in X$, the Cauchy-Riemann equations applied to the analytic function
$F_{x,z}:\Omega^+_\delta\to \D$ 
defined by $F_{x,z}(\sigma+i t)=\mA_{\sigma+it}(x)\cdot z$ shows that
$(\partial_\sigma F_{x,z}(\theta),\partial_t F_{x,z}(\theta))$ is a
directed orthogonal basis with vectors of equal length when $\theta\in\R/\Z$.
Since $\partial_\sigma F_{x,z}(\theta)$ is tangent to
the circle $\partial\D$, and of length at least $2 \epsilon$ (by
monotonicity), we see that
$\partial_t F_{x,z}(\theta)$ is radial (pointing to the origin  because
of the sign assumption) and of length at least
$2 \epsilon$ for any $\theta\in\R/\Z$. This proves the above
estimate (\ref{contr}). 

By the previous discussion, we can choose $\delta>0$ so that for $0<t
\leq \delta$ we have $\mA_{\sigma+it,s}(x) \cdot \overline \D \subset \D$
and we  are thus in a situation where the discussion of section \ref{sec:variationofthe} applies since $A_{\sigma+it,s}(\cdot)\in C^0(X,\Upsilon)$:  there exists a function $\zeta:\tilde\Omega_\delta^+\to\C$ such that on the closure of  $\tilde \Omega_\delta^+$  the map $\theta\mapsto\rho(\theta)=\Re\zeta(\theta)$ is continuous  and $\Im\zeta(\theta)=-\frac{1}{2\pi}L(A_\theta)$; also, since $\theta\mapsto A_\theta$ is holomorphic on  $\Omega_\delta^+$, the same is true for   $\zeta$. We will use the notation   $m^+(z,x)$  for the $\D$-valued invariant section of the cocycle $A_z$.

Now if $(\theta,s)\mapsto A_{\theta,s}$ is analytic and  if $s$ is in some fixed  neighborhood of  $s_0$, we can choose $\delta>0$ so that for $0<t
\leq \delta$ we have $\mA_{\sigma+it,s}(x) \cdot \overline \D \subset \D$.
Let then
\be \label {eq:2.28}
U(t,s)=\int_{\R/\Z} L(A_{\sigma+it,s}) d\sigma,
\ee
so that for $t \neq 0$ we have
\be \label {eq:2.29}
U(t,s)=\int_{X \times \R/\Z} \ln |\tau_{A_{\sigma+it,s}(x)}(m^+_s(\sigma+it,x))| d\mu(x)
d\sigma.
\ee

Then for $0<t<\delta$ the map $s \mapsto U(t,s)$ is analytic.
Moreover, the map $(\sigma+it,z) \mapsto
\tau_{A_{\sigma+it,s}(x)}(z)$ (defined on $\Omega^+_\delta \times \D$
is holomorphic and non zero and so
$\int_X \ln |\tau_{A_{\sigma+it,s}(x)}(m^+_s(\sigma+it,x))| d\mu(x)$ is harmonic w.r.t
$\sigma+i t\in \Omega_\delta^+$; its integral w.r.t
$\sigma$, given by $U(t,s)$,  is thus an affine
function of $0<t<\delta$ since it is harmonic in $\sigma+i t$
and does not depend on
$\sigma$.  

On the other hand, since $\sigma+it\mapsto A_{\sigma+it,s}(x)$
is holomorphic, the functions $\sigma+it\mapsto \ln\|(A_{\sigma+it,s})_n(x)\|$ are subharmonic and the same is true for $\sigma+it\mapsto L(A_{\sigma+it,s})$.
Notice that
$U(t,s)=U(-t,s)$ since $\sigma+it\mapsto L(A_{\sigma+it,s})$ is real symmetric, so by
subharmonicity, $U(t,s)$ is an affine function of
$|t|$ for $0 \leq |t|<\delta$
($t \mapsto U(t,s)$ is convex and thus  continuous in $t$).
Thus, for $0<t<\delta$ we have
\be
\int_{\R/\Z} L(A_{\theta,s}) d\theta=2 U(\frac {t} {2},s)-U(t,s),
\ee
is analytic on $s$, as desired.

\begin{rem}

With a little bit more work, one can get the formula
\be \label {form}
\int_{\R/\Z} L(A_{\theta,s}) d\theta=U(t,s)-2\pi t \deg,
\ee
for $0<t<\delta$, where $\deg$ is the variation of the fibered rotation number
as $\theta$ runs once around $\R/\Z$.\footnote {It is easy to see that
$\deg$ is just the $\mu$-average of $\deg(x)$, where $\deg(x)$ is the
topological degree of $\theta \mapsto A_{\theta,s}(x)$ as a map
$\R/\Z \to \SL(2,\R)$ (in particular, if $X$ is connected, $\deg$ is an
integer).}
Indeed, for fixed $s$, let
$\rho(\sigma+it,s):\tilde \Omega^+_\delta \to \R$ be a
continuous determination of $\rho_{A_{\sigma+it,s}}$, so that
$\deg=\rho(\sigma+it+1,s)-\rho(\sigma+it,s)$.  Then the function
\be
\int_\sigma^{\sigma+1} 
-\frac{i}{2\pi}L(A_{y+it,s})+\rho(y+it,s) dy
\ee
is holomorphic in $\sigma+it \in \tilde \Omega^+_\delta$
and its real part is
an affine function of $\sigma$ of slope $\deg$.  Thus the
function $U(t,s)$ defined above is an affine function of
$0<t<\delta$ (Cauchy-Riemann) 
with slope $2\pi \deg$, and (\ref{form}) follows.

\end{rem}

Let us now describe another application, the analogous of a basic
derivative bound in Kotani Theory.  To set it up, let us state
an obvious consequence of
monotonicity (following directly from the definitions).

\begin{lemma}

Let $A_\theta(\cdot) \in C^0(X,\SL(2,\R))$,
be a one-parameter family
monotonic decreasing in $\theta$.  Then the fibered rotation number
is non-increasing as a function of $\theta$.

\end{lemma}

\begin{thm} \label {simple1}

Let $A_\theta(\cdot)\in C^0(X,\SL(2,\R))$,
$\theta \in \R/\Z$, be analytic and monotonic decreasing in $\theta$.
For Lebesgue a.e  $\theta \in \R/\Z$, if $L(A_\theta)=0$ then
\be
-\frac {d} {d\theta} \rho(\theta) \geq \frac {\epsilon} {2 \pi}>0,
\ee
where $\epsilon$ is the monotonicity constant of
$\theta \mapsto A_\theta(\cdot)$.

\end{thm}

\begin{pf}

Using the analyticity in $\theta$ of $A_\theta$, we can conclude that
$\theta \mapsto \zeta(\theta)$ can be defined as a holomorphic function
$\tilde \Omega^+_\delta \to \H$.
We know that the real part of $\zeta(\theta)$ (the ``fibered
rotation number'')  is continuous up to the closure.  For $\sigma\in  
\R$,
\be
\Im(\zeta({{\sigma+it}}))=\Im(\zeta({\sigma+i0^+}))+\int_{0^+}^t\partial_s \Im(\zeta({\sigma+is}))ds,
\ee
and using the Cauchy-Riemann equations
\be
\Im(\zeta({{\sigma+it}})=\Im(\zeta({\sigma+i0^+}))+\int_{0^+}^t\partial_\sigma \Re\zeta({\sigma+is})ds,
\ee
Since the  map $\Re\zeta(\cdot)$ is harmonic on $\Omega_\delta^+$,
continuous on the closure of $ \tilde \Omega_\delta^+$
and its restriction to $\Im\theta=0$ is non-increasing, one can
say\footnote{A continuous harmonic function $f$ on the disk is the
Poisson integral of its restriction $\rho$ to the boundary of the disk.
If $\rho$ is of bounded variation, the tangential derivative
$\partial_\sigma f$ is the Poisson integral of the  measure $d\rho$.
Fatou theorem asserts that for a.e point on the boundary, the radial
limit of $\partial_\sigma f$ is $d\rho/d\sigma$. The situation on the
strip is easily reduced to the one on the disk.}
that for Lebesgue a.e $\sigma\in\R/\Z$
\be
\lim_{s\to 0}\partial_\sigma \Re\zeta({\sigma+is})=
\frac{d}{d\sigma}\rho(\sigma).
\ee

Since the Lyapunov exponent is upper semicontinuous (it is by subadditivity the infimum of the continuous functions $\theta\mapsto (1/n)\int_X\ln\|(A_\theta)_n(x)\|d\mu(x)$),
if we know additionally that $L(A_\sigma)=0$ it becomes continuous at the point $\sigma$ and  we have
\be
\label{2.26}
-\lim \frac {1} {2 \pi t} L(A_{\sigma+it})=\frac {d} {d \sigma} \rho(\sigma),
\ee
almost surely,
and the result follows (since $L(A_{\sigma+it}) \geq \epsilon t-\frac {1} {2}
\kappa(t)$).
\end{pf}

\subsection{General framework}

After the motivation above, we are ready to introduce a more general
framework for the complexification argument.

It may look like analyticity is crucial in
order to exploit the complexification approach.  This is not the case: in
the non-analytic case, we can still complexify the problem using
asymptotically holomorphic extensions (this idea is inspired from the work
of Lyubich on smooth unimodal maps \cite {Ly1}).

Let $\Delta_\delta$ be the space of all continuous families
$A_{\sigma+it}(\cdot) \in C^0(X,\SL(2,\C))$, $\sigma+it \in
\Omega_\delta$,
which are $C^1$ and real-symmetric in $\sigma+it$, satisfying
$A_{\sigma+it}
\in \inter \Upsilon$, $\sigma+it \in \Omega^+_\delta$,
\be \label {asymptoti1}
\op_z A_z=0, \quad {\rm if}\ \Im(z)=0,
\ee
and such that $\sigma \mapsto A_\sigma(\cdot)$ is monotonic in $\sigma$. 
Condition (\ref {asymptoti1}) is an asymptotic holomorphicity assumption,
some stronger forms of which we will later introduce.

Let us fix $A \in \Delta_\delta$.
Then we have functions $m^+(\sigma+it,x) \in \D$,
$\tau^+(\sigma+it,x) \in \C \setminus \{0\}$,
$\sigma+it \in \Omega^+_\delta$, $x \in X$, characterized by
\be \label {+}
\mA_{\sigma+it}(x) \cdot \bm m^+(\sigma+it,x)\\1\em=\tau^+(\sigma+it,x) \bm
m^+(\sigma+it,f(x))\\1\em.
\ee

Notice that $A_{\sigma+it}(x)^{-1} \in \inter \Upsilon$ for $\sigma+it \in
\Omega^-_\delta$.  Thus we have also functions $m^-(\sigma+it,x)$,
$\tau^-(\sigma+it,x)$, $\sigma+it \in \Omega^-_\delta$, $x \in X$,
characterized by
\be \label {-}
\mA_{\sigma+it}(x) \cdot \bm m^-(\sigma+it,x)\\1\em=\tau^-(\sigma+it,x) \bm
m^-(\sigma+it,f(x))\\1\em.
\ee

Since $A_{\sigma+it}$ is real-symmetric in $\sigma+it$, letting
\be
m^+(\sigma+it,x)=\frac {1} {\overline {m^+(\sigma-it,x)}}, \quad
\tau^+(\sigma+it,x)=\frac {1} {\overline {\tau^+(\sigma-it,x)}}, \quad
\sigma+it \in \Omega^-_\delta
\ee
and
\be
m^-(\sigma+it,x)=\frac {1} {\overline
{m^-(\sigma-it,x)}}, \quad \tau^-(\sigma+it,x)=\frac {1} {\overline
{\tau^-(\sigma-it,x)}}, \quad \sigma+it \in \Omega^+_\delta,
\ee
we have that
(\ref {+}) and (\ref {-}) are valid for $\sigma+it \in \Omega_\delta
\setminus \R/\Z$.


The following key computation generalizes estimates of
Kotani \cite {K} (see also \cite {S}) and Deift-Simon \cite {DeS}.

\begin{lemma} \label {key computa}

Let $A \in \Delta_\delta$ and  let $\sigma_0 \in \R/\Z$.  Then

\begin{enumerate}

\item If
\be
\liminf_{t \to 0+} \frac {L(\sigma_0+it)} {t}<\infty
\ee
then
\be
\liminf_{t \to 0+} \int_X
\frac {1} {1-|m^+(\sigma_0+it,x)|^2} d\mu(x)+
\int_X \frac {1}
{1-|m^-(\sigma_0-it,x)|^2} d\mu(x)<\infty.
\ee

\item If
\be
\limsup_{t \to 0+} \frac {L(\sigma_0+it)} {t}<\infty
\ee
then
\be
\limsup_{t \to 0+} \int_X
\frac {1} {1-|m^+(\sigma_0+it,x)|^2} d\mu(x)+\int_X \frac {1}
{1-|m^-(\sigma_0-it,x)|^2} d\mu(x)<\infty
\ee
and
\be
\liminf_{t \to 0+} \int_X
|m^+(\sigma_0+it,x)-m^-(\sigma_0-it,x)|^2 d\mu(x)=0.
\ee

\end{enumerate}

\end{lemma}

\begin{pf}

Let us assume that
\be \label {one}
(\pa_t \mA_{\sigma_0+it}(x)) \mA_{\sigma_0+it}(x)^{-1}=u(x) \bm 
1&0\\ 0&-1\em+C(\sigma_0+it,x),
\ee
where $u(x)<0$ and $\lim_{t \to 0} \sup_{x \in X} \|C(\sigma_0+it,x)\|=0$.
The general case can be reduced to this one by
conjugacy.  Indeed, if $B:X \to \SL(2,\R)$ is continuous, then it is
indifferent to prove the results for the original $A$ or for its conjugate
$B(f(x)) A_{\sigma_0+it}(x) B(x)^{-1}$, so it is enough to
select $B$ such that
\be
\mB(f(x)) (\pa_t \mA_{\sigma_0}(x))
\mA_{\sigma_0}(x)^{-1}
\mB(f(x))^{-1}
\ee
is a matrix of the form $\bm u(x) & 0 \\ 0 & -u(x) \em $,
$u(x)<0$, for each $x$.
Such a $B$ can be found due to the monotonicity hypothesis as we will now
show.

Since $\bar\partial_z A_z=0$ at $\Im z=0$, we have
$(\partial_t \mA) \mA^{-1}=i (\partial_\sigma \mA) \mA^{-1}$
at $\sigma+i t=\sigma_0$, and
since $\mA_{\sigma}$  is  in  $\SU(1,1)$  for any  $\sigma\in{\R/\Z}$,
we can write
\be
(\partial_t \mA) \mA^{-1}=
\bm a& -i\bar\nu\\ -i\nu& -a\em,\qquad a\in\R,\quad \nu\in\C
\ee
at $\sigma+i t=\sigma_0$.  Moreover,
since $\mA_{\sigma+it}\cdot\D\subset \D_{e^{-\epsilon t+o(t)}} $  and $\mA_{\sigma}\cdot \D=\D$, we have
\be
\biggl|\frac{(1+ta)m-i\overline \nu t+o(t)}{-i\nu m
t+(1-at)+o(t)}\biggr|<1-\epsilon t+o(t)
\ee
for any $m\in\partial \D$ and any small $t>0$;
this implies that $a-|\nu|\sin\phi\leq -\epsilon/2$ for any
$\phi$ and consequently $a<0$ and
$\det((\pa_t \mA)\mA^{-1})=|\nu|^2-a^2\leq -\epsilon^2/4$.
It is then clear that
$(\partial_t \mA) \mA^{-1}$ is conjugated by a matrix in ${\rm SU}(1,1)$
to a matrix of the form $\bm u&0\\ 0&-u\em$ with $u<0$.
(the sign of $u$ cannot be changed to be positive because the conjugacy
is in ${\rm SU}(1,1)$).


Let us denote for simplicity $m^+$ for $m^+(\sigma_0+it,x)$, $m^-$ for
$m^-(\sigma_0+it,x)$, $\tilde m^+$ for $m^+(\sigma_0+it,f(x))$,
$\tilde m^-$ for $m^-(\sigma_0+it,f(x))$,
$\tau^+$ for $\tau^+(\sigma_0+it,x)$,
$\tau^-$ for $\tau^-(\sigma_0+it,x)$, $A$ for $A_{\sigma_0+it}(x)$, $L$ for
$L(\sigma_0+it)$ and $u$ for $u(x)$.

We now
estimate the Lyapunov exponent using (\ref{2.17}),  (\ref{2.16}) and the fact that $A_{\sigma+it}=\bm e^{tu}&0\\0 &e^{-tu}\em A_{\sigma}+o(t)$ by  evaluating the
contraction coefficient $q$ in the
Poincar\'e metric of $\D$.  A straightforward computation yields
\be
q^{-1}=e^{-2 t u+o(t)} \frac {1-|\tilde m^+|^2} {1-e^{-4 t u} |\tilde
m^+|^2}=e^{2 t u+o(t)} \frac {e^{-4 t u} (1-|\tilde m^+|^2)}
{1-e^{-4 t u} |\tilde m^+|^2}.
\ee
Using that
for $r>0$ and $0 \leq s<e^{-r}$ we have
\be
\ln \left (\frac {e^r (1-s)} {1-e^r s} \right ) \geq \frac {r} {1-s},
\ee
we get
\be
\ln q^{-1} \geq 2 t u+o(t)-\frac {4 t u} {1-|\tilde m^+|^2}=-2tu\frac
{1+|\tilde m^+|^2} {1-|\tilde m^+|^2}+o(t).
\ee
Since $L=\frac {1} {2} \int_X \ln q^{-1} d\mu$, we get
\be
L \geq -t \int_X u \frac {1+|\tilde m^+|^2} {1-|\tilde m^+|^2}
d\mu+o(t).
\ee
An analogous argument yields
\be
L \geq -t \int_X u \frac {1+|\tilde m^-|^{-2}} {1-|\tilde m^-|^{-2}}
d\mu+o(t).
\ee
We conclude that
\be \label {four}
\liminf_{t \to 0+} \frac {L} {t}+
\frac {1} {2} \int_X u \left (\frac {1+|\tilde m^+|^2} {1-|\tilde m^+|^2}+
\frac {1+|\tilde m^-|^{-2}} {1-|\tilde m^-|^{-2}} \right ) d\mu \geq 0.
\ee
Since $u<0$,
this gives the first item and the first part of the second item.

Differentiating
\be
\mA^{-1} \cdot \begin{pmatrix}\tilde m^+\\1\end{pmatrix}=
\frac {1} {\tau^+} \begin{pmatrix} m^+\\1\end{pmatrix}
\ee
with respect to $t$, and applying $-\mA$ to both sides, we get
\be
(\pa_t \mA) \mA^{-1}
\cdot \begin{pmatrix} \tilde m^+\\1\end{pmatrix}-\pa_t \tilde m^+ 
\begin{pmatrix}1\\0\end{pmatrix}=\frac {\pa_t \tau^+} {(\tau^+)^2}
\mA \cdot \bm m^+\\1\em-\frac {1} {\tau^+}
\pa_t m^+ \mA \cdot \bm 1\\0\em.
\ee
Using that
\be
\bm 1\\0\em=\frac {1} {m^+-m^-} \left (\bm m^+\\1\em-\bm m^-\\1\em \right 
)=\frac {1} {\tilde m^+-\tilde m^-} \left (\bm \tilde m^+\\1 \em-\bm \tilde
m^-\\1 \em \right ),
\ee
we get
\begin{align} \label {threeagain}
(\pa_t \mA)\mA^{-1} \cdot \bm \tilde m^+\\1\em-&\frac {\pa_t \tilde m^+}
{\tilde m^+-\tilde m^-}
\left (\bm \tilde m^+\\1\em-\bm 
\tilde m^-\\1\em \right )\\
\nonumber
&=\frac {\pa_t \tau^+} {\tau^+} \bm \tilde m^+\\1 \em-\frac {\pa_t 
m^+} {m^+-m^-} \left (\bm \tilde m^+\\1\em-\frac {\tau^-}
{\tau^+} \bm \tilde m^-\\1\em \right ).
\end{align}

On the other hand, we can compute using (\ref {one})
\be \label {fouragain}
(\pa_t \mA)\mA^{-1} \cdot \bm \tilde m^+\\1\em=u \left ( \frac {\tilde
m^++\tilde m^-}
{\tilde m^+-\tilde m^{-}}
\bm \tilde m^+\\1 \em-\frac {2\tilde m^+} {\tilde m^+-\tilde m^-}
\bm \tilde m^-\\1\em \right )
+c^+\bm \tilde m^+\\1\em+c^- \bm \tilde m^-\\1 \em,
\ee
where $c^+ \equiv c^+(\sigma_0+it,x)$, $c^- \equiv c^-(\sigma+it,x)$ satisfy
\be \label {two}
|c^\pm(\sigma_0+it,x)| \leq K \|C(\sigma_0+i t,x)\| \left (\frac
{1} {1-|\tilde m^+|^2}+\frac {1} {1-|\tilde m^-|^{-2}} \right ).
\ee
for some constant $K>0$.

Putting together (\ref {threeagain}) and (\ref {fouragain}),
taking the coefficient of $\bm \tilde m^+\\1 \em$ and
integrating with respect to $\mu$ we get
\be
\int_X u \frac {\tilde m^++\tilde m^-} {\tilde m^+-\tilde m^-} d\mu+
\int_X c^+ d\mu=\int_X \frac {\pa_t \tau^+} {\tau^+} d\mu.
\ee
We can now consider the real part, which gives
\be
\int_X u \frac{|\tilde m^+|^2 |\tilde
m^-|^{-2}-1} {\left |\frac {\tilde m^+} {\tilde m^-}-1
\right |^2} d\mu+\int_X \Re c^+ d\mu=\pa_t L.
\ee
Using (\ref {two}) we conclude
\be \label {five}
\lim_{t\to 0+} -\int_X u \frac{|\tilde m^+|^2 |\tilde m^-|^{-2}-1}
{\left |\frac {\tilde m^+} {\tilde m^-}-1
\right |^2} d\mu+\pa_t L=0.
\ee

Write
\be
I=\frac {1} {2} \frac {1+|\tilde m^+|^2} {1-|\tilde m^+|^2}+
\frac {1} {2} \frac {1+|\tilde m^-|^{-2}} {1-|\tilde m^-|^{-2}}+
\frac{|\tilde m^+|^2 |\tilde m^-|^{-2}-1} {\left |\frac {\tilde m^+} {\tilde
m^-}-1 \right |^2}.
\ee
Using (\ref {four}) and (\ref {five}) we get
\be
\liminf_{t\to 0+} \biggl(\frac {L} {t}-\pa_t L+
\int_X u I d\mu \biggr)\geq 0.
\ee

Notice that
\be \label {seven}
I \geq \left |\tilde m^+-\frac {1} {\overline {\tilde m^-}}
\right |^2 \geq 0,
\ee
and we conclude
\be
\liminf_{t\to0+}\biggl( \frac {L} {t}-\pa_t L\biggr)
\geq \liminf_{t\to 0+} -\int_X u I d\mu \geq 0.
\ee
Since $\limsup_{t\to 0+} \frac {L} {t}<\infty$, we must have
\be
\liminf_{t\to 0+}\biggl( \frac {L} {t}-\pa_t L\biggr)=-\limsup_{t\to
0+}\biggl(t \pa_t \frac {L} {t}\biggr) \leq 0,
\ee
so $\liminf_{t\to 0+} -\int_X u I d\mu=0$, and since $-u$ is positive
and bounded away from $0$ we have $\liminf_{t\to 0+} \int_X I d\mu(x)=0$,
which gives the second part of the second item by (\ref {seven}).
\end{pf}

The following estimates will allow us to work in the higher order
asymptotically holomorphic setting:

\begin{lemma}\label{lemma:2.5}

Let $A_{z,s} \in \Delta_\delta$ be a one-parameter family.
Assume that $s \mapsto A_{z,s}(x)$ is $C^r$, $1 \leq r<\infty$ and
\be
\|\pa^k_s A_{z,s}(x)\|=O(1), \quad 0 \leq k \leq r.
\ee
Then
\be \label {bla21}
|\pa^k_s m^+_s(z,x)|=O(|\Im(z)|^{-2k+1}), \quad 1 \leq k \leq r.
\ee
Moreover, if additionally $s \mapsto \op_z A_{z,s}(x)$ is $C^{r-1}$ and
we have the estimate
\be
\|\pa^k_s \op_z A_{z,s}(x)\|=o(|\Im(z)|^{\eta-k-1}), \quad 0 \leq k \leq r-1,
\ee
for some $\eta \in \R$ then
\be \label {bla22}
|\pa^k_s \op_z m^+_s(z,x)|=o(|\Im(z)|^{\eta-2k-2}), \quad 0 \leq k \leq r-1.
\ee

\end{lemma}

\begin{rem} \label {existdiffe}

Implicit in the statement of Lemma \ref {lemma:2.5}
is the existence of the derivatives taken in the left
hand sides of (\ref {bla21}) and (\ref {bla22}).  This is just a consequence
of usual normal hyperbolicity theory \cite {HPS}, but naturally the
estimates depend on the strength of the uniform hyperbolicity, and hence
may degenerate as $\Im z \to 0$, so our main concern below is to
get the bounds (\ref {bla21}) and (\ref {bla22}).

\end{rem}

\begin{pf}

Let $F^s_z(x,w)=\mA_{z,s}(x) \cdot w$, $m^s_z(x)=m^+_s(z,x)$.  Our estimates
will come from the study of the hyperbolicity of $F$ with respect to the
variable $w$, as measured in the Poincar\'e metric.  The way we exploit this
hyperbolicity is contained in the following.

\begin{prop}\label{prop:2.6}

There exists $K>0$ such that if $(s,z,x)\mapsto u^s_z(x)$ is any
bounded
function then
\be
|u^s_z(x)| \leq K |\Im(z)|^{-1} \sup_{y \in X} |u^s_z(f(y))-(\pa_w
F^s_z)(y,m^s_z(y)) u^s_z(y)|.
\ee

\end{prop}

\begin{pf}

For $s$ and $z$ fixed, let $x$ satisfy
\be
\frac {|u^s_z(f(x))|}
{1-|m^s_z(f(x))|^2}=M=\sup_{y \in X} \frac {|u^s_z(y)|}
{1-|m^s_z(y)|^2}
\ee
(we assume that the suppremum is achieved to keep the argument transparent,
the general case is obtained by approximation).
Then for every $y$,
\be
|u^s_z(y)| \leq \frac {|u^s_z(y)|} {1-|m^s_z(y)|^2} \leq M,
\ee
so it is enough to estimate
\be \label {M1}
M \leq K |\Im(z)|^{-1} |u^s_z(f(x))-(\pa_w F^s_z)(x,m^s_z(x)) u^s_z(x)|.
\ee
We have
\begin{align}
u^s_z(f(&x))-(\pa_w F^s_z)(x,m^s_z(x))
u^s_z(x)=(1-|m^s_z(f(x))|^2)\\
\nonumber
&\cdot \left (\frac {u^s_z(f(x))} {1-|m^s_z(f(x))|^2}
-(\pa_w F^s_z)(m^s_z(x)) \frac {1-|m^s_z(x)|^2}
{1-|m^s_z(f(x))|^2} \frac {u^s_z(x)} {1-|m^s_z(x)|^2} \right ).
\end{align}
Noticing that by the Schwarz Lemma
\be
\left | (\pa_w F^s_z)(x,m^s_z(x)) \frac {1-|m^s_z(x)|^2}
{1-|m^s_z(f(x))|^2} \right |<1,
\ee
we get
\begin{align} \label {M2}
|u^s_z(f(x))-&(\pa_w F^s_z)(x,m^s_z(x)) u^s_z(x)| \geq
M (1-|m^s_z(f(x))|^2)\\
\nonumber
&\cdot
\left (1-|(\pa_w F^s_z)(x,m^s_z(x))| \frac {1-|m^s_z(x)|^2}   
{1-|m^s_z(f(x))|^2} \right ).
\end{align}
The Schwarz Lemma hyperbolicity bound (cf. (\ref{2.20}))
\be
|\pa_w F^s_z(x,m^s_z(x))| \frac {1-|m^s_z(x)|^2} {1-|m^s_z(f(x))|^2} \leq
e^{-\epsilon \Im(z)} \frac {1-e^{2 \epsilon \Im(z)} |m^s_z(f(x))|^2}
{1-|m^s_z(f(x))|^2},
\ee
for some constant $\epsilon>0$, gives
\begin{align}
(1-|m^s_z(f(x))|^2)
&\left (1-|(\pa_w F^s_z)(x,m^s_z(x))| \frac {1-|m^s_z(x)|^2}   
{1-|m^s_z(f(x))|^2} \right )\\
\nonumber
&\geq
1-e^{-\epsilon \Im z}+(e^{\epsilon \Im z}-1) |m^s_z(f(x))|^2 \geq
1-e^{-\epsilon\Im(z)},
\end{align}
which together with (\ref {M2}) implies (\ref {M1}).
\end{pf}

Let us complete the proof of lemma~\ref{lemma:2.5}.
Differentiating (taking $\pa^k_s$)
\be
m^s_z(f(x))=F^s_z(x,m^s_z(x))
\ee
(as we are allowed to do, see Remark \ref {existdiffe}),
we get
\begin{align} \label {pa^k_s}
\pa^k_s m^s_z(f(x))=&\pa_w F^s_z(x,m^s_z(x)) \cdot \pa^k_s
m^s_z(x)\\
\nonumber
&+\sum_{\ntop {l \geq 0, 1 \leq i_1 \leq ... \leq i_l<k,} {i_1+...+i_l=j
 \leq k}}
C \cdot
\pa^{k-j}_s \pa^l_w F^s_z(x,m^s_z(x)) \cdot \prod_{n=1}^l
\pa^{i_n}_s m^s_z(x),
\end{align}
where $C \equiv C(k,i_1,...,i_l)>0$ is a constant.  Observe that 
$F^s_z(x,w)$ is a homography in $w$ whose denominator is bounded from below by $(10\sup_{x\in X}\|A_z(x)\|)^{-1}$ provided $\delta$ is small enough. This comes from the fact that if $z\in\R/\Z$ then  $\mA_z(x)\in SU(1,1)$, say equal to $\bm u& \bar v\\v&\bar u\em$ with $|u|^2-|v|^2=1$ and therefore for any $w\in\D$, $|vw+\bar u|\geq |u|-|v|\geq 1/(2\max(|u|,|v|))\geq (1/4\|A_z(x)\|)$. If $\delta$ is small enough, the claim follows. Consequently, all the quantities of the form $\pa_w^l\pa_s^jF^s_z(x,w)$ are uniformly bounded.

Thus, if
\be
|\pa^j_s m^s_z(x)|=O(|\Im(z)|^{-(2j-1)}), \quad 1 \leq j \leq k-1,
\ee
we get (either  $j=k$ and $l\geq 2$ or $j<k$) 
\be
|(\pa^k_s m^s_z(f(x))-\pa_w F^s_z(x,m^s_z(x)) \cdot \pa^k_s
m^s_z(x)|=O(|\Im(z)|^{-2(k-1)}),
\ee
which implies by the previous proposition
\be
|(\pa^k_s m)(s,z,x)|=O(|\Im(z)|^{-(2k-1)}).
\ee
The first estimate then follows by induction.

Differentiating (taking an extra $\op_z$) (\ref {pa^k_s}), we get
\begin{align}
&\pa^k_s \op_z m^s_z(f(x))=\pa_w F^s_z(x,m^s_z(x)) \cdot
\pa^k_s \op_z m^s_z(x)\\
\nonumber
&+\sum_{\ntop {l \geq 0, 1 \leq i_1 \leq ... \leq i_l\leq k,} {i_1+...+i_l=j
\leq k}}
C \cdot
(\pa^{k-j}_s \op_z \pa^l_w F^s_z(x,m^s_z(x)) \cdot
\prod_{n=1}^l \pa^{i_n}_s m^s_z(x)\\
\nonumber
&+\sum_{\ntop {l \geq 0, 1 \leq i_1 \leq ... \leq i_l<k,} {i_0 \geq 0,
i_0+i_1+...+i_l=j \leq k}}
D \cdot
\pa^{k-j}_s \pa^l_w F^s_z(x,m^s_z(x)) \cdot \pa^{i_0}_s \op_z m^s_z(x)
\cdot \prod_{n=1}^l \pa^{i_n}_s m^s_z(x),
\end{align}
where $D \equiv D(k,i_0,...,i_l)>0$ is a constant.  Since, 
\be
| \op_z m^s_z(f(x))-\pa_w F^s_z(x,m^s_z(x)) \cdot 
\op m^s_z(x)|=|\op_z  F^s_z(x,m^s_z(x))|= o(|\Im(z)|^{\eta-1}),
\ee
we get 
\be|\op m^s_z(x)|=o(|\Im(z)|^{\eta-2}),\ee
Now assuming by induction
\be
|\pa^j_s \op_z m^s_z(x)|=o(|\Im(z)|^{\eta-2j-2}), \quad 0 \leq j \leq k-1,
\ee
we get
\be
|\pa^k_s \op_z m^s_z(f(x))-\pa_w F^s_z(x,m^s_z(x)) \cdot \pa^k_s
\op m^s_z(x)|=o(|\Im(z)|^{\eta-2k-1}),
\ee
which implies as before
\be
|\pa^k_s \op_z m^s_z(x)|=o(|\Im(z)|^{\eta-2k-2}).
\ee
The second estimate then follows by induction.
\end{pf}

\begin{rem}

The estimates above are still valid if the parameter space is allowed to be
multidimensional, or, more generally, a Banach manifold, but the notation is
more cumbersome.

\end{rem}

\begin{rem}

As a particular case of the previous estimates (zero-dimensional
parameter space), if $A \in \Delta^+_\delta$ satisfies
\be
\|\op_z A_z(x)\|=o(|\Im(z)|^{\eta-1})
\ee
then
\be
|\op_z m^+(z,x)|=o(|\Im(z)|^{\eta-2}).
\ee

\end{rem}

\subsection{Asymptotically holomorphic extensions} \label {asymptotically}

In order to apply the estimates obtained in the previous section to
one-parameter families of $\SL(2,\R)$ cocycles, we need to consider
appropriate asymptotic holomorphic extensions.

For $\eta \in [1,\infty)$, a $C^\eta$ function defined in some neighborhood
of $\R/\Z$ satisfying
\be
\frac {d^k} {dt^k} \op F(\sigma+it)=0, \quad \sigma+it \in \R,\,
k=0,...,[\eta-1]
\ee
is called $\eta$-asymptotically holomorphic.

Let $AH^\eta$ be the set of $\eta$-asymptotically holomorphic functions
defined on the whole $\C/\Z$.
it is easy to see that one can define (linear) sections $\Phi_\eta$
of the restriction operator $AH^\eta \to C^\eta(\R/\Z,\C)$.
For instance, one can let
\be
\Phi_\eta(f)(\sigma+it)=\int K(x) f(\sigma+tx) dx,
\ee
where $K:\R \to \C$ is a $C^\infty$ function with compact support satisfying
\be
\int x^k K(x) dx=i^k, \quad k=0,...,[\eta+1].
\ee

We can also define $\eta$-asymptotically $\SL(2,\C)$-valued
functions by requiring each coefficient to be
asymptotically holomorphic.
In order to obtain asymptotically holomorphic extensions of a matrix valued
function $A=\bm a&b\\c&d\em
\in C^\eta(\R/\Z,\SL(2,\R))$, it is enough to consider
\be
\Phi_\eta(A)=(\Phi_\eta(a) \Phi_\eta(d)-\Phi_\eta(b)\Phi_\eta(c))^{-1/2}
\bm \Phi_\eta(a)&\Phi_\eta(b)\\\Phi_\eta(c)&\Phi_\eta(d)\em,
\ee
which is a well defined function $\Omega_\delta \to \SL(2,\R)$, where
$\delta$ only depends on the $C^1$-norm of $A$.

\begin{lemma}

Let $A_\theta(\cdot) \in C^0(X,\SL(2,R))$, $\theta \in \R/\Z$,
be monotonic and $C^\eta$ in $\theta$.  Then there exists
$\delta>0$ and an extension $A \in \Delta_\delta$ which is
$\eta$-asymotitucally holomorphic in $\theta$.

\end{lemma}

\begin{pf}

Consider the $\eta$-asymptoticaly holomorphic extension given above defined
on some $\Omega_\delta$.  We
just need to check that $A(\sigma+it,x) \in \inter \Upsilon$ for $t>0$
sufficiently small.  In the holomorphic case, this was concluded in section
\ref {simple2} using the Cauchy-Riemann equation, but the argument just uses
it for real values, and thus is valid in the asymptotically holomorphic
setting.
\end{pf}

\comm{
In order to illustrate the asymptotically holomorphic technique
(we refer to  Appendix \ref{appendix:A} for more details) ,
we generalize Theorems \ref {simple1} and \ref {simple2} to the
smooth setting.
}

The following result will illustrate
the use of higher order asymptotically holomorphic extensions:

\begin{lemma} \label {0leqt}

If $A_z \in \Delta_\delta$ is $1+\epsilon$-asymptotically holomorphic,
then $L(A_{\sigma+it})$
is a continuous function of $0\leq t<\delta$ for almost every $\sigma \in
\R/\Z$.

\end{lemma}

The proof will involve the following decomposition
technique which will play a role in several other arguments.
Given a function $u:\Omega^+_\delta \to \C/\Z$
which is continuous with derivatives in $L^1_\mathrm {loc}$ and satisfies
\be
|\op u(z)|=O(|\Im(z)|^{\epsilon-1}),
\ee
for some $\epsilon>0$, let us write a canonical decomposition
$u=u_h+u_c$ where $u_h:\Omega^+_\delta \to \C/\Z$
is holomorphic and $u_c:\C/\Z \to \C$ is a real-symmetric
continuous function given by the Cauchy transform 
\be
u_c(z)=\lim_{t \to \infty} \frac {-1} {\pi} \int_{[-t,t] \times
[-\delta,\delta]} \frac {\phi(w)} {z-w} dx dy,
\ee
where $\phi(z)=\op u(z)$ if $0<\Im(z)<\delta$ and $\phi(z)=\overline
{\op u(\overline z)}$ for $0<-\Im(z)<\delta$.

Notice that if
\be
|\op u(z)|=O(|\Im(z)|^{k+\epsilon}),
\ee
then $u_c(z)$ is complex differentiable at each $z \in \R/\Z$, and
$u_c:\R/\Z \to \R$ is $C^{k+1}$.

\noindent{\it Proof of Lemma \ref {0leqt}.}
By Lemma \ref {lemma:2.5}, $\|\op_z m^+(z,x)\|=O(|\Im z|^{\epsilon-1})$ in
$\Omega_+$.  Thus we can write a decomposition
$m^+=m^+_h+m^+_c$ where $z \mapsto m^+_c(z,x)$ is continuous up to $\R/\Z$
(uniformly in $x$),
and $m^+_h$ is of course uniformly bounded for $z$ near $\R/\Z$.  Thus
$m^+_h$ admits non-tangential limits as $\Im z \to 0$ and hence $m^+$ also
does.  By Fubini's Theorem, for almost every $\sigma$,
the right hand side of the expression
\be
L(A_{\sigma+it})=\int_X \ln
|\tau_{A_{\sigma+it}(x)}(m^+(\sigma+it,x))| d\mu(x),
\ee
originally defined for $0<t<\delta$, makes sense up to $t=0$ and defines a
continuous function of $0 \leq t<\delta$.  Thus we just
need to show that we can identify the right hand side (when it makes sense)
with $L(A_\theta)$ for $t=0$.  Indeed, if $\theta$ is such that
the non-tangential limits $m^+(\theta,x)$ exist, then they provide an
invariant section for the cocycle $A_\theta$.
Assuming ergodicity, by the
Oseledets Theorem the left hand side must be $\pm L(A_\theta)$, and (since
it is non-negative (by continuity),
it must be $L(A_\theta)$.  The general case reduces to this
one by ergodic decomposition.
\qed

\comm{
Given a function $u:\Omega^+_\delta \to \C/\Z$
which is continuous with derivatives in $L^1_\mathrm {loc}$ and satisfies
\be
|\op u(z)|=O(|\Im(z)|^{\epsilon-1}),
\ee
for some $\epsilon>0$, let us write a canonical decomposition
$u=u_h+u_c$ where $u_h:\Omega^+_\delta \to \C/\Z$
is holomorphic and $u_c:\C/\Z \to \C$ is a real-symmetric
continuous function given by the Cauchy transform 
\be
u_c(z)=\lim_{t \to \infty} \frac {-1} {\pi} \int_{[-t,t] \times
[-\delta,\delta]} \frac {\phi(w)} {z-w} dw \wedge d\overline w,
\ee
where $\phi(z)=\op u(z)$ if $0<\Im(z)<\delta$ and $\phi(z)=\op u(\overline
z)$ for $0<-\Im(z)<\delta$.

Notice that if
\be
|\op u(z)|=O(|\Im(z)|^{k+\epsilon}),
\ee
then $u_c(z)$ is complex differentiable at each $z \in \R/\Z$, and
$u_c:\R/\Z \to \R$ is $C^{k+1}$.
}

\subsection{Derivative bound}

Another simple application of
the asymptotically holomorphic technique is a
generalization of Theorem \ref {simple1}.

\comm{
Given a function $u:\Omega^+_\delta \to \C/\Z$
which is continuous with derivatives in $L^1_\mathrm {loc}$ and satisfies
\be
|\op u(z)|=O(|\Im(z)|^{\epsilon-1}),
\ee
for some $\epsilon>0$, let us write a canonical decomposition
$u=u_h+u_c$ where $u_h:\Omega^+_\delta \to \C/\Z$
is holomorphic and $u_c:\C/\Z \to \C$ is a real-symmetric
continuous function given by the Cauchy transform 
\be
u_c(z)=\lim_{t \to \infty} \frac {-1} {\pi} \int_{[-t,t] \times
[-\delta,\delta]} \frac {\phi(w)} {z-w} dw \wedge d\overline w,
\ee
where $\phi(z)=\op u(z)$ if $0<\Im(z)<\delta$ and $\phi(z)=\op u(\overline
z)$ for $0<-\Im(z)<\delta$.

Notice that if
\be
|\op u(z)|=O(|\Im(z)|^{k+\epsilon}),
\ee
then $u_c(z)$ is complex differentiable at each $z \in \R/\Z$, and
$u_c:\R/\Z \to \R$ is $C^{k+1}$.
}

\comm{
If $A_\theta \in C^0(X,\SL(2,\R))$ is $C^\eta$ in $\theta$, $\eta \geq 1$,
we can consider a
$\eta$-asymptotically holomorphic extension $A_z$.
Notice that estimate (\ref{contr}) obtained in section \ref {simple2} in the
analytic case
is still
true since the asymptotically holomorphic extension is complex
differentiable on $\R/\Z$.  In particular
$A_z \in \Delta_\delta$ for $\delta$ sufficiently small.

\begin{lemma}

If $A_z \in \Delta_\delta$ is $\eta$-asymptotically holomorphic with
$\eta>1$, then $\int_{\R/\Z} L(A_{\sigma+it}) d\sigma$
is a continuous function of $0\leq t<\delta$.

\end{lemma}

\begin{pf}

By Lemma \ref {lemma:2.5}, $\|\op_z m^+(z,x)\|=O(|\Im z|^{\epsilon-1})$ in
$\Omega_+$, where $\epsilon=\eta-1$.  Thus we can write a decomposition
$m^+=m^+_h+m^+_c$ where $z \mapsto m^+_c(z,x)$ is continuous up to $\R/\Z$
(uniformly in $x$,
and $m^+_h$ is of course uniformly bounded for $z$ near $\R/\Z$.  Thus
$m^+_h$ admits non-tangential limits as $\Im z \to 0$ and hence $m^+$ also
does.  It follows that the right hand side of the expression
\be
\int_{\R/\Z} L(A_{\sigma+it}) d\sigma=\int_{X \times \R/\Z} \ln
|\tau_{A_{\sigma+it}(x)}(m^+(\sigma+it,x))| d\sigma d\mu(x),
\ee
originally defined for $0<t<\delta$, makes sense up to $t=0$ and defines a
continuous function of $0 \leq t<\delta$.  Thus we just
need to identify the right hand side with
$\int_{\R/\Z} L(A_\theta) d\theta$ when $t=0$.

By Fubini's Theorem, for almost every
$\theta$, the non-tangential limits $m^+(\theta,x)$ exist, and provide an
invariant section for the cocycle $A_\theta$.  We claim that in this case
$\int_X \ln
|\tau_{A_\theta(x)}(m^+(\theta,x))| d\mu(x)=L(A_\theta)$.
Notice that the left hand side is non-negative.  Assuming ergodicity, by the
Oseledets Theorem the left hand side must be $\pm L(A_\theta)$, and by
non-negativity it must be $L(A_\theta)$.  The general case reduces to this
one by ergodic decomposition.
\end{pf}
}

\comm{
then $\int_X \ln
|\tau_{A_\theta(x)}(m^+(\theta,x))| d\mu(x)$ must be zero, or

Since we shall need this later, we observe that this allows us to identify for
Lebesgue a.e $\sigma\in\T$, the limit
$\lim_{t\to 0^+}\Im\zeta(\sigma+it)$ with $L(\sigma)$\footnote{In case the dependence
of $A$ on $z$ is holomorphic, we have already seen this to be true and if the dependence
is only continuous a simple case where this is also true is when $L(\sigma)=0$.}.
Indeed, we know that the $m$-functions $m^+(\sigma+it,x)$ are, uniformly in $x$,
asymptotically holomorphic: $ |\op_z m^+(\sigma+it,x)|=O(|\Im(z)|^{\eta-2})$.  We can
hence write $m^+(z,x)$ as a sum $m^+_h(z,x)+m^+_s(z,x)$ where $m^+_h(z,x)$ is holomorphic,
and if $\eta>1$, $m^+_s(z,x)$ is continuous and real symmetric. Observe that since
$m^+_h(z,x)$ is bounded, we can apply Fatou's theorem and deduce that for Lebesgue
a.e $\sigma\in\T$, the limit $\lim_{t\to 0^+}m^+(\sigma+it,x)$ exists.  Now we know that
for $\mu$-a.e. $x$ and all $t\geq 0$
\be L(\sigma+it)=\int_X|\ln\tau_{A_{\sigma+it}(x)} (m^+(\sigma+it,x))|d\mu(x)
\ee
We can hence apply Lebesgue dominated convergence Theorem to conclude. 
Let us mention a consequence of this fact: the map $t\mapsto\int_\T L(\sigma+it)d\sigma$ is continuous for all $t\geq 0$.
}

\begin{thm} \label {derivative bound}

Let $A_\theta \in C^0(X,\SL(2,\R))$,
$\theta \in \R/\Z$, be $C^{2+\epsilon}$ and monotonic decreasing
in $\theta$.
For almost every $\theta \in \R/\Z$, if $L(A_\theta)=0$ then
\be
-\frac {d} {d\theta} \rho(\theta) \geq \frac {\epsilon_0} {2 \pi}>0,
\ee
where $\epsilon_0$ is the monotonicity constant of $\theta \mapsto A_\theta$.

\end{thm}

\begin{pf}

For $\delta>0$ small, let us denote by $A_z \in \Delta_\delta$,
some fixed $2+\epsilon$-asymptotically holomorphic extension of
$A_\theta$, thus in particular
\be
|\op_z A_z(x)|=O(|\Im(z)|^{1+\epsilon}).
\ee
Notice that estimate (\ref{contr}) obtained in section \ref {simple2} in the
analytic case
is still true since the asymptotically holomorphic extension is complex
differentiable on $\R/\Z$.
Let us show that our hypothesis imply that for almost
every $\sigma \in \R$,
\be \label {L}
\partial_\sigma \rho(\sigma)=\lim_{t \to 0} \frac
{L(\sigma+it)-\lim_{t \to 0} L(\sigma+it)} {t},
\ee
since the result then follows as in Theorem \ref{simple1}.

We have
\be
|\op_z m^+(z,x)|=O(|\Im(z)|^\epsilon),
\ee
which implies from equation (\ref{eq:2.20})
\be
|\op_z \zeta(z)|=O(|\Im(z)|^\epsilon)\label{2.90}
\ee
as well, when $0<\Im z<\delta$.
Like in Theorem \ref{simple1} we get using the fact that $\zeta$ asymptotically
satisfies Cauchy-Riemann equations
\be
\Im(\zeta({{\sigma+it}}))=\Im(\zeta({\sigma+i0^+}))+\int_{0^+}^t\partial_\sigma \Re\zeta({\sigma+is})ds+o(|t|^\eta).
\ee
Notice that by Lemma \ref {0leqt}, $\Im \zeta(\sigma)=\Im
\zeta(\sigma+i0^+)$ for almost every $\sigma \in \R/\Z$.
From equation (\ref{2.90}), decomposing $\zeta=\zeta_h+\zeta_c$,
$\zeta_c(z)$ is complex differentiable at $z \in \R/\Z$ and
$\zeta_c:\R/\Z \to \R$ is $C^1$.  Since $\Im \zeta>0$ and
$\sigma \mapsto \rho(\sigma)=\lim_{t \to 0+} \Re \zeta(\sigma+it)$ is monotonic,
this is enough to conclude that (\ref {L}) holds for almost every $\sigma$.
\end{pf}

For further use, let us remark that an argument analogous to the proof of
Theorem~\ref {derivative bound} also gives:

\begin{prop} \label {derivativebound}

Let $A \in \Delta_\delta$ satisfy
\be
\|\op_z A_z\|=O(|\Im(z)|^{1+\epsilon}).
\ee
Then, for every $\sigma_0 \in \R/\Z$, if
\be
\limsup_{\sigma \to \sigma_0}
\frac {|\rho_{A_\sigma}-\rho_{A_{\sigma_0}}|}
{|\sigma-\sigma_0|}<\infty
\ee
then
\be
\limsup_{t \to 0} \frac
{|L(A_{\sigma_0+it})-L(A_{\sigma_0})|}
{|t|}<\infty.
\ee

\end{prop}

\subsection{Proof of Theorem \ref {dependence}}

\comm{
\begin{thm}

Let $A_{\theta,s} \in C^0(\R/\Z,\SL(2,\R))$, $\theta \in \R/\Z$, $s$ a
one-dimensional parameter, be monotonic in $\theta$ and
$C^{2r+1+\epsilon}$, $1 \leq r<\infty$, $C^\infty$ or $C^\omega$,
in $(\theta,s)$.
Then
\be
s \mapsto \int_{\R/\Z} L(A_{\theta,s}) d\theta
\ee
is $C^r$, $C^\infty$, or $C^\omega$.

\end{thm}
}

\comm{
We will need the following:

\begin{lemma}

If $A_z \in \Delta_\delta$ is $\eta$-asymptotically holomorphic with
$\eta>1$, then $\int_{\R/\Z} L(A_{\sigma+it}) d\sigma$
is a continuous function of $0\leq t<\delta$.

\end{lemma}

\begin{pf}

By Lemma \ref {lemma:2.5}, $\|\op_z m^+(z,x)\|=O(|\Im z|^{\epsilon-1})$ in
$\Omega_+$, where $\epsilon=\eta-1$.  Thus we can write a decomposition
$m^+=m^+_h+m^+_c$ where $z \mapsto m^+_c(z,x)$ is continuous up to $\R/\Z$
(uniformly in $x$,
and $m^+_h$ is of course uniformly bounded for $z$ near $\R/\Z$.  Thus
$m^+_h$ admits non-tangential limits as $\Im z \to 0$ and hence $m^+$ also
does.  It follows that the right hand side of the expression
\be
\int_{\R/\Z} L(A_{\sigma+it}) d\sigma=\int_{X \times \R/\Z} \ln
|\tau_{A_{\sigma+it}(x)}(m^+(\sigma+it,x))| d\sigma d\mu(x),
\ee
originally defined for $0<t<\delta$, makes sense up to $t=0$ and defines a
continuous function of $0 \leq t<\delta$.  Thus we just
need to identify the right hand side with
$\int_{\R/\Z} L(A_\theta) d\theta$ when $t=0$.

By Fubini's Theorem, for almost every
$\theta$, the non-tangential limits $m^+(\theta,x)$ exist, and provide an
invariant section for the cocycle $A_\theta$.  We claim that in this case
$\int_X \ln
|\tau_{A_\theta(x)}(m^+(\theta,x))| d\mu(x)=L(A_\theta)$.
Notice that the left hand side is non-negative.  Assuming ergodicity, by the
Oseledets Theorem the left hand side must be $\pm L(A_\theta)$, and by
non-negativity it must be $L(A_\theta)$.  The general case reduces to this
one by ergodic decomposition.
\end{pf}
}

It is enough to consider the case of finite differentiability, since we have
already proved the analytic case in section \ref {simple2}.

Let $A_{z,s} \in \Delta_\delta$ be
an asymptotically holomorphic extension of $A_{\theta,s}$
satisfying
\be
\|\pa^k_s A_{z,s}(x)\|=O(1), \quad 0 \leq k \leq 2r,
\ee
\be
\|\pa^k_s \op_z A_{z,s}(x)\|=O(|\Im(z)|^{2r-k+\epsilon}),
\quad 0 \leq k \leq 2 r.
\ee

\comm{
Define
\be
U(t,s)=\int_{\R/\Z} L(A_{\sigma+it,s}) d\sigma=
\int_{X \times \R/\Z} \ln |\tau_{A_{\sigma+it,s}(x)}(m^+_s(\sigma+it,x))| d\mu(x)
d\sigma.\label{eq:2.28}
\ee
}
Define $U(t,s)$ by formula (\ref{eq:2.28}).  It still satifies (\ref
{eq:2.29}) for $t \neq 0$.

We must show that $s \mapsto U(0,s)$ is $C^k$.
By Lemma \ref {0leqt}, $t \mapsto U(t,s)$ is continuous up to $t=0$, for
each $s$.  For $t \neq 0$, the functions $s \mapsto U(t,s)$ are in $C^k$,
so we just
need to show that as $t \to 0$, those functions converge uniformly in $C^k$
if $0 \leq k \leq r$.

We have the estimate
\be \label {lntau}
|\pa^k_s \op_z \zeta_s(z)|=O(|\Im(z)|^{2r-2k-1+\epsilon}),
\quad 0 \leq k \leq r-1.
\ee
Fix $0<\delta'<\delta''<\delta$.  Use Stoke's Theorem
to integrate $\zeta_s(z) dz$
on $[0,1] \times [\delta',\delta'']$, and take the imaginary part to get
\begin{align}
\frac {1} {2 \pi} U(\delta',s)=&\frac {1} {2 \pi} U(\delta'',s)+
\int_{X \times \{\delta'<\Im z<\delta''\}}
\op_z \zeta_s(z)
d\overline z \wedge dz\\
\nonumber
&+\int_X \int_{\delta'}^{\delta''} \Re \zeta_s(1+it)-
\Re \zeta_s(it) dt d\mu(x)
\end{align}
The first term on the right hand side is a fixed $C^k$ function of $s$.  The
second is also in $C^k$ and converges in $C^k$ as $\delta' \to 0$ by (\ref
{lntau}).  Notice that
$\Re \zeta_s(1+it)-\Re \zeta_s(it)$ is independent of $s$ and $t$:
it is the integral over $x$ of the topological degree of $A_{z,s}$ as $z$
runs once around $\Omega^+_\delta$.
Thus the third term is a linear
function $\delta''-\delta'$ (and independent of $s$).
This shows that $U(\delta',s)$ converges uniformly in $C^k$, as desired.
\qed

\comm{
We have the estimate
\be \label {lntau}
|\pa^k_s \op_z \ln \tau_{A_{z,s}(x)}(m^+_s(z,x))|=O(|\Im(z)|^{r-2k-1}),
\quad 0 \leq k \leq r-1.
\ee
Fix $0<\delta''<\delta$.  Use Stoke's Theorem
to integrate $\ln \tau_{A_{z,s}(x)}(m^+_s(z,x)) dz$
on $[0,1] \times [\delta',\delta'']$, and take the real part to get
\begin{align}
U(\delta',s)=&U(\delta'',s)+
\int_{X \times \{\delta'<\Im z<\delta''\}}
\op_z \ln \tau_{A_{z,s}(x)}(m^+_s(z,x))
d\overline z \wedge dz\\
\nonumber
&+\int_X
\int_{\delta'}^{\delta''} \Im (\ln
\tau_{A_{1+it,s}(x)}(m^+_s(1+it,x))-
\ln \tau_{A_{it,s}(x)}(m^+_s(1+it,x))) dt d\mu(x)
\end{align}
The first term on the right hand side is a fixed $C^k$ function of $s$.  The
second is also in $C^k$ and converges in $C^k$ as $\delta' \to 0$ by (\ref
{lntau}).  Since
$\tau_{A_{1+it,s}(x)}(m^+_s(1+it,x))=
\tau_{A_{it,s}(x)}(m^+_s(1+it,x))$
the integrand of the third term belongs to $2 \pi \Z$, so it
is independent of $(s,t)$.
Thus the third term is a linear
function of $\delta''-\delta'$ (and independent of $s$).\footnote {In fact
it is $2 \pi \deg (\delta''-\delta')$
where $\deg$ is the variation of the fibered rotation
number for the family $\theta \mapsto A_{\theta,s}$ as $\theta$ runs once
over $\R/\Z$.}
This shows that $U(\delta',s)$ converges uniformly in $C^k$, as desired.
}

\comm{
as $\Im \ln
\frac {\tau_{A_{1+it,s}(x)}(m^+_s(1+it,x))}
{\tau_{A_{it,s}(x)}(m^+_s(1+it,x))}$ and using that

\be
\pa^k_s U(\delta',s)=\par^k_s U(\delta'',s)+
\int_{X \times \{\delta'<\Im z<\delta\}} \pa^k_s
\op_z \ln |\tau_{A_{z,s}(x)}(m^+_s(z,x))|
\frac{d\overline z \wedge dz}{\theta-z}+\pa^k_s (2 \pi t \deg)
\ee
converge as continuous functions of $s$
when $\delta' \to 0$, for $0 \leq k \leq r$.

\bignote{I removed a factor 2 and replaced Re by Im; check consistency}
}

\subsection{$L^2$-estimates}

We continue with some crucial $L^2$-estimates, which play an important role
in renormalization theory, obtaining a somewhat more precise form of Theorem
\ref {l2ae}, Theorem \ref {L^2}.

First, let us recall the basic connection between estimates for
$\D$-valued invariant sections
and the existence of conjugacy to rotations.

\begin{lemma}

Let $A:X \to \SL(2,\R)$ be measurable.
The following are equivalent:
\begin{enumerate}
\item There exists a measurable $B:X \to \SL(2,\R)$ such that
$\int_X \|B(x)\|^2 d\mu(x)<\infty$ and
$B(f(x))A(x)B(x)^{-1} \in \SO(2,\R)$ for almost every $x$,
\item There exists a measurable $m:X \to \D$ such that $\int_X
\frac {1} {1-|m(x)|^2} d\mu(x)<\infty$ and $\mA(x) \cdot m(x)=m(f(x))$ for
almost every $x$.
\end{enumerate}

\end{lemma}

\begin{pf}

Let $B(x)$ be such that $\mB(x)=\frac {1} {(1-|m(x)|^2)^{1/2}}
\begin{pmatrix}1&-m(x)\\ \overline {-m(x)}&1\end{pmatrix}$.
\end{pf}

If $A \in C^0(\R/\Z,\SL(2,\R))$ satisfies the equivalent
conditions of the previous lemma, we will say that the corresponding cocycle
is $L^2$-conjugate to a cocycle of rotations.

Our aim in this section is to prove the following.

\begin{thm} \label {L^2}

Let $A_\theta \in C^0(X,\SL(2,\R))$,
$\theta \in \R/\Z$, be $C^{2+\epsilon}$ and monotonic in $\theta$.
For every $\theta \in \R/\Z$, if
\be \label {l2condi}
\limsup_{\theta' \to \theta} \frac
{|\rho_{A_{\theta'}}-\rho_{A_\theta}|}
{|\theta'-\theta|}<\infty
\ee
(in particular if $\theta \mapsto \rho_{A_{\theta}}$
is Lipschitz) and $L(A_\theta)=0$
then $A_\theta$ is $L^2$-conjugate to a cocycle of rotations.

\end{thm}

We will need a simple compactness result:

\begin{prop}\label{prop:2.12}

Let $A_k \in C^0(X,\Upsilon)$ be a sequence
converging to $A$.  Assume there exists
measurable functions $m_k:X \to \D$ satisfying $\mA_k(x) \cdot
m_k(x)=m_k(f(x))$, such that $\liminf \int_X \frac {1}
{1-|m_k(x)|^2} d\mu(x)<\infty$.
Then there exists a measurable
$m:X \to \D$ such that $\mA(x) \cdot m(x)=m(f(x))$
and $\int_X \frac {1} {1-|m(x)|^2} d\mu(x)<\infty$.\footnote {More
generally, we can also let the dynamics vary, thus considering a sequence of
$\mu$-preserving homeomorphisms $f_k:X \to X$ converging to $f$.}

\end{prop}

The proof uses the notion of conformal barycenter \cite {DE}, and we leave
it for the Appendix~\ref {conformal barycenter}.

\begin{cor}

Let $A \in \Delta_\delta$.  If $\sigma_0 \in \R/\Z$ satisfies
\be
\liminf_{t \to 0} \frac {L(\sigma_0+it)} {t}<\infty
\ee
then $A_{\sigma_0}$ is $L^2$-conjugate to a cocycle of rotations.

\end{cor}

\begin{pf}

The corollary follows from the previous Proposition \ref{prop:2.12} and Lemma~\ref {key computa}.
\end{pf}

\noindent {\it Proof of Theorem \ref {L^2}.}
It is enough to apply Proposition~\ref {derivativebound} and the corollary
above.
\qed

\noindent {\it Proof of Theorem \ref {l2ae}.}
Since $\theta \mapsto \rho_{A_\theta}$ is monotonic (\ref {l2condi}) holds
for almost every $\theta$.  We can thus apply Theorem \ref {L^2}.
\qed

\begin{rem}

If one is only concerned with a result valid for almost every $\theta$ (as
in the case of Theorem \ref {l2ae}), one
can bypass the use of the conformal barycenter argument.  Indeed, the most
usual argument in such situations is to apply the Lemma of Fatou to
guarantee convergence of $m^+(\sigma+it,x)$ as $t \to 0+$ for almost every
$x$, and then apply Fubini's Lemma to obtain a set of $\sigma$ of full
Lebesgue measure for which $\lim_{t \to 0+} m^+(\sigma+it,x)$ exists for
almost every $x$.

\end{rem}

\subsection{Proof of Theorem \ref {rigidity}}

\comm{
\begin{thm}

Let $A_\theta \in C^0(X,\SL(2,\R))$ be monotonic and
$C^{r+1+\epsilon}$, $0 \leq r<\infty$, $C^\infty$ or $C^\omega$,
in $\theta$.  If $L(A_\theta)=0$
for every $\theta$ in some open interval $J$
then there exists $B_\theta \in C^0(X,\SL(2,\R))$, $\theta \in J$
depending $C^r$, $C^\infty$ or $C^\omega$
on $\theta$ and conjugating $A_\theta$ to a
cocycle of rotations.

\end{thm}
}

We will consider only the finitely differentiable case, the other cases
following basically the same argument.
We shall assume that $J=\R/\Z$ for simplicity.  Consider an asymptotically
holomorphic extension of $A_\theta$ satisfying
\be
\|\op_z A_z\|=O(|\Im(z)|^{r+\epsilon}).
\ee
Then we have
\be
\|\op_z m^+(z,x)\|=O(|\Im(z)|^{r-1+\epsilon}), \quad \Im(z)>0
\ee
and analogously
\be
\|\op_z m^-(z,x)\|=O(|\Im(z)|^{r-1+\epsilon}), \quad \Im(z)<0.
\ee
Let
\be
\phi(z,x)=\op_z m^\pm(z,x), \quad \pm \Im(z)>0,
\ee
and let $u:\C/\Z \times X \to \C$ be given by
\be
u(z,x)=\lim_{t \to \infty} \frac {-1} {\pi} \int_{[-t,t] \times
[-\delta,\delta]} \frac {\phi(w,x)} {z-w} dx dy.
\ee
A compactness argument shows that $u(z,x)$ is continuous on both variables. 
Moreover, $\R/\Z \ni y \mapsto u(y,x)$ is $C^r$ (uniformly in $x$).
Let
\be
m(z,x)=m^\pm(z,x), z \in \Omega^\pm_\delta.
\ee
Then $\lim_{t \to 0} m(\sigma+it,x)$
exists for almost every $\sigma$ and almost every $x$ by Lemma \ref {key
computa}.  Thus for almost every $x \in X$, $z \mapsto m(z,x)-u(z,x)$
extends to a holomorphic function defined on $\Omega_\delta$.  A compactness
argument shows that this holds indeed for all $x \in X$, and that the
function $\Omega_\delta \times \R/\Z \ni
(z,x) \mapsto m(z,x)-u(z,x)$ is continuous (as in the classical De
Concini-Johnson argument \cite {CJ}).  It also follows that $\R/\Z \ni
y \mapsto m(y,x)$ is $C^r$ (uniformly on $x$).  To conclude, it is enough to
show that $m(y,x)$ takes values on $\D$. 

Let us assume first that $f:X \to X$ is minimal.  If $y$ is such that
$m(y,x_0) \in \partial \D$ for some $x_0 \in X$,
then for every $x \in X$
we also have $m(y,x) \in \partial \D$ (by
invariance).
However, since $L(A_y)=0$ for every $y$,
$\rho_{A_y}$ is $C^1$ (by Schwarz Reflection and $r \geq 1$),
so by Proposition \ref {derivativebound} and Lemma \ref {key
computa}, for every $y$
we have
\be
\limsup_{t \to 0+} \int_{\R/\Z} \frac {1} {1-|m^+(y+it,x)|^2}
d\mu(x)<\infty,
\ee
so by continuity $m(y,x) \in \D$ for almost every $x$.

Let us now consider the general case.  Notice that if $\mu'$ is any ergodic
invariant probability measure, the Lyapunov exponent $L'(A_y)$ with respect
to $\mu'$ is still $0$ for every $y \in \R/\Z$.  Indeed,
\be
L'(A_z)=\pm \int \ln |\tau_{A_z}(m(z,x))| d\mu'(x) \geq 0, \quad z \in
\Omega^\pm_\delta,
\ee
and since $m$ is continuous, we have $\int
\ln |\tau_{A_y}(m(y,x))|d\mu'(x)=0$.  On the other hand, since $m$ is an
invariant section, if $L'(A_y) \neq 0$ then
$\int \ln |\tau_{A_y}(m(y,x))|d\mu'(x)$ coincides with $L'(A_y)$ up to sign, as
desired.

Since any minimal set supports an ergodic invariant measure, we
can apply the previous argument to conclude that $m(y,x) \in \D$ whenever
$x$ belongs to a minimal set.  To conclude, notice that
for each $y$, the set of all $x$ such that
$m(y,x) \in \partial \D$ is compact and invariant, so if it is non-empty it
must contain a minimal set.\qed

\comm{
Similar arguments yield the analytic and infinitely differentiable cases, so
we will not get into their proof:

\begin{thm}

Let $A_\theta \in C^0(X,\SL(2,\R))$ be
$C^r$, $r=\omega,\infty$
and monotonic in $\theta$.  If $L(A_\theta)=0$
for every $\theta$ in some open interval $J$
then there exists $B_\theta \in C^0(X,\SL(2,\R))$, $\theta \in J$
depending $C^r$ on $\theta$ and conjugating $A_\theta$ to a
cocycle of rotations.

\end{thm}
}

\section{Monotonic cocycles}\label{sec:3}

We now turn to the study of quasiperiodic cocycles presenting monotonicity
in phase space.  In this section, $d \geq 1$ is a fixed integer and the
dynamics $f:\R^d/\Z^d \to \R^d/\Z^d$ is a translation
$f(x)=x+\alpha$, where $\alpha$ is assumed to be fixed except when otherwise
noted.  The underlying probability measure will be Lebesgue
measure.

Given $w \in \R^d$, we shall say that
$A \in C^0(\R^d/\Z^d,\SL(2,\R))$ is $w$-{\it monotonic}
if $A^w_\theta(x)=A^w(x+\theta w)$ is monotonic decreasing.

We say that $A$ is {\it monotonic} if it is $w$-monotonic for some $w$ \footnote{We apologize for this collision of definition with the notion of $\epsilon$-monotonic cocycle.}.
If $A$ is $C^1$ then the set of
$w$ such that $A$ is $w$-monotonic is an open convex cone, hence if
non-empty it contains primitive vectors of $\Z^d$, and up to linear
automorphism of $\R^d/\Z^d$ it contains $(1,0,...,0)$.  Notice that the
monotonicity condition implies that for fixed $(x_2,...,x_d)$,
$x_1 \mapsto A(x_1,x_2,...,x_d)$ is a map $\R/\Z \to \SL(2,\R)$
with negative topological degree.  
Thus monotonic cocycles are never
homotopic to a constant.  Indeed the phenomena we will uncover for monotonic
cocycles collide often with the intuition developed for
Schr\"odinger cocycles, which are homotopic to a constant.


A straightforward application of our previous results yields:

\begin{thm} \label {3.1}

Let $A \in C^r(\R^d/\Z^d,\SL(2,\R))$, $r=\omega,\infty$
be monotonic.  If $L(A)=0$ then $A$ is
$C^r$-conjugate to a cocycle of rotations.

\end{thm}

\begin{pf}

If $A$ is $w$-monotonic, then the family $A^w_\theta(x)$ is
obviously $C^r$ in $\theta$ and satisfies $L(A^w_\theta)=L(A)$ for every
$\theta$.  If $L(A)=0$, the proof of Theorem \ref {rigidity} gives
$m^w_\theta \in C^0(\R^d/\Z^d,\D)$, $C^r$ in $\theta$, such that
$\mA^w_\theta(x) \cdot m^w_\theta(x)=m^w_\theta(x+\alpha)$.

We claim that $m^w_\theta(x)=m^w(x+\theta w)$, where $m^w=m^w_0$.
Recall that $m^w_\theta(x)$
is obtained, almost everywhere, as a non-tangential limit
$\lim_{\epsilon \to 0} m^w_{\theta+\epsilon i}(x)$ where
$m^w_{\theta+\epsilon i}$ is the unique $\D$-valued invariant section
of an asymptotically holomorphic extension of $A^w_\theta$.  This
asymptotically holomorphic extension can be chosen here to satisfy
$A^w_{\theta+\epsilon i}(x)=A^w_{\epsilon i}(x+\theta w)$: this is in
fact automatic if we follow the procedure described in section
\ref{asymptotically}
for the construction of the asymptotically holomorphic extension.
In this case, we get
$m^w_{\theta+\epsilon i}(x)=m^w_{\epsilon i}(x+\theta w)$, and
hence this equality
is satisfied almost surely by the non-tangential limit.  Since it is
continuous, the claim follows.

Let us now consider another monotonicity vector $w'$.  We claim that
$m^{w'}=m^w$.  Let us first assume ergodicity.  In this case, if
$m^w \neq m^{w'}$
at some $x \in \R^d/\Z^d$ then this must happen at every $x$, and in fact
the hyperbolic distance between $m^{w'}(x)$ and $m^w(x)$ must be some
constant $c$,
by ergodicity of $f$ and the fact that the projective action preserves the Poincar\'e distance of the disk.
Then we can define $B:\R^d/\Z^d \to \PSL(2,\R)$ by $\mB(x) \cdot m^w(x)=0$ and
$\mB(x) \cdot m^{w'}(x) \in t i$, where $t>0$ is such that the hyperbolic
distance between $0$ and $t i$ is $c$.  It follows that
$\mB(f(x)) \mA(x) \mB(x)^{-1}$ takes $0$ to $0$ and $t i$ to $t i$, so
$B(f(x)) A(x) B(x)^{-1}$ is the identity in $\PSL(2,\R)$.  This is clearly
impossible since $A$ is non-homotopic to a constant.
In the non-ergodic case, the open set $U$ where $m^{w'}(x) \neq m^w(x)$
is not necessarily everything, but it is foliated by periodic subtori where
the dynamics is ergodic.  The previous argument shows that there exists a
continuous $B:U \to \PSL(2,\R)$ such that $B(f(x))A(x)B(x)^{-1}=\id$.
By monotonicity, for each $x \in U$ there exists $C=C(x)>0$ such that for
every $\theta>0$, for any line $l \in \P\R^2$, and for every $k \in \Z$,
there exists $\gamma=\gamma(x,k,\theta,l)$ with
$C^{-1} \theta<\gamma<C \theta$ and
$B(f^{k+1}(x)) A(f^k(x+\theta w)) B(f^k(x))^{-1} \cdot l=R_\gamma \cdot l$. 
This shows that for every $\epsilon>0$, $x \in U$, if $\theta>0$ is
sufficiently small then for any line $l \in \P\R^2$ the sequence
$(B(f^n(x))A_{n}(x+\theta w)B(x)^{-1})_{n \in \Z}$ is $\epsilon$-dense in
$\P\R^2$.  But this is impossible since for small $\theta>0$ we have that
$B(f^n(x))A_n(x+\theta
w)B(x)^{-1}=B(f^n(x)) B(f^n(x+\theta w))^{-1} B(x+\theta w) B(x)^{-1}$ is
close to the identity for every $n \in \Z$.


\comm{
Consider an homotopy between $B$ and a map
$B' \in C^0(\R^d/\Z^d,\PSL(2,\R))$ such that
$B'(x)$ is a rotation for every $x$, with angle depending affinely on $x$.
It gives an homotopy between $A$ and $\pm B'(f(x)) B(x)^{-1}$, which is
constant.  But since $A$ is monotonic it can not be homotopic to a
constant, and the claim follows.
}

We conclude that there exists a single $m \in C^0(\R^d/\Z^d,\D)$ which
coincides with all $m^w$'s.  Thus for each $w$ in an open cone,
$\theta \mapsto m(x+\theta w)$ is $C^r$.
By Journ\'e's Theorem \cite {Jo},
it is $C^r$ as a function of $x$.
\end{pf}

\comm{
\begin{thm}

Let $A \in C^{r+1+\epsilon}(\R/\Z,\SL(2,\R))$, $0 \leq r <\infty$
be a monotonic cocycle.  If $L(\alpha,A)=0$ then $(\alpha,A)$ is
$C^r$-conjugate to a cocycle of rotations.

\end{thm}
}

\begin{thm} \label {3.2}

Let $A_s \in C^0(\R^d/\Z^d,\SL(2,\R))$ be a one-parameter family which is
$C^r$ in $x$ and $s$, $r=\omega,\infty$.
If $A_{s_0}$ is monotonic then $s \mapsto L(A_s)$ is $C^r$ in a neighborhood of $s_0$.

\end{thm}

\begin{pf}

Let $s_0$ and $w$ be such that $A_{s_0}$ is $w$-monotonic.  We may assume
that $w$ is a primitive vector of $\Z^d$.  Consider the two-parameter family
$A_{\theta,s}$, $\theta \in \R/\Z$ given by $A_{\theta,s}(x)=A_s(x+w \theta)$. 
By Theorem \ref {dependence}, the $\theta$-average of
$L(A_{\theta,s})$ is $C^r$ in $s$
near $s_0$.  But $L(A_{\theta,s})=L(A_s)$ for every $s$, which gives the
result.
\end{pf}

\comm{
analytic
Let us consider a one-parameter family $A^s \in
\R \times C^\omega(\R/\Z,\SL(2,\R))$ of monotonic cocycles.  If $(s,x)
\mapsto A^s(x)$ is $C^\omega$ then
\be
s \mapsto L(\alpha(s),A^s)
\ee
is $C^\omega$.

\end{thm}
}

\begin{rem} \label {skewshift}

Much of our analysis generalizes to some other dynamical systems, including
the usual skew-shift $(x,y) \mapsto (x+\alpha,y+x)$. 
Consider a skew-product $f:X \times \R^d/\Z^d \to X \times
\R^d/\Z^d$, $f(x,y)=(\phi(x),y+\psi(x))$, and define a cocycle over
$f$ to be monotonic if there exists $w \in \R^d$ such that
$y \mapsto A(x,y)$ is $w$-monotonic for every $x \in X$.  Our arguments
imply that if $A$ is monotonic and $C^r$, $r=\omega,\infty$,
with respect to the second coordinate
then $L(A)=0$ implies that $A(x,y)$ is $C^0$ conjugated
to rotations (the conjugacy being $C^r$ in $y$). 
We can also show that if $A_s(x,y)$ is a family of monotonic cocycles
which is $C^r$ with respect to $(s,y)$ then
$s \mapsto L(A_s)$ is $C^r$.

\end{rem}

\subsection{Varying the frequency}

Let us briefly allow the dynamics to vary, in order to obtain a result about
the regularity of the Lyapunov exponent with respect to such more general
perturbations.

We first need a replacement for Lemma \ref {lemma:2.5}.
Given a one-parameter family $A_{z,s} \in
\Delta_\delta$ which is monotonic, and a $C^1$
one-parameter family of frequencies $\alpha(s) \in \R^d$, define
invariant sections $m^+_s(z,x) \in \D$ so that $\mA_{z,s}(x) \cdot
m^+_s(z,x)=m^+_s(z,x+\alpha(s))$.

\begin{lemma}

Assume that $\alpha(s)$ and
$(s,x) \mapsto A_{z,s}(x)$ are $C^r$, $1 \leq r<\infty$ and
\be
\|\pa^i_s \pa^j_x A_{z,s}(x)\|=O(1), \quad 0 \leq k=i+j \leq r.
\ee
Then
\be
|\pa^i_s \pa^j_x m^+_s(z,x)|=O(|\Im(z)|^{-2 k-i+1)}),
\quad 1 \leq k=i+j \leq r.
\ee
Moreover, if additionally $s \mapsto \op A_{z,s}(x)$ is $C^{r-1}$ and
we have the estimate
\be
\|\pa^i_s \pa^j_x \op_z A_{z,s}(x)\|=O(|\Im(z)|^{\eta-k-1}), \quad
0 \leq k=i+j \leq r-1,
\ee
for some $\eta \in \R$ then
\be
|\pa^i_s \pa^j_x \op_z m^+_s(z,x)|=O(|\Im(z)|^{\eta-2k-i-2}), \quad
0 \leq k=i+j \leq r-1.
\ee

\end{lemma}

\begin{pf}

One can basically repeat the proof of Lemma \ref {lemma:2.5},
since the added complications are not very serious.  Write
$m^s_z(x)$ for $m^+_s(z,x)$ and
$F^s_z(x,w)$ for $\mA_{z,s}(x) \cdot w$.

Consider first the first estimate.  Notice that when $i=0$ things reduce to
Lemma \ref {lemma:2.5}.  Let us assume by
induction (first on $k$ and then on $i$)
that we have the desired bounds when $i'+j'<k$ and also for
$i'+j'=k$, $0 \leq i'<i$.
Consider the derivative $\partial^i_s \partial^j_x$ of
$m^s_z(x+\alpha(s))=F^s_z(x,m^s_z(x))$.  The left hand side has a main
term of the form $(\pa^i_s \pa^j_x m^s_z)(x+\alpha(s))$
and a lower order term of the form
\be
\sum_{l<i,l+n \leq k}^i C \cdot (\pa^l_s (\pa^n_x m^s_z)(x+\alpha(s)),
\ee
with $C$ polynomials (depending on the indices)
on the derivatives of $\alpha$ of
order at most $i$.  The lower order term can be estimated by induction to
be $O(|\Im z|^{-(2k+i-2)})$.  The right hand side has a
main term of the form
$\pa_w F^s_z(x,m^s_z(x)) \cdot \pa^i_s \pa^j_x m^s_z(x)$ and a lower order
term of the form
\be
\sum C \cdot \prod_{n=1}^l \pa^{i_n}_s \pa^{j_n}_x m^s_z(x),
\ee
where the sum runs over all $0 \leq l \leq k$, and $(i_1,j_1) \leq ... \leq
(i_l,j_l)$ (lexicographic order) with $i_n,j_n \geq 0$, $1 \leq n \leq l$,
$\sum_{n=1}^l i_n=i' \leq i$,
$\sum_{n=1}^l j_n=j' \leq j$ and 
$i'+j'<k$, 
and $C$ are now constant
multiples of $\partial^{i-i'}_s \partial^{j-j'}_x \pa^l_w F$.  The lower
order term can be estimated by induction to be $O(|\Im z|^{-(2k+i-2)})$.
Rearranging we get
\be
|(\pa^i_s \pa^j_x) m^s_z(f+\alpha(s))-
\pa_w F^s_z(x,m^s_z(x)) \cdot \pa^i_s \pa^j_x m^s_z(x)|=O(|\Im
z|^{-(2k+i-2)}),
\ee
and the first estimate follows from Proposition \ref {prop:2.6}.

Consider now the second estimate.
Let us assume by
induction that we have the desired bounds when $i'+j'<k$ and also for
$i'+j'=k$, $0 \leq i'<i$.
We consider the derivative $\partial^i_s \partial^j_x \op_z$ of
$m^s_z(x+\alpha(s))=F^s_z(x,m^s_z(x))$.  The left hand side has a main term
of the form $(\pa^i_s \pa^j_x \op_z m^s_z)(x+\alpha(s))$
and a lower order term
\be
\sum_{l<i,l+m \leq k}^i C\cdot 
(\pa^l_s \pa^m_x \op_z m^s_z)(x+\alpha(s)),
\ee
with $C$ polynomials (depending on the indices)
on the derivatives of $\alpha$
of order at most $i$.  The lower order term can be estimated by induction to
be $O(|\Im z|^{\eta-2k-i-1)})$.
The right hand side has one main term
$\pa_w F^s_z(x,m^s_z(x)) \cdot \pa^i_s \pa^j_x \op_z m^s_z(x)$ and two lower
order terms.  The first has the form
\be
\sum C \cdot \prod_{n=1}^l \pa^{i_n}_s \pa^{j_n}_x m^s_z(x),
\ee
where the sum runs over all $0 \leq l \leq k$, and $(i_1,j_1) \leq ... \leq
(i_l,j_l)$ (lexicographic order) with $i_n,j_n \geq 0$, $1 \leq n \leq l$,
$\sum_{n=1}^l i_n=i' \leq i$,
$\sum_{n=1}^l j_n=j' \leq j$, and $C$ are now constant
multiples of $\partial^{i-i'}_s \partial^{j-j'}_x \op_z
\pa^l_w F$.  Using the first estimate, we see that it is $O(|\Im
z|^{\eta-2k-i})$ except when $k=0$, where it is $O(|\Im z|^{\eta-1})$.
The second has the form
\be
\sum D\cdot  \pa^{i_0}_s \pa^{j_0}_x \op_z m^s_z(x)
\prod_{n=1}^l \pa^{i_n}_s \pa^{j_n}_x m^s_z(x),
\ee
where the sum runs over all $0 \leq l \leq k$, $0 \leq i_0 \leq i$, $0 \leq
j_0 \leq j$, $i_0+j_0<k$, and $(i_1,j_1) \leq ... \leq
(i_l,j_l)$ (lexicographic order) with $i_n,j_n \geq 0$, $1 \leq n \leq l$,
$\sum_{n=0}^l i_n=i' \leq i$,
$\sum_{n=0}^l j_n=j' \leq j$, and $D$ are constant
multiples of $\pa^{i-i'}_s \pa^{j-j'}_x \pa^l_w F$.  This term can be
estimated using the first estimate and the induction hypothesis to be
$O(|\Im z|^{\eta-2k-i-1})$.
Rearranging we get
\be
|(\pa^i_s \pa^j_x \op_z) m^s_z(f+\alpha(s))-
\pa_w F^s_z(x,m^s_z(x)) \cdot \pa^i_s \pa^j_x \op_z m^s_z(x)|=O(|\Im
z|^{\eta-2k-i-1}),
\ee
and the second estimate follows from Proposition \ref {prop:2.6}.
\end{pf}

\comm{
order term can be estimated by induction to be $O(|\Im z|^{-(2k+i-2)})$.
\be
\sum C \prod_{n=1}^l \pa^{i_n}_s \pa^{j_n}_x m^s_z(x),
\ee
with $C$ polynomials (depending on the indices)
on the derivatives of $\alpha$
of order at most $i$.  The lower order term can be estimated by induction to
be $O(|\Im z|^{\eta-(2k+i-2)})$.

The right hand side has one ``main
term'' $\partial_w F^s_z(x,m^s_z(x)) \cdot
(\partial^i_s \partial^j_x m^s_z)(x+\alpha(s))$ and ``lower order'' terms
which are estimated by induction to be $O(|\Im z|^{-k})$.  The left hand
side is slightly different from the previous computation:
it has one ``main term''
$(\partial^i_s \partial^j_x m^s_z)(x+\alpha(s))$, a secondary term
$(\partial_s \alpha(s))^l
\sum_{l=1}^{s-1} (\partial^{i-l}_s \partial^{j+l}_x m^s_z)(x+\alpha(s))$,
and ``lower order'' terms which are estimated by induction to be $O(|\Im
z|^{-k}$.  If $i=0$, the secondary term disappears, and we get the estimate
just as before, as a consequence of Proposition.  This allows us to start a
new

$\partial_w F^s_z(x,m^s_z(x)) \cdot
(\partial^i_s \partial^j_x m^s_z)(x+\alpha(s))$ and ``lower order'' terms
\end{pf}
}

\begin{thm} \label {3.4}

Let $A_{\theta,s} \in C^0(\R^d/\Z^d,\SL(2,\R))$ be monotonic on $\theta$
and $C^\infty$ on $\theta,s$, and let $s \mapsto \alpha(s)$ be $C^\infty$.
Then the average with respect to $\theta$ of the
Lyapunov exponent of $A_{\theta,s}$ over
$x \mapsto x+\alpha(s)$ is a $C^\infty$ function of $s$.

\end{thm}

\begin{pf}

Fixing $k \geq 0$, choose $\eta$ large and
let $A_{z,s} \in \Delta_\delta$ be an $\eta$-asymptotically holomorphic
extension of $A_{\theta,s}$.
As in the proof of Theorem \ref {dependence}, we define for every $s$ and
$0<t<\delta$
\be
U(t,s)=\int_{\R/\Z} L(A_{\sigma+it,z}) d\sigma.
\ee
Then for $0<t<\delta$ the
functions $s \mapsto U(t,s)$ are $C^k$ and as $t \to 0$ they converge
uniformly in $C^k$.  For each fixed $s$, the limit is seen to be the
$\theta$-average of the Lyapunov exponent of $A_{\theta,s}$:
since $s$ is fixed, we can just apply Lemma \ref {0leqt}.
\end{pf}

\begin{rem}

Even if everything is analytic, we do not, in general,
get analytic dependence when varying the frequency by this argument. 
Analytic dependence should not be expected, since when the frequency is
complexified the domain of analyticity does not remain invariant by the
dynamics.
A special case where 
 analytic dependence holds is when   $A_{\theta,s}=R_\theta A_s$, since the dynamics
does not influences the $\theta$-average of the Lyapunov exponent in this
case.

\end{rem}

\begin{thm} \label {3.5}

If $A_s \in C^\infty(\R^d/\Z^d,\SL(2,\R))$ is monotonic and $C^\infty$ on
$x$ and $s$, and $s \mapsto \alpha(s)$ is $C^\infty$, then
the Lyapunov exponent of $A_s$ as a cocycle over
$x \mapsto x+\alpha(s)$ is a $C^\infty$ function of $s$.

\end{thm}

\begin{pf}

Assume that $A$ is $w$-monotonic with $w$ a primitive vector of $\Z^d$ and
consider the family $A_{\theta,s}=A_s(x+\theta w)$.
\comm{
Fixing $k \geq 0$, choose $\eta$ large and
let $A^s \in \Delta_\delta$ be an $\eta$-asymptotically holomorphic
extension of $A_{\theta,s}(x)=A_s(x+\theta w)$.  As in the proof of
Theorem, we define
\be
U_{t,s}=\int_{\R^d/\Z^d \times \R/\Z} \ln
|\tau_{A_{\sigma+it,s}(x)}(m^s_z(x))| dx d\sigma,
\ee
defined for $0<t<\delta$ is such that for each $0<t<\delta$, the
functions $s \mapsto U(t,s)$ are $C^k$ and as $t \to 0$ they converge
uniformly in $C^k$.  For each fixed $s$, the limit is seen to be the
$\theta$-average of the Lyapunov exponent of $A_{\theta,s}$, regarded as a
cocycle over $x \mapsto x+\alpha(s)$: since $s$ is fixed, the argument is
the same as in the proof of Theorem.
}
The Lyapunov exponents of
$A_{\theta,s}$ and of $A_s$ (both considered as cocycles over the same $x
\mapsto x+\alpha(s)$) are obviously equal  for every $\theta$.  The
result follows by the previous theorem.
\end{pf}

\comm{
$s \mapsto U_{0,s}=\lim_{t \to 0} U(t,s)$ exists
uniformly
\be
\lim \pa^k_s U_{t,s}

$|\partial^k \overline \partial_z \zeta_s| \leq |\Im z|^{-2}$ for
the function $\zeta_s$,
which allows us to conclude using Stokes Theorem that
$\partial^k_s \int_{\R/\Z} L(A_{\theta+\epsilon i,s}) d\theta$ converges
uniformly to a continuous function of $s$.
\end{pf}

Let us consider the $X=(\R^d/\Z^d)^2$ and let $f:X \to X$
$f(\alpha,x)=(\alpha,x+\alpha)$
and $A^w_\theta:X \to \SL(2,\R)$, $A^w_\theta(\alpha,x)=A(x+\theta w)$.

Let us consider the cocycle generated by $f$ and $A^w_\theta$.  Taking an
asymptotically holomorphic extension of $A^w_\theta$, and applying Lemma
with zero-dimensional parameter space,
we obtain functions
$m^+:\Omega^+ \times X \to \D$ such that $\mA^w_z \cdot
m^+(z,\alpha,x)=m^+(z,\alpha,x+\alpha)$ satisfying
\be
|\partial^k m^+(z,\alpha,x)|=O(|\Im z|^{-k}), \quad 0 \leq k \leq n,
\ee
\be
|\partial^k \partial_z m^+(z,\alpha,x)|=O(|\Im z|^n), \quad 0 \leq k \leq n,
\ee
where $n$ is fixed but arbitrarily large.  This allows us to apply 

\begin{thm}

Let us consider a one-parameter family $(\alpha(s),A^s) \in
\R \times
C^{2r+1+\epsilon}(\R/\Z,\SL(2,\R))$, $1 \leq r <\infty$
of monotonic cocycles.  If $s \mapsto
\alpha(s)$ is $C^{2r+1+\epsilon}$ and $(s,x)
\mapsto A^s(x)$ is $C^{2r+1+\epsilon}$ then
\be
s \mapsto L(\alpha(s),A^s)
\ee
is $C^r$.

\end{thm}
}

\subsection{Low regularity considerations}

We return to the consideration of a fixed translation dynamics, but now
focus on trying to obtain conclusions at low regularity.

Consider some $A \in C^0(\R^d/\Z^d,\SL(2,\R))$ which is $w$-monotonic, and
let $A^w_\theta(x)=A(x+\theta w)$.  It
follows directly from the definitions that
$\rho_{A^w_\theta}$ is an affine function of $\theta$ with
negative slope $\deg^w$.  In fact $\deg^w$
can be explicitly given as a linear function
of $w$ which only depends on topological data: $\deg^w=\langle l,w \rangle$
where $l \in \Z^d$ is the unique integer such that $A$ is homotopic to $x
\mapsto R_{\langle l,x \rangle}$.\footnote {To see this, define
$\delta_{v,w} \xi_n=\delta_\gamma \xi_n$ and
$\delta_{v,w} \zeta=\delta_\gamma \zeta$ where
$\gamma$ is any path homotopic to $\gamma(t,x)=A(x+v+t w)$.  Notice that
$\rho_{A^w_\theta}$ is given,
up to additive constant, by $\delta_{0,\theta w} \zeta$.  Let us show that
$\delta_{v,w} \zeta=\langle l,w \rangle$.  It is clear
that $\delta_{v,w} \xi_n(x,z_0,z_1)=\delta_{0,w}(x+v,z_0,z_1)$, so that
$\delta_{v,w} \zeta$ does not depend on $v$.  By Remark \ref {concate},
$\delta_{v,w} \zeta+\delta_{v+w,w'} \zeta=\delta_{v,w+w'} \zeta$.  This
shows that $\delta_{v,w} \zeta$ is a linear function of $w$.  Moreover, for
$w \in \Z^d$ we have exactly
$\delta_{v,w} \xi_n=\langle l,w \rangle n$, so that $\delta_{v,w}
\zeta=\langle l,w \rangle$.} 

In particular, $\theta \mapsto \rho_{A^w_\theta}$ is Lipschitz.
More generally, we have the following result.

\begin{lemma}

Let us consider a one-parameter family
$A_\theta \in C^0(\R^d/\Z^d,\SL(2,\R))$.  If for
some $\theta_0$, $A_{\theta_0}$ is a $w$-monotonic cocycle, and
\be
K=\limsup_{\theta \to \theta_0} \frac {1} {|\theta-\theta_0|}
\sup_x \|A_\theta(x)-A_{\theta_0}(x)\|<\infty
\ee
then
\be
\limsup_{\theta \to \theta_0} \frac {1} {|\theta-\theta_0|}
|\rho_{A_\theta}-\rho_{A_{\theta_0}}| \leq K',
\ee
where $K'$ depends on $K$, the monotonicity constant of $A^w_{\theta_0}$,
$\|A_{\theta_0}\|_{C^0}$ and $\deg^w$.

\end{lemma}

\begin{pf}

The hypothesis imply that for $h$ close to $0$ and $z \in \partial \D$,
$\mA_{\theta_0+h}(x) \cdot z$ lies in the shortest segment of $\partial \D$
determined by $\mA_{\theta_0}(x-C h w) \cdot z$ and
$\mA_{\theta_0}(x+C h w) \cdot z$, for some $C>0$.  This implies
that $\rho_{A_{\theta_0+h}}$ lies between
$\rho_{A_{\theta_0}(\cdot+C h w)}$ and
$\rho_{A_{\theta_0}(\cdot-C h w)}$,\footnote {Indeed
we can construct a monotonic decreasing family $\tilde A_t \in
C^0(\R^d/\Z^d,\SL(2,\R))$, $t \in [0,1]$, such that $\tilde
A_0(\cdot)=A_{\theta_0}(\cdot-C h w))$,
$\tilde A_1(\cdot)=A_{\theta_0}(\cdot+C h w)$ and $\tilde
A_{1/2}=A_{\theta_0+h}$, and $\tilde A_t$ is close to $A_{\theta_0}$ for
every $t$ (so that we remain in a region where there is a continuous
determination of $\rho$).  Monotonicity of this family gives $\rho_{\tilde
A_1} \leq \rho_{\tilde A_{1/2}} \leq \rho_{\tilde A_0}$.}
that is, in the segment
$[\rho_{A_{\theta_0}}+C h \deg^w,\rho_{A_{\theta_0}}-C h
\deg^w]$, and the result follows.
\end{pf}

Thus, if $A$ is monotonic then $\theta \mapsto R_{-\theta} A$ is a monotonic
family with Lipschitz rotation number.  In low regularity, it may be
preferable to work with this family, because it is always analytic in
$\theta$.  As an application, we have the following result  which
is a direct consequence of Theorem \ref{L^2} (if we were to
use only the family $\theta \mapsto A^w_\theta$, we would need
to assume further regularity).

\begin{thm}

Let $A \in C^0(\R/\Z,\SL(2,\R))$ be monotonic.  If $L(A)=0$
then $A$ is $L^2$-conjugate to a cocycle of rotations.

\end{thm}

It also allows us to get continuity results in Lipschitz open sets of
cocycles.  For the family $R_{-\theta} A$, $e^{2 \pi i \theta}
\mapsto e^{-2 \pi i \rho(\theta)-L(\theta)}$, $\Im \theta>0$,
defines an univalent function from $\D \setminus \{0\}$ to $\D \setminus
\{0\}$ (see \cite {CJ}).
This gives an harmonic conjugacy relation between
the Lyapunov exponent and the fibered rotation number.  If the fibered
rotation number turns out to be Lipschitz, then the
Lyapunov exponent (as a function of the circle) has derivative in $L^1$
and in fact $\partial_\theta L$ is a zero-average
BMO function (since it is basically the Hilbert trasform of the
derivative of the fibered rotation number, which is in $L^\infty$).  This
argument also shows that the
BMO norm of $\partial_\theta L$ can be bounded in terms of the Lipschitz
constant of $\rho$.

\begin{thm}

Let $\epsilon>0$ be fixed.
The Lyapunov exponent is a continuous function of $\epsilon$-monotonic
cocycles (the frequency may be varied as well).

\end{thm}

\begin{pf}

Let $A^{(n)} \to A$ be a sequence of
$\epsilon$-monotonic cocycles converging in $C^0$.  We allow
$A^{(n)}$ to be regarded as cocycles over $x
\mapsto x+\alpha_n$ and $A$ over $x \mapsto x+\alpha$, as long as $\alpha_n
\to \alpha$.
\comm{
Then $\zeta_n(\theta)$ can be defined on the fundamental
cover of $\R/\Z$ to satisfy $\zeta_n(\theta)=\zeta_n(\theta+1)$.  Since
$\Re(\zeta_n)$ is non-decreasing, its
extension to the upper half plane is seen to be a $1$-periodic
univalent map $\H$ to some subset of $\H$, with $\zeta_n(z)=z+O(1)$ at
infinity.  Identifying $\H \setminus \{0\}$
}
By the  previous lemma, there exists $C>0$ such that the
fibered rotation number of $\theta \mapsto
R_{-\theta} A^{(n)}$ is a 
$C$-Lipschitz function for all $n$, and by the discussion
above, $\theta \mapsto L_n(\theta)=L(R_{-\theta} A^{(n)})$
is uniformly equicontinuous.
\comm{
$\theta \mapsto L_n(\theta)=
L(\alpha_n,R_\theta A^{(n)})$ are Hilbert
transforms of (uniformly bounded) Lipschitz functions
(the rotation numbers), so they belong to a
compact set of continuous functions (their derivatives being uniformly in
BMO).
}
Thus we may assume $L_n \to L_\infty$ in $C^0(\R/\Z,\R)$.
By upper semicontinuity of the Lyapunov exponent,
$\theta \mapsto L(R_{-\theta} A)-L_\infty(\theta)$ is a
non-negative continuous
function, which we must show to be identically zero.
By \cite {AB},
\begin{align}
\int_{\R/\Z} L(&R_{-\theta} A)-L(R_{-\theta} A^{(n)}) d\theta\\
\nonumber
&=\int_{\R/\Z} \ln \frac
{\|A(x)\|+\|A(x)\|^{-1}} {2}-\ln \frac {\|A^{(n)}(x)\|+\|A^{(n)}(x)\|^{-1}}
{2} dx,
\end{align}
so $\int_{\R/\Z}
L(R_{-\theta} A)-L_\infty(\theta)d\theta=
\lim_{n \to \infty} \int_{\R/\Z}
L(R_{-\theta} A)-L(R_{-\theta} A^{(n)}) d\theta=0$.
\end{pf}

\subsection{Minimality} \label {minim}

Here we are interested in considering the dynamics of cocycles
from the topological point of view.  For this, one considers
the cocycle given by some $A$
as a map, still denoted by $(f_\alpha,A)$, from
$\R^d/\Z^d \times \partial \D \to \R^d/\Z^d \times
\partial \D$ (a $d+1$-dimensional torus) given by $(x,w) \mapsto
(x+\alpha,\mA(x) \cdot w)$, where $\partial \D$ is identified with $\P\R^2$
in the usual way.
Below we consider only the case where  $x \mapsto x+\alpha$ is ergodic.

It can be shown (see \cite {KKHO}) that if $A \in C^0(\R/\Z,\SL(2,\R))$
is not homotopic to the identity then for every $\alpha \in \R \setminus
\Q$, $(f_\alpha,A)$ is transitive, and an adapted argument for the
multidimensional case is given in Appendix \ref {transitivity}.
Though proofs of transitivity
just involve simple topological arguments, the following question seems
much harder:

\begin{problem}

Is $(f_\alpha,A)$ minimal whenever
$A \in C^0(\R^d/\Z^d,\SL(2,\R))$ is non-homotopic to the identity?

\end{problem}

Of course the same problem
still makes sense under additional smoothness assumptions.
In this section we will show that the complexification
methods allow one to address the local case, at least if one assumes enough
smoothness.


Let us first discuss some known results
on the minimal sets of non-uniformly hyperbolic cocycles (we follow the
presentation of Herman \cite {H}, but the results are due to Johnson
\cite {J}).
If $A \in C^0(\R/\Z,\SL(2,\R))$
and $L(A)>0$ then it follows from Oseledets Theorem that there exist two
measurable functions $u,s:\R^d/\Z^d \to \partial \D$ (the unstable and
stable directions) such that
$\mA(\theta) \cdot u(\theta)=u(\theta+\alpha)$ and $\mA(\theta)
\cdot s(\theta)=s(\theta+\alpha)$ and for almost every $\theta$, for every
$w \in \overline \C$, if $w \neq s(\theta)$ then $|\mA_n(\theta) \cdot
w-u(\theta+n\alpha)| \to 0$ exponentially fast, and if $w
\neq u(\theta)$ then $|\mA_n(\theta-n\alpha)^{-1} \cdot w-s(\theta-n\alpha)|
\to 0$ exponentially fast as $n\to\infty$.  It follows (from unique ergodicity of $\theta
\mapsto \theta+\alpha$) that there are exactly two ergodic
invariant measures on $\R^d/\Z^d \times \partial \D$,
the push-forwards of
Lebesgue measure on $\R^d/\Z^d$ by
$\theta \mapsto (\theta,u(\theta))$ and
$\theta \mapsto (\theta,s(\theta))$, which we denote by $\mu_u$ and $\mu_s$. 
Let us denote their (compact) support by $K_u$ and $K_s$.  It follows that
any minimal set for $(f_\alpha,A)$ coincides with either $K_u$ or
$K_s$.  Moreover, if $K_u \neq K_s$ then $A$ would necessarily be uniformly
hyperbolic.
So assuming $A$ to be non-homotopic to the identity with
$L(A)>0$ we get that $K_u=K_s$ is the unique minimal set of $(f_\alpha,A)$.

\comm{
they are
automatically minimal and there can be no further minimal sets for
$(\alpha,A)$.

\begin{prop}[Herman, \cite {H}]

Let $(\alpha,A) \in (\R \setminus \Q) \times C^0(\R/\Z,\SL(2,\R))$ be
non-uniformly hyperbolic.  Then $K_u=K_s$.

\end{prop}
}


\begin{thm} \label {3.9}

Let $A \in C^{1+\epsilon}(\R^d/\Z^d,\SL(2,\R))$ be monotonic.
Then $(f_\alpha,A)$ is minimal.

\end{thm}

\begin{pf}

If $L(A)=0$ then $(f_\alpha,A)$ is $C^0$ conjugate
to a cocycle of rotations, by the argument of the proof of Theorem \ref
{rigidity}.\footnote {Theorem \ref {rigidity} gives a $C^1$ conjugacy under
a slightly stronger $C^{2+\epsilon}$ condition.  Under $C^{1+\epsilon}$, it
still gives a continuous
invariant section $m:\R^d/\Z^d \to \overline \D$, which
implies the $C^0$ conjugacy, unless the invariant section is real, i.e., it
lies in $\partial \D$.  However, this last possibility is impossible here,
since $A$ is non-homotopic to a constant so it
can not admit invariant continuous sections in $\partial \D$.}

For a cocycle of rotations, transitivity
obviously implies minimality, so the result follows from
Proposition \ref {tra}.

Let now $L(A)>0$.  We consider the analytic case, the smooth case
being analogous.  Up to coordinate change, we may assume that $A$ is
$w$-monotonic where $w=(1,0,...,0)$.  Let $m:\Omega^+_\delta \times
\R^{d-1}/\Z^{d-1} \to \D$ satisfy $\mA(z) \cdot
m(z)=m(z+\alpha)$.  Then for every $x_2,...,x_n$ and for almost every $x_1$,
$m(x_1,...,x_n)=\lim_{t \to 0} m(x_1+t i,x_2,...,x_n)$ exists.
Since $L(A)>0$, for almost every
$x \in \R^d/\Z^d$, $m(x) \in \partial \D$
and $m|\R^d/\Z^d$
coincides with either the unstable or stable directions $u$, $s$
defined above.\footnote {It is easy to see that it actually coincides with
$u$, but this will not play a role here.}

We claim
that for every open set of the form $J \times U \subset \R^d/\Z^d$ with $J
\subset \R/\Z$, and any interval
$J' \subset \partial \D$, there exists a positive measure set of
$x \in J \times U$ such
that $m(x) \in J'$: by the previous discussion, the unique minimal set
$K_u=K_s$ of $(f_\alpha,A)$ must intersect $J \times U \times J'$, and since it is
arbitrary we must have $K_u=K_s=\R^d/\Z^d \times \partial \D$.

Suppose by contradiction that the claim does not hold.  Up to a change of
coordinates, we may assume that $-1\in  J'$.  Then
for almost every $y \in U$,
$z \mapsto (1-m(z,y))(1+m(z,y))$ is a holomorphic function on
$\Omega^+_\delta$ with positive real part, bounded near $J$ and whose
non-tangential limits are purely imaginary on
$J$; by the Schwarz Reflection Principle there
is a holomorphic extension to $\Omega_\delta \setminus (\R/\Z \setminus
J)$.  Thus $x_1 \mapsto m(x_1,...,x_n)$ is analytic on $J$
for almost every $(x_2,...,x_n) \in U$.  By invariance we have $x_1
\mapsto m(x_1,...,x_n)$ analytic on $\R/\Z$ for
almost every $(x_2,...,x_n)$.  The topological
degree of $x_1 \mapsto m(x_1,...,x_n)$
is an integer valued measurable
function $\deg(x_1,x_2,...,x_n)$ which does not depend on $x_1$.
But $\deg(x+\alpha)=\deg(x)+\deg$, where $\deg$ is the topological
degree of $x_1 \mapsto \mA(x_1,0,...,0) \cdot 1$.
Since $A$ is $w$-monotonic with $w=(1,0,...,0)$, $\deg<0$.
But by Poincar\'e recurrence
$\deg(x+n \alpha)$ takes the same value infinitely many times, for almost every
$x$, so $\deg=0$, contradiction.
\end{pf}

\subsection{Premonotonic cocycles}

As remarked in the introduction, the concept of monotonicity is not
dynamically natural.
The easiest way to extend the concept of monotonicity is the following.
We say that a cocycle $A \in C^1(\R^d/\Z^d,\SL(2,\R))$
is {\it premonotonic} if some iterate is $C^1$ conjugate to a
monotonic cocycle: there exist $n \geq 1$ and
$B \in C^1(\R^d/\Z^d,\SL(2,\R))$ such that
$B(x+n\alpha)A_n(x)B(x)^{-1}$ is monotonic.
This happens if and only if some iterate of $(f_\alpha,A)$ is
real-analytic conjugate to a monotonic cocycle (any $C^1$-perturbation of
$B$ which is real analytic will do).  Notice that premonotonic cocycles are
$C^1$-stable, and there is even stability with respect to
perturbations of the frequency vector defining the dynamics in the basis.

Cocycles of rotations over ergodic translations
which are not homotopic to a constant provide the
simplest examples of premonotonic cocycles.
Given a cocycle $A \in
C^0(\R^d/\Z^d,\SL(2,\R))$, let $[A]$ be the unique cocycle homotopic to $A$
of the form $[A](x)=R_{\langle l,x \rangle}$ with $l=l^A \in \Z^d$.  Notice
that $A$ is not homotopic to a constant if and only if $l \neq 0$, and in
this case $[A]$ is $l$-monotonic.

\begin{lemma} \label {conjpre}

Let $A \in C^r(\R^d/\Z^d,\SL(2,\R))$, $1 \leq r \leq \infty$ or $r=\omega$,
and let $x \mapsto x+\alpha$ be ergodic on $\R^d/\Z^d$.  If $(f_\alpha,A)$ is
$C^r$-conjugate to a cocycle of
rotations then there exists a sequence $B^{(n)} \in
C^r(\R^d/\Z^d,\SL(2,\R))$ such that
$A^{(n)} \to [A]$ in the $C^r$-topology, where
$A^{(n)}(x)=B^{(n)}(x+\alpha) A(x) B^{(n)}(x)^{-1} \in \SO(2,\R)$.

\end{lemma}

\begin{pf}

By definition, we may assume that $A$ is itself a cocycle or rotations (the
homotopy class being clearly conjugacy invariant).  Thus let
$A(x)=[A](x) R_{\phi(x)}$, where $\phi:\R^d/\Z^d \to \R$ is $C^r$.
Let $c=\int \phi(x) dx$.  Since $\alpha$ is irrational, we can consider a
sequence $l_k \in \Z^d$ such that $\langle \alpha,l_k \rangle \to c \mod 1$.

Let us consider a sequence $\phi^{(n)}:\R^d/\Z^d \to \R$ of trigonometric
polynomials converging to $\phi$ in $C^r$ and with $\int \phi^{(n)}(x)
dx=c$.
Since $x \mapsto x+\alpha$ is an
ergodic translation, it is easy to define, using Fourier series,
trigonometric polynomials $\psi^{(n)}:\R^d/\Z^d \to \R$
such that
$\phi^{(n)}(x)=-\psi^{(n)}(x+\alpha)+\psi^{(n)}(x)+c$.  Letting
$C^{(n)}(x)=R_{\psi^{(n)}(x)}$,
we see that $C^{(n)}(x+\alpha) R_{\phi(x)} C^{(n)}(x)^{-1}$ is $C^r$
close to $R_c$.  Consider now
$B^{(n)}(x)=R_{-\langle l_{k_n},x \rangle} C^{(n)}(x)$.  If we choose $k_n
\to \infty$ very slowly, we will have $B^{(n)}(x+\alpha) A(x)
B^{(n)}(x)^{-1} \to [A](x)$ in $C^r$, as desired.
\end{pf}


\comm{
\begin{lemma}

Let $A \in C^1(\R^d/\Z^d,\SL(2,\R))$ be non-homotopic to a constant and let
$x \mapsto x+\alpha$ be ergodic on $\R^d/\Z^d$.
If $(f_\alpha,A)$ is $C^1$-conjugated to a cocycle of rotations then
$(f_\alpha,A)$ it is $C^1$-conjugated to a monotonic cocycle.

\end{lemma}

\begin{pf}

By the previous lemma, $A$ can be $C^1$-conjugated arbitrarily close to
$[A]$.
But $[A]$ is
monotonic (in fact $w$-monotonic for every $w$ with $\langle w,l^A
\rangle>0$), so any $C^1$ nearby cocycle is also monotonic.
\end{pf}
}

\comm{
By definition, we may assume that $A$ is itself a cocycle or rotations (the
homotopic class being clearly conjugacy invariant).  Thus let
$A(x)=R_{\phi(x)}$, where $\phi:\R^d/\Z^d \to \R/\Z$ is $C^1$.  We can write
$\phi(x)=\langle l,x \rangle+P$ where $l \in
\Z^d$, and $P:\R^d/\Z^d \to \R$ has zero average.  Moreover, since $A$
is not homotopic to a constant, $l \neq 0$.

Let us consider a sequence $P^{(n)}:\R^d/\Z^d \to \R$ of trigonometric
polynomials converging to $P$ in $C^1$.
Since $x \mapsto x+\alpha$ is an
ergodic translation, it is easy to define, using Fourier series,
trigonometric polynomials $\psi^{(n)}:\R^d/\Z^d \to \R$
such that
$P^{(n)}(x)=\psi^{(n)}(x+\alpha)-\psi^{(n)}(x)$.  Letting
$B^{(n)}(x)=R_{\phi^{(n)}(x)}$,
we see that $B^{(n)}(x+\alpha) A(x) B^{(n)}(x)^{-1}$ is $C^1$
close to $R_{\langle l,x \rangle}$.  Since $l \neq 0$, $x \mapsto R_{\langle
l,x \rangle}$ is monotonic (indeed $w$-monotonic for every $w$ with $\langle
w,l \rangle>0$).  It follows that $B^{(n)}(x+\alpha) A(x) B^{(n)}(x)^{-1}$
is monotonic for every $n$ sufficiently large, so $A$ is premonotonic.
\end{pf}
}


\comm{
Although premonotonicity is (expressily) invariant under smooth conjugacies,
we have not succeeded in proving it is invariant under more general classes
of transformations.  For instance, it is not clear if
a cocycle which admits a monotonic iterate is premonotonic (though it is
easy to show that there exist real-analytic
premonotonic cocycles with any given
frequency which do not admit a monotonic iterate).  Worse, this definition
does not behave well under renormalization.  Thus we were led to study some
more general classes of cocycles which can be shown to admit a description
similar to monotonic cocycles.  Our results are, at the moment, not
completely satisfactory: we have identified a natural class which is
invariant under renormalization, but we did not succeed in
showing that this class is actually bigger than the class of premonotonic
cocycles.  In order not to distract from the normal flow of our arguments,
we have left this discussion for Appendix~\ref {premonotonicity}.  In the
same appendix we shall also prove the existence of many cocycles which
are not premonotonic.}

While premonotonicity is only
{\it a priori} invariant under $C^1$-conjugacies, we have:

\def\Leb{{\mathrm{Leb}}}

\begin{thm} \label {c1c0}

Let $A \in C^1(\R^d/\Z^d,\SL(2,\R))$ be non homotopic to the identity, and let $f_\alpha:x\mapsto x+\alpha$ be
ergodic on $\R^d/\Z^d$.  If $(f_\alpha,A)$ is $C^0$-conjugate to a cocycle
of rotations, then $(f_\alpha,A)$ is premonotonic.

\end{thm}

\begin{pf}


Up to isometric automorphism of
$\R^d/\Z^d$, we may assume that the first coordinate $l_1$
of $l=l^A$ is positive.

Denote by $B:\R^d/\Z^d\to \SL(2,\R)$ the $C^0$-map such that $R(\cdot):=B(\cdot+\alpha)A(\cdot)B(\cdot)^{-1}$ takes it values in the group of rotations.
We will identify $\P\R^2$ with $\R/\Z$ so that the projective action of
rotations corresponds to translations.
Let $F_n:\R^d/\Z^d \times \P\R^2 \to \R^d/\Z^d \times \P\R^2$ be the
projective action, $F_n(x,y)=(x+n\alpha,A_n(x) \cdot y)$.
For $(x,y) \in \R^d/\Z^d \times \P\R^2$,
let $a_n(x,y)=\partial_{x_1} (A_n(x) \cdot y)$, $b_n(x,y)=\partial_y (A_n(x)
\cdot y)$, $p(x,y)=\partial_y (B(x) \cdot y)$, $q_n=a_n p \circ F^n$,
$q=q_1$.
We claim that  $b_n=\frac {p} {p \circ F^n}$ and
\be
q_n=\sum_{k=0}^{n-1} q \circ F^k.\label{3.16}
\ee
Indeed,
from the definition of $B$ and $R$, one can  write for any $x$, $B(x+n\alpha)\cdot (A_{n}(x)\cdot y)=R_{n}(x)\cdot (B_{n}(x)\cdot y)$ and thus taking derivatives with respect to $y$, $p(x+n\alpha,A_{n}(x)\cdot y)b_{n}(x,y)
=p(x,y)$ which is the relation $b_{n}=\frac{p}{p\circ F^n}$. For (\ref{3.16}) we just write $A_{n+1}(x)\cdot y=A(x+n\alpha)\cdot (A_{n}(x)\cdot y)$ and take derivatives with respect to $x_{1}$ to  get $a_{n+1}=a_{1}\circ F^n+(b_{1}\circ F^n) a_{n}$; since we have just seen that $b_{1}\circ F^n=(p\circ F^n)/(p\circ F^{n+1})$ we have  $a_{n+1}(p\circ F^{n+1})=(a_{1}p\circ F)\circ F^n+a_{n}(p\circ F^n)$ which obviously gives (\ref{3.16}).

Let $e_1=(1,...,0)$.
Notice that $a_n(x,y)>0$ for all $(x,y)$ is equivalent to $-e_1$-monotonicity of
$A_n$.  Below we will prove that $q_n \to \infty$, and hence $a_n \to
\infty$, uniformly in $(x,y)$, giving the premonotonicity of $(f_\alpha,A)$.

For $x_0,x \in \R^d/\Z^d$, $y \in \P\R^2$
\be
d^\epsilon_n(x_0,x,y)=
\int_0^\epsilon p(x_0+n \alpha,A_n(x+t e_1) \cdot y) a_n(x+t e_1,y) dt
\ee
gives the oriented length of the path
$\gamma|[0,\epsilon]$, where $\gamma=\gamma_{n,x_0,x,y}:\R \to \P\R^2$ is
given by $\gamma(t)=B(x_0+n \alpha) A_n(x+t e_1) \cdot y$
(the oriented length can be defined
as $\hat \gamma(\epsilon)-\hat \gamma(0)$ where $\hat \gamma:\R
\to \R$ is a lift to the universal cover).  Especially, for any $y,y'
\in \P\R^2$ we must have
\be
-1<d^\epsilon_n(x_0,x,y)-d^\epsilon_n(x_0,x,y')<1,
\ee
since when $y \neq y'$ we must have $\gamma_{n,x_0,x,y}(t) \neq
\gamma_{n,x_0,x,y'}(t)$ for every $t \in \R$.

Let $H(x,y)=(x,B(x) \cdot y)$ and
let $G=H \circ F \circ H^{-1}$.  Since $G$ is topologically conjugate to
$F$, and
Lebesgue measure on $(x,y)$ is invariant for $G$, $H^{-1}_*
\Leb$ is an invariant measure for $F$ equivalent to $\Leb$.
Thus for Lebesgue almost every $(x,y)$,
\be
\hat q(x,y)=\lim \frac {1} {n} \sum_{k=0}^{n-1} q \circ F^k(x,y)
\ee
exists.  Since $\hat q$ is measurable, Lebesgue almost every $(x,y)$ is a
measurable continuity point along the $x_1$ direction.  Especially, for
almost every $(x,y)$ we have
\be
\lim_{\epsilon \to 0} \frac {1} {\epsilon} \lim_{n \to \infty} \frac {1} {n}
d^\epsilon_n(x,x,y)=\hat q(x,y).
\ee
Thus $\hat q(x,y)$ is almost surely independent of $y$, and since $x \mapsto
x+\alpha$ is ergodic, $\hat q(x,y)$ is almost surely independent of $x$ and
$y$.

Notice that $G$ commutes with shifts in the second coordinate
$T_t(x,y)=(x,y+t)$.  Thus any ergodic
invariant measure $\mu$ for $G$ gives rise to a
one-parameter family of ergodic invariant measures $\mu_t=(T_t)_* \mu$. 
By unique ergodicity of $x \mapsto x+\alpha$,
all those measures project down to Lebesgue measure on $x$.
It follows that $\int_{\R/\Z} \mu_t dt=\Leb$.  By uniqueness of the
ergodic decomposition, it
follows that all ergodic invariant measures are of the form $\mu_t$, for
some $t \in \R$ (for any fixed $\mu$).

Since $\mu_t$ depends continuously on $t \in \R/\Z$ (with respect to the
weak-$*$ topology) and $\int_{\R/\Z} \mu_t
dt=\Leb$, the fact that $\hat q \circ H^{-1}=\lim \frac {1} {n}
\sum_{k=0}^{n-1} q \circ H^{-1} \circ G^k$ is almost everywhere
constant implies that $\int q \circ H^{-1} d\mu_t$ is independent of $t$. 
Since $q \circ H^{-1}$ has constant average with respect to all ergodic
invariant measures, the Birkhoff averages
of $q \circ H^{-1}$ converge uniformly
to a constant limit.  Thus $\frac {q_n} {n} \to \int q \circ H^{-1} d\Leb$
uniformly.

To conclude, we must show that $\int q \circ H^{-1} d\Leb>0$.  If this is
not the case, then for every $n$ sufficiently large we will have
$\frac {a_n} {n}<\frac {1} {2}$.  But this is impossible because of the
identity $\frac {1} {n} \int_0^1 a_n(t,x_2,...,x_n,y) dt=l_1$ which is a
positive integer.
\end{pf}

The definition of premonotonicity is such that results proved for
monotonic cocycles extend easily to this larger setting.  Let us comment in
more detail on the results stated in the introduction which involve
premonotonicity (except for Theorem \ref {premonotonicpositive}, which we
discuss in the next section).

Theorem \ref {lyapanal} follows from Theorem \ref {3.2}
and Theorem \ref {lyapsmooth} follows from Theorem \ref {3.5} (as the Lyapunov
exponent is well behaved when taking conjugacies and iterates).

In order to derive Theorem \ref {l=0cr} from Theorem \ref {3.1}, it is enough
to notice that if a $C^r$ cocycle $(f_\alpha,A)$ admits an iterate which is
$C^r$ conjugate to a cocycle of rotations, then $(f_\alpha,A)$ is itself
$C^r$-conjugate to a cocycle of rotations.\footnote {This is most easily
seen by working with $C^r$ invariant sections (which arise from
and give rise to
a conjugacy to rotations in the usual way).  If $m \in C^r(\R^d/\Z^d,\D)$
satisfies $\mA_n(x) \cdot m(x)=m(x+n\alpha)$, let $m_j(x)=\mA_j(x-j\alpha)
\cdot m(x-j\alpha)$.  Then $m_{j+n}=m_j$ and $\mA(x) \cdot
m_j(x)=m_{j+1}(x)$.  For each $x \in \R^d/\Z^d$, let $m_*(x)$ minimizes the
sum of the squares of the hyperbolic distances (in $\D$) to
$(m_j(x))_{j=0}^{n-1}$: this is a well defined $C^r$ function of $x$ by
strict convexity.  Then $\mA(x) \cdot m_*(x)=m_*(x+\alpha)$.}

The proof (if not the statement) of Corollary \ref {c0cr} also
involves premonotonicity: it
follows from Theorem \ref {l=0cr} and Theorem \ref {c1c0}.

Theorem \ref {preminimal} follows from Theorem \ref {3.9} (since minimality of any
iterate implies minimality).


\subsection{Non-uniform hyperbolicity for typical premonotonic cocycles}



In this section, we will only consider, for simplicity, the case of
$C^r$ cocycles with $r=\infty$ or $\omega$.  Then the Lyapunov
exponent is indeed a $C^r$ function of
premonotonic cocycles (while we have only carried out the formal arguments
for the dependence of the Lyapunov exponent along one-parameter families, it
is clear the estimates go through to the infinite dimensional
parametrization).
Since the Lyapunov exponent $L$
takes non-negative values, we must have $DL=0$ whenever $L=0$.
Here we are going to show that, in the case of premonotonic cocycle,
if $L=0$ then $D^2 L \neq 0$.  This implies
that $\{L=0\}$ is a subvariety of positive codimension in the space of
premonotonic cocycles and completes the proof of
Theorem \ref {premonotonicpositive}.

If $B \in C^r(\R^d/\Z^d,\SL(2,\R))$, then the conjugacy operator $A \mapsto
A'$, $A'(x)=B(x+\alpha) A(x) B(x)^{-1}$ is a $C^r$ diffeomorphism in
$C^r(\R^d/\Z^d,\SL(2,\R))$.  Since the Lyapunov exponent is clearly
invariant by conjugacy, it suffices to check that any premonotonic cocycle
$A$ with $L(A)=0$
is conjugate to some $A'$ such that $D^2 L(A') \neq 0$.  But
$C^r$ premonotonic cocycles with zero Lyapunov exponent are $C^r$ conjugate
to cocycles of rotations by Theorem \ref {l=0cr}, and in fact, by Lemma \ref
{conjpre},
those may be chosen arbitrarily close to a cocycle of the form
$x \mapsto [A](x)=R_{\langle l,x \rangle}$ with $l \neq 0$.
Since the Lyapunov exponent is $C^2$
near $[A]$, it suffices to show $D^2 L([A]) \neq 0$.
We will in fact give a simple estimate implying the existence
of cocycles near $[A]$ with a quadratic lower bound on the Lyapunov
exponent.

For a matrix $s \in \mathrm{sl}(2,\R)$, let
$s_1,s_2,s_3$ be such that $s=\bm s_1&s_2+s_3\\s_2-s_3 & -s_1 \em$.

\begin{lemma}

Let $l \in \Z^d \setminus \{0\}$, $s \in C^0(\R/\Z,\mathrm{sl}(2,\R))$, and
define $A_{\theta,t} \in C^0(\R^d/\Z^d,\SL(2,\R))$, $\theta \in \R/\Z$, $t \in
\R$ by
$A_{\theta,t}(x)=R_{\langle l,x \rangle} e^{t s(\langle l,x \rangle-
\theta)}$.
Then
\be
\lim_{t \to 0} \frac {2} {t^2} \int_{\R/\Z} L(A_{\theta,t}) d\theta=\int_{\R/\Z}
s_1^2(\theta)+s_2^2(\theta) d\theta.
\ee
In particular, the limit is zero if and only if $s$ takes values in
$\so(2,\R)$.

\end{lemma}

\begin{pf}

Let $C_{t,\theta}(x)=R_\theta C_t$ where
$C_t(x)=R_{\langle l,x \rangle} e^{t s(\langle l,x \rangle)}$.
Notice that $A_{\theta,t}(x+\theta \frac {l} {\|l\|^2})=C_{t,\theta}(x)$.  So
$L(A_{\theta,t})=L(C_{t,\theta})$.  By \cite {AB},
\be
\int_{\R/\Z} L(C_{t,\theta}) d\theta=\int_{\R^d/\Z^d} \ln \frac
{\|C_t(x)\|+\|C_t(x)\|^{-1}} {2}
dx=\int_{\R/\Z} \ln \frac {\|e^{t s(\theta)}\|+\|e^{t s(\theta)}\|^{-1}} {2}
d\theta.
\ee
On the other hand, a direct computation shows that
\be
\lim_{t \to 0} \frac {2} {t^2} \int_{\R/\Z}
\ln \frac {\|e^{t s(\theta)}\|+\|e^{t s(\theta)}\|^{-1}} {2}
d\theta=
\int_{\R/\Z}
s_1^2(\theta)+s_2^2(\theta) d\theta.
\ee
The result follows.
\end{pf}

Choosing, say, $s_1(\theta)=\cos 2 \pi \theta$, $s_2=s_3=0$, we see that
the family $A_{\theta,t}$ is an analytic family (on $\theta$ and $t$) of
analytic cocycles such that $A_{0,t}$ is constant equal to $x \mapsto
R_{\langle l,x \rangle}$.  The previous lemma then implies that $D^2 L$
(in either setting, analytic or smooth)
does not vanish on $x \mapsto R_{\langle l,x \rangle}$, as desired.

\comm{
\begin{lemma}\label{lemma:4.1}

Let $B \in C^0(\R^d/\Z^d,\SL(2,\R))$.  Let $A_\theta(x)=R_{n x}
B(x-\theta)$.  Then
\be
\int_0^1 L(\alpha,A_\theta(x)) d\theta=\int^1_0 \ln \left
(\frac {\|B(x)\|+\|B(x)\|^{-1}} {2} \right ) dx.
\ee

\end{lemma}

\begin{pf}

Let $C_\theta(x)=R_{n \theta} R_{n x} B(x)$.
Notice that $A_\theta(x+\theta)=C_\theta(x)$.  In
particular, $L(\alpha,A_\theta)=L(\alpha,C_\theta)$.  The result follows by
\cite {AB}.
\end{pf}

Let $s \in C^0(\R/\Z,\Sl(2,\R))$, that is,
\be
s(x)=\begin{pmatrix}a(x)&b(x)+c(x)\\b(x)-c(x)&-a(x)\end{pmatrix},
\ee
where $a,b,c:\R/\Z \to \R$ are continuous functions.
Let $A_{\theta,t}=R_{n x} e^{t s(x-\theta)}$.  Then the previous lemma
implies that
\be
\lim_{t \to 0} \frac {1} {2 t^2} \int_0^1 L(\alpha,A_{\theta,t})
d\theta=\int_0^1 a^2(x)+b^2(x) dx.
\ee
In particular, the limit above is zero if and only if $a$ and $b$ vanish
identically, that is, if and only if $s$ takes values in $\so(2,\R)$.
}

\section{One-frequency cocycles: renormalization and rigidity} \label {sec4}

We continue our investigations of quasiperiodic cocycles, but now specify to
the case of one frequency.  Though the number of frequencies is quite
irrelevant in the analysis of monotonic cocycles, in the one-frequency case
we will be able to obtain global consequences from our local analysis,
by means of renormalization, a tool that is not as effective when several
frequencies are involved.

Below we will only consider cocycles over irrational rotations.
To highlight the dependence on the base dynamics, through this section
a cocycle will be specified by a pair $(\alpha,A) \in (\R \setminus \Q)
\times C^0(\R/\Z,\SL(2,\R))$.

After defining the renormalization operator, we are going to show that if
$(\alpha,A) \in (\R \setminus \Q) \times C^1(\R/\Z,\SL(2,\R))$
is $L^2$-conjugate to rotations, then it admits a
``renormalization representative'' $(\alpha',A') \in (\R \setminus \Q)
\times C^1(\R/\Z,\SL(2,\R))$
(seen as a cocycle over some different irrational rotation)
with $A'$ $C^1$-close to
$x \mapsto R_{\theta+\deg x}$ for some $\theta$
(here $\deg$ is the topological
degree of $A$).  Moreover, if $A$ is $C^r$,
$A'$ can be chosen to be $C^r$.
The dynamics
of $A$ and $A'$ can
be related, in particular if $L(\alpha,A)=0$ then $L(\alpha',A')=0$
and if $(\alpha',A')$ is $C^r$ conjugate to rotations then $(\alpha,A)$
is also $C^r$-conjugate to rotations.

Now, if $A$ is not
homotopic to a constant, $\deg \neq 0$, so $A'$ is monotonic.  This leads to
our main global rigidity result in the one-dimensional case, Theorem \ref
{globall2}.




Let us note that
by our analysis of one-parameter families, Theorem \ref {globall2} implies:

\begin{thm}

Let $(\alpha,A_\theta)$, $\alpha \in \R \setminus \Q$, $A \in
C^r(\R/\Z,\SL(2,\R))$, $r=\infty,\omega$, be a
one-parameter family which is monotonic and
$C^{2+\epsilon}$ in $\theta$.
If the $A_\theta$ are non-homotopic to a constant then for almost every
$\theta$, either $L(A_\theta)>0$ or $A_\theta$ is $C^r$ conjugate to
rotations.

\end{thm}

\begin{pf}

By Theorem \ref {l2ae},
for almost every $\theta$ with $L(\alpha,A_\theta)=0$,
$(\alpha,A)$ is $L^2$ conjugate to rotations.  By Theorem \ref {globall2}, 
they
must be actually $C^r$ conjugate to rotations.
\end{pf}

\subsection{Renormalization}

In this section we recall some basic facts on renormalization.
We refer to \cite{AK} for the proofs and further details.

Let $(\alpha,A) \in ((0,1) \setminus \Q) \times C^r(\R/\Z,\SL(2,\R))$ be a
cocycle
Let $p_n/q_n$ be the continued fraction
approximants of $\alpha$ and let $\beta_n=(-1)^n (q_n \alpha-p_n)$,
$\alpha_n=\beta_n/\beta_{n-1}$.  Thus $\alpha_n=G^n(\alpha)$ where
$G(\alpha)=\{\alpha^{-1}\}=\alpha^{-1}-[\alpha^{-1}]$ is the Gauss map.

Classically, the dynamical systems
$x \mapsto x+\alpha_n$ can be interpreted as
the sequence renormalization of $x \mapsto x+\alpha$.
We would like to produce, starting from $(\alpha,A)$,
a sequence $A^{(n)} \in C^0(\R/\Z,\SL(2,\R))$ such that
$(\alpha_n,A^{(n)})$ can be interpreted as the sequence of renormalizations
of $(\alpha,A)$.  However, this can not be done canonically, and to define
renormalization one must introduce {\it commuting pairs}.

Fixing $x_* \in \R/\Z$, we associate to $(\alpha,A)$ a sequence of pairs
$(A^{(n,0)},A^{(n,1)}) \in C^0(\R,\SL(2,\R))$, by
\be
A^{(n,0)}(x)=A_{(-1)^{n-1} q_{n-1}}(x_*+\beta_{n-1} x),
\ee
\be
A^{(n,1)}(x)=A_{(-1)^n q_n}(x_*+\beta_{n-1} x).
\ee
We should regard $A^{(n,0)}$ and $A^{(n,1)}$ as defining cocycles
over the dynamics on $\R$ given by
$x \mapsto x+1$ and $x \mapsto x+\alpha_n$.
It is easy to see that $A^{(n,1)}(x+1)
A^{(n,0)}(x)=A^{(n,0)}(x+\alpha_n) A^{(n,1)}(x)$, which expresses  the
commutation of the cocycles.  We call $((1,A^{(n,0)}),(\alpha_n,A^{(n,1)}))$
the $n$-th {\it renormalization} of $(\alpha,A)$ around $x_*$.

The dynamics of $A^{(n,0)}$ (and of $A^{(n,1)}$ as well)
is trivial, since all orbits go to infinity.
In fact we can always define (\cite {AK}, Lemma 4.1)
a (non-canonical) {\it normalizing
map} associated to $(1,A^{(n,0)})$, that is, some
$B^{(n)} \in C^0(\R,\SL(2,\R))$ such that $B^{(n)}(x+1) A^{(n,0)}(x)
B^{(n)}(x)^{-1}=\id$.

Because of the commutation relation, it follows that
if $B^{(n)}$ is a normalizing map for $(1,A^{(n,0)})$, then
$A^{(n)}(x)=B^{(n)}(x+\alpha_n) A^{(n,1)}(x) B^{(n)}(x)^{-1}$ satisfies
$A^{(n)}(x+1)=A^{(n)}(x)$.  Thus $A^{(n)}$ can be seen as an element of
$C^0(\R/\Z,\SL(2,\R))$, and $(\alpha_n,A^{(n)})$ is
called a {\it representative} of the $n$-th
renormalization of $(\alpha,A)$.

Of course, choosing a different normalizing map $\tilde B^{(n)}$
leads to a possibly different $\tilde A^{(n)}$.  But it is easy to see that
$C=\tilde B^{(n)} (B^{(n)})^{-1}$ is $1$-periodic,
which implies that $\tilde A^{(n)}(x)=C(x+\alpha_n) A^{(n)}(x) C(x)^{-1}$,
expressing the fact that all renormalization representatives are conjugate,
and in fact any element of the conjugacy class of $A^{(n)}$ arises as a
renormalization representative.

Now, if $A$ is $C^r$, $1 \leq r \leq \infty$ or $r=\omega$,
the normalizing maps may be chosen to be $C^r$ as
well (\cite {AK}, Lemma 4.1).
Hence we may restrict considerations to renormalization
representatives obtained by the use of a $C^r$ normalizing map, which we
call $C^r$-renormalization representatives.  Such $C^r$-renormalization
representatives are defined up to $C^r$-conjugacy.

The dynamics of $(\alpha,A)$ and of its renormalization representatives are
of course intimately related.  For instance:

\begin{prop} \label {rencr}

If a $C^r$-renormalization representative
$(\alpha_n,A^{(n)})$ is $C^r$ conjugate to rotations, then $(\alpha,A)$ is
$C^r$ conjugate to rotations.

\end{prop}

\begin{pf}

Let $B^{(n)}$ be a $C^r$-normalizing map for $(1,A^{(n,0)})$ such that we
have, for every $x \in \R$,
$B^{(n)}(x+\alpha_n)
A^{(n,1)} B^{(n)}(x)^{-1} \in \SO(2,\R)$, and let
$B'(x_*+\beta_{n-1} x)=B^{(n)}(x)$.  Note that $A^{(n,1)}_{q_{n-1}}(x+q_n)
A^{(n,0)}_{q_n}(x)=\id$ so that $B^{(n)}(x+q_n+\alpha_n q_{n-1})
B^{(n)}(x)^{-1} \in \SO(2,\R)$ for every $x \in \R$.
Writing $\tilde A^{(0)}(x)=B'(x+1) B'(x)^{-1}$ and using that that $\frac
{1} {\beta_{n-1}}=q_n+\alpha_n q_{n-1}$, we see that $\tilde A^{(0)}(x) \in
\SO(2,\R)$ for every $x \in \R$.  An analogous argument shows that
$\tilde A^{(1)}(x)=B'(x+\alpha) A(x) B'(x)^{-1} \in \SO(2,\R)$ for
every $x \in \R$.  As remarked in the beginning of the proof of Lemma 4.4 of
\cite {AK}, a simpler version of Lemma 4.1 of \cite {AK} shows the existence
of an $\SO(2,\R)$-valued $C^r$-normalizing map $\tilde B$
for $(1,\tilde A^{(0)})$.  Then $B(x)=\tilde B(x) B'(x)$ is
$1$-periodic and $B(x+\alpha) A(x) B(x)^{-1} \in \SO(2,\R)$ for every $x \in
\R$.
\end{pf}


\subsection{Convergence of renormalization}

A weak version of convergence of renormalization can be stated as follows:

\begin{thm} \label {weak}

Let $(\alpha,A)$ be a $C^r$ cocycle,
$1 \leq r \leq \infty$ or
$r=\omega$, over an irrational rotation, and let $\deg$ be the
topological degree of $A$.
If $(\alpha,A)$ is $L^2$-conjugate to rotations then
there exist a sequence of $C^r$-renormalization representatives
$(\alpha_n,A^{(n)})$ and $\theta_n \in \R$,
such that $R_{-\theta_n-(-1)^n \deg x} A^{(n)}(x) \to \id$ in $C^r$.\footnote{In fact, as the proof will show the convergence holds uniformly on any compact subsets of larger and larger complex strips.}

\end{thm}

\noindent{\it Proof of Theorem \ref {globall2}.}
If $(\alpha,A)$ is non-homotopic to a constant, then $\deg \neq 0$.  By
Theorem \ref {weak}, it admits a monotonic $C^r$-renormalization
representative, which is $C^r$-conjugate to rotations by Theorem \ref
{rigidity}.  By Proposition \ref {rencr}, $(\alpha,A)$ is
$C^r$-conjugate to rotations as well.
\qed

We call the convergence given by Theorem \ref {weak}
weak because it does not say anything
about the normalizing map leading to the ``nice''
renormalization representative.  The strong form of convergence is the
following:

\begin{thm} \label {strong}

Let $(\alpha,A)$ be a $C^r$ cocycle, $1 \leq r \leq \infty$ or
$r=\omega$, over an irrational rotation.
If $(\alpha,A)$ is $L^2$-conjugate to rotations then for almost every $x_*
\in \R/\Z$ there exists $B(x_*) \in \SL(2,\R)$, and a sequence of affine
functions with bounded linear coefficients $\phi^{(n,0)},\phi^{(n,1)}:\R \to \R$
such that
\be
R_{-\phi^{(n,0)}(x)}
B(x_{*}) A^{(n,0)}(x) B^{-1}(x_{*}) \to \id
\ee
and
\be
R_{-\phi^{(n,1)}(x)} B(x_{*}) A^{(n,1)}(x) B^{-1}(x_{*}) \to \id,
\ee
in $C^r$.

\end{thm}

In \cite {AK} is is shown that if $A$ is $C^r$,
then there exists $x_*$
and $B(x_{*})\in \SL(2,\R)$ such that $B(x_{*}) A^{(n,i)}(x) B^{-1}(x_{*})$, $i=0,1$, approaches
$\SO(2,\R)$-valued functions in the $C^r$ topology for $r=\infty,\omega$,
or $C^{r-1}$ if $1 \leq r<\infty$.  While computations in \cite {AK} are
``local'', the more precise version obtained, based on the recent work \cite
{A}, takes into account global aspects of the (asymptotically) holomorphic
extensions of matrix products.  This complex variables proof turns out to
be simpler and more powerful than our original real variables approach,
which shows that if $A$ is $C^1$ then the oscillations of the
derivative of $B(x_{*}) A^{(n,i)}(x) B^{-1}(x_{*})$ become less pronounced as $n \to
\infty$ (due to cancellations appearing through the Ergodic Theorem).

We will prove Theorem \ref {strong}
in the next section.  For the moment, we will just
relate it to Theorem \ref {weak}.

\noindent{\it Proof of Theorem \ref {weak}.}

Let $B(x_{*})$, $\phi^{(n,0)}=a_{n,0} x+b_{n,0}$ and
$\phi^{(n,1)}=a_{n,1} x+b_{n,1}$ be as in Theorem \ref {strong}.
Let $n$ be large and let $\tilde B(x)=R_{-(a_{n,0}\frac {x^2-x} {2}+b_{n,0} x)} B(x_{*})$.
Then $\tilde A(x)=\tilde B(x+1) A^{(n,0)}(x) \tilde
B(x)^{-1}$ is $C^r$-close to the identity and
$\tilde B(x+\alpha_n) A^{(n,1)}(x) \tilde
B(x)^{-1}$ is $C^r$-close to $R_{\psi^{(n)}(x)}$, where
$\psi^{(n)}(x)=(-a_{n,0} \alpha_n+a_{n,1})
x+(b_{n,1}-\alpha_{n}b_{n,0}-a_{n,0} \frac {\alpha_n^2-\alpha_n} {2}
)$.
By Lemma 4.1 of \cite {AK},
there exists $C \in C^r(\R,\SL(2,\R))$ which is $C^r$-close to the identity
such that $C(x+1) \tilde A(x) C(x)^{-1}=\id$.  Set $B^{(n)}=C \tilde B$. 
Then $B^{(n)}$ is a $C^r$ normalizing map for $A^{(n,0)}$ and
$A^{(n)}(x)=B^{(n)}(x+\alpha_n) A^{(n,1)} B^{(n)}(x)^{-1}$ is a $C^r$
renormalization representative close to $R_{\psi_n(x)}$.  Since
$(\alpha_n,A^{(n)})$ is a renormalization representative of $(\alpha,A)$,
the topological degree of
$A^{(n)}$ is $(-1)^n \deg$ (compute directly the degree of an $n$-th
renormalization representative of $(\alpha,R_{\deg x})$, which will be
automatically homotopic to $(\alpha,A^{(n)})$, or
see \cite {AK}, Appendix A).  Thus the linear
coefficient of $\psi_n$ must be close to $(-1)^n \deg$ and $A^{(n)}(x)$ must
be $C^r$-close $R_{\theta_n+(-1)^n \deg x}$ for some $\theta_n \in \R$.
\qed

\subsection{Proof of Theorem \ref {strong}}

The complex variables proof given below follows basically \cite {A},
which uses ``renormalization in parameter space'' as an approach to the
local distribution of zeros of orthogonal polynomials, originally treated
in \cite {ALS} with a different technique.  We translate the argument of
\cite {A} to the usual renormalization operator in the analytic case, and
then use asymptotically holomorphic extensions to address the non-analytic
case.

We consider first the analytic case.  Let $B:\R/\Z \to \SL(2,\R)$ be a
measurable map with $\|B\|^2 \in L^2$ and for any $x \in \R/\Z$,  $B(x+\alpha) A(x) B(x)^{-1} \in
\SO(2,\R)$.  Let $S(x)=\sup_{n \geq 1} \frac {1} {n} \sum_{k=0}^{n-1}
\|B(x+k\alpha)\|^2$, which is finite almost everywhere by the Maximal
Ergodic Theorem.
Assume that $A$ has a holomorphic extension which is
Lipschitz in $\Omega_\delta$.

\begin{lemma} \label {mproducts}

There exists $C>0$ such that if $x_0 \in \R/\Z$
then for every $k \geq 1$ and $z \in \Omega_\delta$ we have
\be
\|A_k(x_0)^{-1}(A_k(x)-A_k(x_0))\| \leq e^{C \|B(x_0\|^2
S(x_0) k |x-x_0|}-1.
\ee

\end{lemma}

\begin{pf}

The proof is the same as the proof of
Lemma 3.1 of \cite {AK} (there only the case $x \in \R/\Z$
is considered, but the proof works equally well for the complex extension).
\end{pf}

Suppose now that $x_*$ is a measurable continuity point of $S$ and $B$
(this means that $x_*$ is a Lebesgue density point of
$\{|S(x)-S(x_*)|<\epsilon\}$ and of $\{\|B(x)-B(x_*)\|<\epsilon\}$ for
every $\epsilon>0$).  Then we get the estimate
\be
\|A_{(-1)^n q_n}(x)\| \leq \inf_{x_0'-x_* \in [-\frac {d} {q_n},+\frac {d}
{q_n}]} C(x_*) e^{C(x_*) q_n |x-x_0'|},
\ee
for every $d>0$, as long as $n>n_0(d)$.
The argument is as in Lemma 3.3 of \cite {AK}: if $n$ is large, the
measurable continuity hypothesis implies that for every $x_0' \in
[x_*-\frac {d} {q_n},x_*+\frac {d} {q_n}]$
we can locate $x_0$ with $|x_0'-x_0| \leq \frac {1} {q_n}$
and such that $B(x_0),B(x_0+\beta_n)$ are close to $B(x_*)$ and
$S(x_0),S(x_0+\beta_n)$ are close to $S(x_*)$, and then apply
Lemma \ref {mproducts}
to estimate either $\|A_{q_n}(x_0)^{-1} (A_{q_n}(x)-A_{q_n(x_0)}\|$
(if $n$ is even), or
$\|A_{-q_n}(x_0) (A_{-q_n}(x)^{-1}-A_{-q_n}(x_0)^{-1})\|$
(if $n$ is odd), using also the bound $\|A_{(-1)^n q_n}(x_0)\| \leq
\|B(x_0)\| \|B(x_0+\beta_{n-1})\|$.

This estimate implies, since $q_{n-1}<q_n<\beta_{n-1}^{-1}$,
\be \label{4.7}
\|A^{(n,i)}(x)\| \leq \inf_{x_0 \in [-d,d]} C e^{C |x-x_0|}, \quad i=0,1,
\quad x \in \Omega_{\delta/\beta_{n-1}}, \quad x_0 \in [-d,d].
\ee
It follows that the sequences $A^{(n,i)}$ are precompact in
$C^\omega$, and the limits are entire functions $\tilde A$ with
$\|\tilde A(z)\| \leq C e^{C |\Im z|}$.  We now  show that
$B(x_*) \tilde A(x) B(x_*)^{-1}$ must be of the form
$R_{\tilde \phi(x)}$ with
$\tilde \phi$ affine with bounded linear coefficient.

Indeed, Lemma 3.4 of \cite {AK} shows that limits $\tilde A$ of the
$A^{(n,i)}$ satisfy $B(x_*) \tilde A(x) B(x_*)^{-1} \in \SO(2,\R)$,
$x \in \R$.  It follows that we can
write $\tilde A(z)=B(x_*)^{-1} R_{\tilde \phi(x)} B(x_*)$
for some entire function $\tilde \phi:\C \to \C$, satisfying the
estimate
 $|\Im \tilde \phi(z)| \leq C+C |\Im  z|$ (since $\Im \tilde\phi=0$ on the real axis).
This implies that   $\tilde \phi$ is affine with bounded linear coefficient.

We consider now the $C^r$ case, $1 \leq r<\infty$, since it implies the
$C^\infty$ case.  Consider an $r$-asymptotically holomorphic extension of
$A$ to some $\Omega_\delta$,
and let $x_*$ be selected as in the analytic case.  The asymptotically
holomorphic extension is in particular Lipschitz in $\Omega_\delta$, thus
estimate (\ref{4.7})  still holds.  For $0<\epsilon \leq \delta$,
let us denote by $\| \cdot \|_{C^{r-1}_\epsilon}$ the
$C^{r-1}$ norm of the restriction to $\Omega_\epsilon$ of a function defined
on $\Omega_\delta$.

\begin{lemma}

Suppose that $x_0 \in \R/\Z$ and $k \geq 1$ satisfy
$S(x_0), \|B(x_0)\|, \|B(x_0+k \alpha)\| \leq C_0$.  Then there exists
$C>0$ (depending on $C_0$ and $\|A\|_{C^r_\delta}$, but not on $x_0$),
such that if $z \in \Omega_\epsilon$ with $0<\epsilon \leq \delta$ then
\be
\max_{0 \leq s \leq r-1}
\|D^s \op_z A_k(z)\| \leq C k^r e^{C k |z-x_0|}
\|\op_z A\|_{C^{r-1}_\epsilon}
\ee
($D$ stands for the full derivative).

\end{lemma}

\begin{pf}

The proof is the same as that of
Lemma 3.2 of \cite {AK} which estimates the
real derivatives of matrix products: the consideration of the
complex extension is again harmless, and the incorporation of a $\op_z$ in
the estimates is straightforward.
\end{pf}

By the same measurable continuity argument given above, we obtain
\be
\max_{0 \leq s \leq r-1}
\|D^s \op_z A_{(-1)^n q_n}(z)\| \leq
\inf_{x_0-x_* \in [-\frac {d} {q_n},\frac {d}
{q_n}]} C q_n^r e^{C q_n |z-x_0|}
\|\op_z A\|_{C^{r-1}_\epsilon}, \quad z
\in \Omega_\epsilon,
\ee
which yields
\be
\max_{0 \leq s \leq r-1}
\|D^s \op_z A^{(n,i)}(z)\| \leq
\inf_{x_0 \in [-d,d]} C e^{C |z-x_0|}
\|\op_z A\|_{C^{r-1}_\epsilon}, \quad z
\in \Omega_{\epsilon/\beta_{n-1}},
\ee
This implies that we can write $A^{(n,i)}=A^{(n,i)}_c+A^{(n,i)}_h$ where
each term is defined in an increasing sequence of disks $D_n$ with $\cup
D_n=\C$, $A^{(n,i)}_h$ are matrix valued (not necessarily
$\SL(2,\C)$) holomorphic functions and form precompact sequences
with limits satisfying $\|\tilde A(z)\| \leq C e^{C |z|}$, and 
$A^{(n,i)}_c$ are $C^r$ matrix valued functions with $C^r$ norm going to
$0$.  It follows that $A^{(n,i)}$ are precompact in $C^r$ and the limits are
entire functions (necessarily $\SL(2,\C)$ valued now) satisfying
$\|\tilde A(z)\| \leq C e^{C |\Im z|}$.  By the same argument of the
analytic case, the limits have the form $B(x_*)^{-1} R_{\tilde \phi(z)}
B(x_*)$, where $\tilde \phi$ is affine with  bounded linear coefficient.

\comm{
Now the proof of Lemma 3.3 of \cite {AK} shows that Lemmas \ref {3.1}
and \ref {3.2} of \cite {AK} imply that for all $d>0$ there exists
$n_0(d)>0$ such that for $n>n_0(d)$ we have
\be
\|A_{(-1)^n q_n}(x_*+x)\| \leq \inf_{x_0-x_* \in [-\frac {d} {q_n},\frac {d}
{q_n}]} C e^{C q_n |x-x_0|}, \quad x \in \Omega_1,
\ee
\be
\sup_{0 \leq k \leq r-1}
\|\partial_k \op_z A_{(-1)^n q_n}(x_*+z)\| \leq
\inf_{x_0-x_* \in [-\frac {d} {q_n},\frac {d}
{q_n}]} C \|\op_z A\|_{C^{r-1}_\epsilon}
q_n^r e^{C q_n |x-x_0|}, \quad x
\in \Omega_\epsilon,
\ee
where $\| \cdot \|_{C^{r-1}_\epsilon}$ stands for the $C^{r-1}$ norm of a
function defined in $\Omega_\epsilon$.
It follows that
\be
\|A^{(n,i)}(x)\| \leq \inf_{x_0 \in [-d,d]} C e^{C |x-x_0|}, \quad i=0,1,
\quad x \in \Omega_{1/\beta_{n-1}}, \quad x_0 \in [-d,d].
\ee
\be
\sup_{0 \leq k \leq r-1}
\|\partial_k \op_z A^{(n,i)}(z)\|
\leq \inf_{x_0 \in [-d,d]} C \|\op_z A\|_{C^{r-1}_\epsilon}
e^{C |x-x_0|}, \quad i=0,1,
\quad x \in \Omega_{\epsilon/\beta_{n-1}}.
\ee
This implies that we can write $A^{(n,i)}=A^{(n,i)}_h+A^{(n,i)}_c$ where
each term is defined in an increasing sequence of disks $D_n$ with $\cup
D_n=\C$, $A^{(n,i)}_h$ is holomorphic matrix valued (not necessarily
$\SL(2,\C)$) and form a precompact sequence
with limits satisfying $\|\tilde A(z)\| \leq C e^{C |Im z|}$, and
$A^{(n,i)}_c$ are $C^r$ matrix valued functions with $C^r$ norm going to
$0$.  It follows that $A^{(n,i)}$ are precompact in $C^r$ and the limits are
entire functions (necessarily $\SL(2,\C)$ valued now) satisfying
$\|\tilde A(z)\| \leq C e^{C |\Im z|}$.  By the same argument of the
analytic case, the limits have the form $B(x_*)^{-1} R_{\tilde \phi(z)}
B(x_*)$, where $\tilde \phi$ has bounded linear coefficient.

Fix some $\delta>0$ such that $A$ is Lipschitz in $\Omega_\delta$.

As in the proof of Lemma 3.3 of \cite {AK}, Lemma 3.1 of \cite {AK}
implies that for every $d>0$, there exists $n_0(d) \geq 0$ such that if
$n>n_0(d)$ then
\be
\|A_{(-1)^n q_n}(x_*+x)\| \leq \inf_{x_0-x_* \in [-\frac {d} {q_n},\frac {d}
{q_n}} C e^{C q_n |x-x_0|}, \quad x
\in \Omega_\delta,
\ee
for some $C>0$ independent of $n$.  Since $q_{n-1}<q_n<\beta_{n-1}^{-1}$,
\be
\|A^{(n,i)}(x)\| \leq \inf_{x_0 \in [-d,d]} C e^{C |x-x_0|}, \quad i=0,1,
\quad x \in \Omega_{\delta/\beta_{n-1}}, \quad x_0 \in [-d,d].
\ee
It follows that the sequences $A^{(n,i)}$ are precompact in
$C^\omega$, and the limits are entire functions $\tilde A$ with
$\|\tilde A(z)\| \leq e^{C |\Im z|}$.  Thus it is enough to show that
$B(x_*) \tilde A(x) B(x_*)^{-1}$ must be of the form
$R_{\tilde \phi(x)}$ with
$\tilde \phi$ affine with bounded linear coefficient.

Taking $x_*$ as a measurable continuity point also of $B$, Lemma 3.4 of
\cite {AK} shows that limits $\tilde A$ of the $A^{(n,i)}$ satisfy $B(x_*)
\tilde A(x) B(x_*)^{-1} \in \SO(2,\R)$, $x \in \R$.  It follows that we can
write $\tilde A(z)=B(x_*)^{-1} R_{\tilde \phi(x)} B(x_*)$
for some entire function $\tilde \phi:\C \to \C$, satisfying the
estimate $|\Im \tilde \phi(z)| \leq C+C |\Im z|$.
This implies that $\tilde \phi$ is linear with bounded linear coefficient.

We consider now the $C^r$ case, $1 \leq r<\infty$, since it implies the
$C^\infty$ case.  Consider an $r$-asymptotically holomorphic extension of
$A$, and let $x_*$ be selected as in the analytic case.
Now the proof of Lemma 3.3 of \cite {AK} shows that Lemmas \ref {3.1}
and \ref {3.2} of \cite {AK} imply that for all $d>0$ there exists
$n_0(d)>0$ such that for $n>n_0(d)$ we have
\be
\|A_{(-1)^n q_n}(x_*+x)\| \leq \inf_{x_0-x_* \in [-\frac {d} {q_n},\frac {d}
{q_n}]} C e^{C q_n |x-x_0|}, \quad x \in \Omega_1,
\ee
\be
\sup_{0 \leq k \leq r-1}
\|\partial_k \op_z A_{(-1)^n q_n}(x_*+z)\| \leq
\inf_{x_0-x_* \in [-\frac {d} {q_n},\frac {d}
{q_n}]} C \|\op_z A\|_{C^{r-1}_\epsilon}
q_n^r e^{C q_n |x-x_0|}, \quad x
\in \Omega_\epsilon,
\ee
where $\| \cdot \|_{C^{r-1}_\epsilon}$ stands for the $C^{r-1}$ norm of a
function defined in $\Omega_\epsilon$.
It follows that
\be
\|A^{(n,i)}(x)\| \leq \inf_{x_0 \in [-d,d]} C e^{C |x-x_0|}, \quad i=0,1,
\quad x \in \Omega_{1/\beta_{n-1}}, \quad x_0 \in [-d,d].
\ee
\be
\sup_{0 \leq k \leq r-1}
\|\partial_k \op_z A^{(n,i)}(z)\|
\leq \inf_{x_0 \in [-d,d]} C \|\op_z A\|_{C^{r-1}_\epsilon}
e^{C |x-x_0|}, \quad i=0,1,
\quad x \in \Omega_{\epsilon/\beta_{n-1}}.
\ee
This implies that we can write $A^{(n,i)}=A^{(n,i)}_h+A^{(n,i)}_c$ where
each term is defined in an increasing sequence of disks $D_n$ with $\cup
D_n=\C$, $A^{(n,i)}_h$ is holomorphic matrix valued (not necessarily
$\SL(2,\C)$) and form a precompact sequence
with limits satisfying $\|\tilde A(z)\| \leq C e^{C |Im z|}$, and
$A^{(n,i)}_c$ are $C^r$ matrix valued functions with $C^r$ norm going to
$0$.  It follows that $A^{(n,i)}$ are precompact in $C^r$ and the limits are
entire functions (necessarily $\SL(2,\C)$ valued now) satisfying
$\|\tilde A(z)\| \leq C e^{C |\Im z|}$.  By the same argument of the
analytic case, the limits have the form $B(x_*)^{-1} R_{\tilde \phi(z)}
B(x_*)$, where $\tilde \phi$ has bounded linear coefficient.
}

\comm{
\begin{rem}

Convergence of renormalization has several applications beyond the one
considered here.  See for instance \cite {FK}.

\end{rem}
}

\comm{
The renormalization operator does not actually act on the space of cocycles
over irrational rotations, but in a closely related space.

Let $\Gamma^r$ be the space of triples $(\alpha,\Phi_0,\Phi_1)$ where
$\alpha \in (0,1) \setminus \Q$, $\Phi_0 \in C^r(\R,\SL(2,\R))$ and $\Phi_1
\in C^r(\R,\SL(2,\R))$ such that $\Phi_0(x+

To each
$(\alpha,A)$, $\alpha\in\R$, $A\in C^r(\R/\Z,\SL(2,\R))$ one can
associate an ${\Z}^2$-action $\Phi$ on the set $\Omega^r=\R \times
C^r(\R,\SL(2,\R))$ by letting $\Phi(1,0)=(1,Id)$ and $\Phi(0,1)=(\alpha,A)$.

For any action $\Phi$ one can define new actions $M_\lambda \Phi$
($\lambda>0$),  $N_U\Phi$ ($U\in \GL(2,\Z)$), $T_{x_*}\Phi$ ($x_*\in \R$) by
(the operations $M,N,T$ are named respectively rescaling, translation and
base change):
\be
M_\lambda(\Phi)(n,m)=(\lambda^{-1} \g^\Phi_{n,m},x \mapsto
A^\Phi_{n,m}(\lambda x)).
\ee
\be
T_{x_*}(\Phi)(n,m)=(\g^\Phi_{n,m},x \mapsto A^\Phi_{n,m}(x+x_*)).
\ee
\be
N_U(\Phi)(n,m)=\Phi(n',m'), \quad
\begin{pmatrix} n'\\m' \end{pmatrix}=
U^{-1} \cdot
\begin{pmatrix} n\\m \end{pmatrix},
\ee
where we have denoted $\Phi(n,m)=(\g^\Phi_{n,m},A^\Phi_{n,m})$.
Notice that $M_\lambda M_{\lambda'}=M_{\lambda \lambda'}$,
$T_{x_*} T_{x'_*}=T_{x_*+x'_*}$,
and $N_U N_{U'}=N_{U U'}$ (that is, $M$, $T$, and $N$ are left actions of
$\R^*$, $\R$ and $\GL(2,\Z)$ on $\Lambda^r$).
Moreover, base changes commute with translations and rescalings.

Notice that $C^r(\R,\SL(2,\R))$ acts on $\Omega^r$ by
$\Ad_B(\alpha,A(\cdot))=(\alpha,B(\cdot+\alpha) A(\cdot)
B(\cdot)^{-1})$.  This action extends to an action (still denoted $\Ad_B$)
on $\Lambda^r$.  We will say that $\Phi$ and $\Ad_B(\Phi)$ are
$C^r$-conjugate via $B$.
}

\comm{
We say that $(\alpha,A)$ is premonotonic if there exists an iterate of
$(\alpha,A)$ which is real-analytic conjugate to a monotonic cocycle.

There exist premonotonic cocycles (with any given frequency) which do not
admit any monotonic iterate.

We do not know if a premonotonic cocycle is always real-analytic conjugate
to a monotonic cocycle.

\subsection{Possible extensions of the notion of premonotonicity}

It will be convenient to identify $\R/\Z$ with $\partial \D$ through
$h:\R/\Z \to \partial \D$ given by $h(y)=e^{2 \pi i y}$.
Given a cocycles $(\alpha,A) \in \R \times C^0(\R/\Z,\SL(2,\R))$,
we let $F \equiv
F_{\alpha,A}:\R/\Z \times \R/\Z \to \R/\Z \times \R/\Z$ be its
projective action: $F(x,y)=(x+\alpha,h^{-1}(A(x) \cdot h(y)))$.

Let $\Pi:\R/\Z \times \R/\Z \to \R/\Z$ be the projection in the second
coordinate.
Let us say that a continuous function $u:\R/\Z \times \R/\Z \to \R$ is
transverse (with respect to $F$) if
\be
\Pi DF(x,y)(1,u(x,y))-u(F(x,y)) \neq 0
\ee
for every $(x,y) \in \R/\Z \times \R/\Z$.  We will speak of positively
transverse functions or negatively transverse functions according to the
sign of $\Pi DF(x,y)(1,u(x,y))-u(F(x,y))$.  Transverse functions form an
open set.

We will say that $u:\R/\Z \times \R/\Z \to \R$ is projective if.

We will say that $u:\R/\Z \times \R/\Z \to \R$ is trivial if the vector
field $(1,u)$ is uniquely integrable and its.

For simplicity, we shall restrict ourselves to cocycles which are at least
$C^1$.  In this case, monotonicity can be defined as follows (this is
actually how monotonicity was presented in the introduction): the horizontal
foliation $\{(\cdot,y)\}_{y \in \R/\Z}$ is sent by $F$ into a foliation
which is transverse to the original horizontal foliation.  In other words, a
cocycle is monotonic if the zero function is transverse.

The notion of premonotonicity we gave can also be stated in those terms: a
cocycle is premonotonic if it admits a transverse function which is
projective and trivial.

There are several possible generalizations of premonotonicity in those
lines:
\begin{enumerate}
\item Declare a cocycle $\alpha$-premonotonic if it admits a projective
trivial transverse function (this is just premonotonicity as we defined),
\item Declare a cocycle $\beta$-premonotonic if it admits a
trivial transverse function,
\item Declare a cocycle $\gamma$-premonotonic if it admits a
projective transverse function,
\item Declare a cocycle $\delta$-premonotonic if it admits a
transverse function.
\end{enumerate}

Obviously the last definition includes all others, but it is not clear if
all are equivalent.  All definitions are $C^1$-open.

Let us sketch how the weakes of the definitions above is already enough to
prove all result we obtained for monotonic cocycles.

\subsection{Examples of non-premonotonic cocycles}

Our contruction is based on the following result.

\begin{thm} [Young, \cite {Y1}]

Let $B_t \in
C^1(\R/\Z,\SL(2,\R))$ be a one-parameter family (defined on a neighborhood
of $t=0$), such that $(t,x) \mapsto B_t(x)$ is $C^1$.  Let
$\beta(x,t)=B_t(x)^{-1}(-1)$ and $C=\{x,\, \beta(x,0)=1\}$.  Assume that
$\pa_x \beta(x,t) \neq 0$, $x \in C$ and
$\{\frac {\pa_t \beta(x,t)} {\pa_t \beta(x,t)}\}_{x \in C}$ are all
distinct.  Let $A_{\lambda,t}=\bm\lambda&0\\0&\lambda^{-1}\em B_t$.

Let $\alpha \in \R/\Z$ be a Brjuno number.
Then there exists $\epsilon_0 \equiv \epsilon_0(\alpha)>0$ (small)
and $\lambda_0 \equiv \lambda_0(\alpha,\epsilon_0)>0$ (large) such that
for every $\epsilon<\epsilon_0$ and for every $\lambda>\lambda_0$
there exists $\delta \equiv \delta(\alpha,\epsilon,\lambda)>0$,
$\eta \equiv \eta(\alpha,\lambda)>0$, such
that $\lim_{\epsilon \to 0,\lambda \to \infty}
\delta=0$ and $\lim_{\lambda \to \infty} \eta=0$, and a set
$X \equiv X(\alpha,\epsilon,\lambda) \subset (-\epsilon,\epsilon)$
such that $|X|>2 \epsilon-\eta$ with the following property.
If $t \in X$ then for every $x \in C$, there exists
$c \equiv c(\alpha,\lambda,t,x) \in \R/\Z$, $z \equiv
z(\alpha,\lambda,t,x) \in \partial \D$ such that
$|c-x| \leq \delta$, $|z-1| \leq \delta$ and
for $n \geq 0$ we have
\be
\left \|\left (\prod_{k=n-1}^{0} A_{\lambda,t}(x+k\alpha) \right )
\cdot \bm z_M&1 \em \right \| \leq \lambda^{-2n/3},
\ee
\be
\left \|\left (\prod_{k=-n}^{-1} A_{\lambda,t}(x+k\alpha)^{-1} \right )
\times \bm z_M&1\em \right \| \leq \lambda^{-2n/3}.
\ee

\end{thm}

This result uses an inductive construction inspired by the work of
Benedicks-Carleson \cite {BC}
on H\'enon maps.  The points $c$ that appear in the
description are `critical points''.  In the case of H\'enon maps the
critical set is a Cantor set of tangencies between stable and unstable
manifolds \cite {Y2}.
In the case discussed here the critical set is a finite set of
points displaying coincidence between the stable and unstable directions.
Besides the qualitative aspects of the result, it will be very important
that those results localize very precisely the critical points near the
easily defined set $C$, and gives quantitative estimates for the behavior of
the orbit of the critical points.

\begin{thm}

In the same setting of the previous lemma, assume that $\{\Im \pa_x
\beta(x,t)\}_{x \in C}$ do not have all the same sign.  Then
if $\epsilon$ is sufficiently small and $\lambda$ is
sufficiently big, if $t \in X(\alpha,\epsilon,\lambda)$ then $A_{\lambda,t}$
is not premonotonic.

\end{thm}

\begin{pf}

Let $h:\R/\Z \to \partial \D$ be given by $h(x)=e^{2 \pi i x}$.
Let $F,G,H:\R/\Z \times \R/\Z \to \R/\Z \times \R/\Z$ be given by
$F(x,y)=(x+\alpha,h^{-1}(A_{\lambda,t}(x) \cdot h(x)))$,
$G(x,y)=(x+\alpha,h^{-1}(B_t(x) \cdot h(x)))$, and
$H(x,y)=(x,h^{-1}(\bm\lambda&0\\0&\lambda^{-1}\em \cdot h(x)))$, so that
$F=H \circ G$.  Let also $\Pi:\R/\Z \times \R/\Z \to \R/\Z$
be the projection on the second coordinate.  We have
\be
\Pi \pa_1 F(x,y)=(\Pi \pa_2 H(G(x,y))) (\Pi \pa_1 G(x,y))=
(\Pi \pa_2 F(x,y)) (\Pi \pa_2 G(x,y))^{-1} (\Pi \pa_1 G(x,y)),
\ee
\be
\Pi \pa_1 F^{-1}(F(x,y))=\Pi \pa_1 G^{-1}(G(x,y))=-\frac {\Pi \pa_1 G(x,y)}
{\Pi \pa_2 G(x,y)}.
\ee

By the hypothesis, there
exists $x \in C$ such that $\Pi \pa_1 G(x,0)>0$.
Let $c=c(\alpha,\lambda,t,x)$ and let
$z=z(\alpha,\lambda,t,x)$.  Let $d=h^{-1}(z)$.  Then $\Pi \pa_1 G(c,d)>0$.

Let us consider a continuous function $u:\R/\Z \times \R/\Z \to \R$.
We have
\be
\Pi DF^n(F^{-n}(c,d))(1,u(F^{-n}(c,d)))=
\sum_{k=0}^{n-1} (\Pi \pa_2 F^k(F^{-k}(c,d)))
(\Pi \pa_1 F(F^{-k-1}(c,d)))+\Pi \pa_2 F^n((F^{-n}(c,d))u(F^{-n}(c,d))
\ee
so we have
\be
\gamma^-(c,d) \equiv
\lim_{n \to \infty} \Pi_2 DF^n(F^{-n}(c,d))(1,u(F^{-n}(c,d)))=
\sum_{k=0}^\infty (\Pi \pa_2 F^k(F^{-k}(c,d)))
(\Pi \pa_1 F(F^{-k-1}(c,d))),
\ee
and the series in the right side converges exponentially fast.  We may
rewrite
\be
\gamma^-(c,d)=\sum_{k=0}^\infty
(\Pi \pa_2 F^k(F^{-k-1}(c,d))) (\Pi \pa_2 G(F^{-k-1}(c,d)))^{-1} (\Pi \pa_1
G(F^{-k-1}(c,d))),
\ee
which gives
\be
|\gamma^-(c,d)|=O(\lambda^{-4/3}).
\ee

On the other hand, we have
\be
\Pi DF^{-n}(F^n(c,d))(1,u(F^n(c,d)))=
\sum_{k=0}^{n-1} (\Pi \pa_2 F^{-k}(F^k(c,d)))
(\Pi \pa_1 F^{-1}(F^{k+1}(c,d)))+\Pi \pa_2 F^{-n}(F^n(c,d))u(F^n(c,d))
\ee
so that
\be
\gamma^+(c,d) \equiv
\lim_{n \to \infty} DF^{-n}(F^n(c,d))(1,u(F^n(c,d)))=
\sum_{k=0}^\infty (\Pi \pa_2 F^{-k}(F^k(c,d)))
(\Pi \pa_1 F^{-1}(F^{k+1}(c,d))),
\ee
which we can rewrite
\be
\gamma^+(c,d)=-\sum_{k=0}^\infty (\Pi \pa_2 F^{-k}(F^k(c,d)))
\frac {\Pi \pa_1 G(F^k(c,d))} {\Pi \pa_2 G(F^k(c,d))},
\ee
so that
\be
\left |\gamma^+(c,d)+\frac {\Pi \pa_1 G(c,d)} {\Pi \pa_2 G(c,d)} \right |=
O(\lambda^{-4/3}).
\ee

Those estimates together imply that $\gamma^+(c,d)<\gamma^-(c,d)$.  In
particular, if $n$ is large,
\be
\Pi DF^n(F^{-n}(c,d))(1,u(F^{-n}(c,d)))>
\Pi DF^{-n}(F^n(c,d))(1,u(F^n(c,d))),
\ee
which implies
$\Pi DF^{2n}(F^{-n}(c,d))(1,u(F^{-n}(c,d)))>u(F^n(c,d))$.  Thus $u$ can not
be negatively transverse.

An analogous argument (considering a different critical point)
shows that $u$ can not be positively transverse.  Thus $u$ cannot be
transverse at all.
\end{pf}

\begin{rem}

Let us define the sign of a critical point as
the sign of $\Pi \pa_1 G(c,d)$ in the notation of the proof of the previous
lemma.  In the argument above, we could have used only one critical point
with sign opposite to the degree (the case of degree $0$ being trivial).
But this does not give a result stronger then what
is stated above: indeed, there are always critical points with the same
sign of the degree.  Actually there are always $2 \deg$ critical points
more with the sign of the degree than with opposite sign
(since $x \mapsto \beta(x,0)$ has degree $2d$).

\end{rem}

It would be interesting to know if (in the case of non-zero degree) the
absence of critical points with sign opposite to the degree implies
premonotonicity of the cocycles coming from Young's construction (for
$\lambda$ sufficiently large and $\epsilon$ sufficiently small).

The examples discussed above show that absence of premonotonicity is
non-negligeable in the measure-theoretical sense.  We believe that
premonotonicity is not even dense, but the method above does not answer this
question.

\comm{
\subsection{Examples of cocycles which are not premonotonic}

Let $\gamma:\R/\Z \to \R/\Z$ be $C^1$ of non-zero degree, and assume that
$D\gamma$ changes sign.  Let $A(\theta)=R_{\gamma(\theta)}$.  Then $(0,A)$
is not homotopic to the identity and
not premonotonic, and indeed if $B \in C^1(\R/\Z,\SL(2,\R))$ is
$C^1$-close to $A$, the cocycle $(0,B)$ is not premonotonic.  More
generally, it is easy to construct, for each $\alpha \in \Q$, $C^1$-open
subsets of non-premonotonic, non-homotopic to the identity
cocycles with frequency $\alpha$.

It is harder to give examples of non-premonotonic, non-homotopic to the
identity cocycles with irrational
frequencies.  For instance, any cocycle of rotations with irrational
frequencies is premonotonic (actually it admits a monotonic iterate and is
real-analytic conjugate to a monotonic cocycle), so the example above does
not work.  Indeed, it could be expected that cocycles with irrational
frequency tend eventually to start turning (after a change of coordinates)
in the direction of the degree.  However, the growth of the
degree is only linear, so it might not be enough to overcome the exponential
behavior of cocycles with a positive Lyapunov exponent.  This is the
mechanism exploited in the following.

\begin{thm}

Let $N>0$ be very large and let $\gamma:\R/\Z \to \R/\Z$ be $C^1$ and
of degree $-1$, satisfying $\gamma(\frac {i} {N}+\theta)=\theta$ for $0 \leq
i \leq N-1$ and $0<\theta<\frac {1} {N}-\frac {1} {N^2}$.
Let $A=R_{\gamma(\theta)}\begin{pmatrix}2&0\\0&1/2\end{pmatrix}$.

There exists $\delta>0$ such that if $|\alpha|<\delta$ and $B:\R/\Z \to
\SL(2,\R)$ is $\delta$-close to $A$ in the $C^1$-topology satisfying
$L(\alpha,B)>L(0,A)-\delta$ then $(\alpha,B)$ is not premonotonic.

\end{thm}

We shall use the following criteria.

\begin{lemma}

Let $(\alpha,A) \in (\R \setminus \Q) \times C^1(\R/\Z,\SL(2,\R))$, and let
$F(x,y)$ denote the corresponding projective action.  If there exists
$(x,y) \in \R/\Z \times \R/\Z$ such that
\be
\deg \lim_{n \to \infty} (\pa_1 F^n) F^{-n}(x,y)=-\infty
\ee
then $(\alpha,A)$ is not premonotonic.

\end{lemma}

\begin{pf}

By taking an iterate, we may assume by contradiction that $(\alpha,A)$ is
real-analytic conjugate to a monotonic cocycle.  Let $G$ be the projective
action of the conjugacy.  Then
\be
\pa_1 (G \circ F^n \circ G^{-1})(G(F^{-n}(x,y)))=
\pa_1 G(x,y)+\pa_2 G(x,y) \cdot (\pa_1 F^n)(F^{-n}(x,y))+
\pa_2 G(x,y) \cdot (\pa_2 F^n)(F^{-n}(x,y)) \cdot (\pa_1 G)(F^{-n}(x,y)).
\ee
If $(\alpha,A)$ would be premonotonic then it would be minimal as well

\begin{pf}

To simplify the expressions below, we identify $\R/\Z$ with $\partial \D$ by
the exponential map and consider $(\alpha,B)$ as a map from $\R/\Z \times
\R/\Z$ to itself.

Let us consider a sequence $(\alpha_n,A^{(n)}) \to (0,A)$ in the
$C^1$ topology and satisfying $L(\alpha_n,A^{(n)}) \to L(0,A)$.  Assume that
$(\alpha_n,A^{(n)})$ is premonotonic for every $n$.  Since premonotonicity
is $C^1$-open, we may also assume that $\alpha_n \in \R \setminus \Q$.

Let us denote by $F_n,F:\R/\Z \times \R/\Z \to \R/\Z \times \R/\Z$
the action of $(\alpha_n,A^{(n)})$, $(0,A)$.  To reach a contradiction, it
is enough to find $n$ sufficiently large such that for every $k \geq 1$ (it
is actually sufficient to consider $k$ large)
there exists $(x,y) \in \R/\Z \times \R/\Z$ such that $\partial_1
F^k_n(x,y)>0$.  Notice that
\be
\partial_1 F^k_n(x,y)=\sum_{j=0}^{k-1} (\partial_2 F^{k-j-1}_n)
(F^{j+1}_n(x,y)) \cdot (\partial_1 F_n) (F^j_n(x,y)),
\ee
where $\partial_1$, $\partial_2$ are the partial derivatives with respect
to each coordinate.
Thus the above sum is positive provided the following conditions are
satisfied:
\begin{enumerate}
\item $(\partial_1 F) (F^j_n(x,y))>\frac {1} {10}$, $j=k-C,...,k-1$,
\item $(\partial_2 F) (F^j_n(x,y)) \leq e^{-\frac {1} {10}}$,
$j=k-C,...,k-1$,
\item $(\partial_2 F^{k-j-1}_n)(F^{j+1}_n(x,y)) \leq e^{-\frac {1}
{10} (k-j-1)}$, $j=0,...,k-C-1$,
\end{enumerate}
where $C \geq 1000$ is an upper bound for $2 |\partial_1 F| \geq |\partial_1
F_n|$ and does not depend on $n$.

Let $U=\cup_{i=0}^{N-1} (\frac {i} {N}+\frac {1} {N^2},\frac {i+1} {N}-\frac
{2} {N^2}) \times (-\frac {1} {1000},\frac {1} {1000})$.  Then the first two
conditions are satisfied whenever $F^{k-C}(x,y) \in U$, if $n$ is
sufficiently large.  Thus we only have
to find $(x',y') \in X_n \cap U$ where $X_n$ is the set of points satisfying
\be \label {X}
\sum_{j=1}^l \ln (\partial_2 F_n) F^{-j}_n(x',y') \leq
-\frac {k} {10}, l \geq 1.
\ee

Let $\mu_n$
be a probability measure on $\R/\Z \times \R/\Z$ obtained as a limit of
$\frac {1} {k} \sum_{j=0}^{k-1} (F^j_n)_* ({\mathrm Leb})$.  To conclude the
proof, it is enough to show that
\be
\mu_n(U)+\mu_n(X_n)>1
\ee
for $n$ large enough, since this implies that $X_n \cap U \neq \emptyset$.

The following properties of $\mu_n$ are immediate to verify
\begin{enumerate}
\item The $\mu_n$ are invariant by $F_n$,
\item We have $\lim \int -\ln |\partial_2 F_n| d\mu_n=L(\alpha_n,A^{(n)})$,
\item If $\mu_{n,x}$ is a desintegration of $\mu_n$ ($\int \mu_{n,x}
dx=\mu_n$) then for almost every $x$, $\mu_{n,x}$ is a Dirac mass
supported on the unstable
direction of the Oseledets decomposition corresponding to $x$.
\end{enumerate}

We may assume that $\mu_n$ converges to $\mu$.  Then $\mu$ is invariant
under $F$ and we have
\be
\lim \int -\ln |\partial_2 F| d\mu=L(0,A).
\ee
Thus if $\mu_x$ is a desintegration of $\mu$ then for almost every $x$ such
that $A(x)$ is hyperbolic we must have that $\mu_x$ is a Dirac mass
supported on the unstable eigendirection of $A(x)$.  In particular,
$1-\mu(U) \leq \frac {1} {2000}$.

Let $n$ be large enough so that $1-\mu_n(U) \leq \frac {1} {1000}$.  Then we
have
\be
\int |\ln |\partial_2 F|+\ln 2| d\mu<\frac {1} {100}.
\ee
By the Maximal Ergodic Theorem, the set $Y_n$ of the pairs $(x',y')$ such
that
\be
\sum_{j=1}^l |\ln |(\partial_2 F) (F^{-j}(x',y'))|+\ln 2| \leq \frac {l} {10},
\quad l \geq 1,
\ee
satisfies $\mu_n(Y_n)>\frac {1} {10}$.  Obviously $Y_n \subset X_n$, so
$\mu_n(X_n)+\mu_n(U)>1$ as required.
\end{pf}

Recall that if $A \in C^0(\R/\Z,\SL(2,\R))$ then many cocycles close to
$(0,A)$ have Lyapunov exponent close to
$L(0,A)$: indeed, for every $\delta>0$, there exists $\epsilon>0$,
$0<\rho<\delta$ such that if $|\alpha|<\epsilon$, $\|B-A\|_{C^0}<\epsilon$
then $|\{|\theta|<\rho,\, |L(\alpha,R_\theta B)-L(0,A)|>\delta\}|<\delta
\rho$.  This follows from upper semicontinuity of the Lyapunov exponent and
\cite {AB}.  Moreover, using \cite {BJ} we get the following.

\begin{cor}

The set of $(\alpha,A) \in \R \times
C^\omega(\R/\Z,\SL(2,\R))$ which are non-homotopic to the identity and not
premonotonic has non-empty interior.

\end{cor}

One can obtain examples of non-premonotonic cocycles with any given
frequency by renormalizing the example discussed here.
}

\begin{problem}

Let $(\alpha,A) \in (\R \setminus \Q) \times C^r(\R/\Z,\SL(2,\R))$, $r \geq
\mathrm {Lip}$.
Does $L(\alpha,A)=0$ imply that $(\alpha,A)$ is premonotonic?

\end{problem}

A positive answer to this problem would show that, for irrational
frequencies, the only obstruction to non-uniform hyperbolicity (in $C^r$,
$r>2$) is to be smoothly conjugated to a cocycle of rotations.

\comm{
\section{Premonotonic cocycles}

We say that $(\alpha,A)$ is premonotonic if there exists $n>0$ such that
$x \mapsto A_n(x)$ is monotonic.

\begin{lemma}

Let $(\alpha,A) \in (\R \setminus \Q) \times C^1(\R/\Z,\SO(2,\R))$.
Then there exists $B \in C^\omega(\R/\Z,\SO(2,\R))$ such that
$B(x+\alpha) A(x) B(x)^{-1}$ is monotonic.

\end{lemma}

\begin{lemma}

Let $(\alpha,A) \in \R \times C^1(\R/\Z,\SL(2,\R))$ be $C^1$ conjugate to a
monotonic cocycle.  Then $(\alpha,A)$ is premonotonic.

\end{lemma}

\begin{problem}

Let $(\alpha,A) \in \R \times C^1(\R/\Z,\SL(2,\R))$ be a premonotonic
cocycle.  Is $(\alpha,A)$ $C^1$ (and automatically $C^\omega$)
conjugate to a monotonic cocycle?

\end{problem}

\subsection{Examples of cocycles which are not premonotonic}

Let $\gamma:\R/\Z \to \R/\Z$ be $C^1$ of non-zero degree, and assume that
$D\gamma$ changes sign.  Let $A(\theta)=R_{\gamma(\theta)}$.  Then $(0,A)$
is not homotopic to the identity and
not premonotonic, and indeed if $B \in C^1(\R/\Z,\SL(2,\R))$ is
$C^1$-close to $A$, the cocycle $(0,B)$ is not premonotonic.  More
generally, it is easy to construct, for each $\alpha \in \Q$, $C^1$-open
subsets of non-premonotonic, non-homotopic to the identity
cocycles with frequency $\alpha$.

It is harder to give examples of non-premonotonic, non-homotopic to the
identity cocycles with irrational
frequencies.  For instance, any cocycle of rotations with irrational
frequencies is premonotonic, so the example above does not work.  Indeed, it
could be expected that cocycles with irrational frequency tend eventually
to start turning in the direction of the degree.  However, the growth of the
degree is only linear, so it might not be enough to overcome the exponential
behavior of cocycles with a positive Lyapunov exponent.  This is the
mechanism exploited in the following.

\begin{thm}

Let $N>0$ be very large and let $\gamma:\R/\Z \to \R/\Z$ be $C^1$ and
of degree $-1$, satisfying $\gamma(\frac {i} {N}+\theta)=\theta$ for $0 \leq
i \leq N-1$ and $0<\theta<\frac {1} {N}-\frac {1} {N^2}$.
Let $A=R_{\gamma(\theta)}\begin{pmatrix}2&0\\0&1/2\end{pmatrix}$.

There exists $\delta>0$ such that if $|\alpha|<\delta$ and $B:\R/\Z \to
\SL(2,\R)$ is $\delta$-close to $A$ in the $C^1$-topology satisfying
$L(\alpha,B)>L(0,A)-\delta$ then $(\alpha,B)$ is not premonotonic.

\end{thm}

\begin{pf}

To simplify the expressions below, we identify $\R/\Z$ with $\partial \D$ by
the exponential map and consider $(\alpha,B)$ as a map from $\R/\Z \times
\R/\Z$ to itself.

Let us consider a sequence $(\alpha_n,A^{(n)}) \to (0,A)$ in the
$C^1$ topology and satisfying $L(\alpha_n,A^{(n)}) \to L(0,A)$.  Assume that
$(\alpha_n,A^{(n)})$ is premonotonic for every $n$.  Since premonotonicity
is $C^1$-open, we may also assume that $\alpha_n \in \R \setminus \Q$.

Let us denote by $F_n,F:\R/\Z \times \R/\Z \to \R/\Z \times \R/\Z$
the action of $(\alpha_n,A^{(n)})$, $(0,A)$.  To reach a contradiction, it
is enough to find $n$ sufficiently large such that for every $k \geq 1$ (it
is actually sufficient to consider $k$ large)
there exists $(x,y) \in \R/\Z \times \R/\Z$ such that $\partial_1
F^k_n(x,y)>0$.  Notice that
\be
\partial_1 F^k_n(x,y)=\sum_{j=0}^{k-1} (\partial_2 F^{k-j-1}_n)
(F^{j+1}_n(x,y)) \cdot (\partial_1 F_n) (F^j_n(x,y)),
\ee
where $\partial_1$, $\partial_2$ are the partial derivatives with respect
to each coordinate.
Thus the above sum is positive provided the following conditions are
satisfied:
\begin{enumerate}
\item $(\partial_1 F) (F^j_n(x,y))>\frac {1} {10}$, $j=k-C,...,k-1$,
\item $(\partial_2 F) (F^j_n(x,y)) \leq e^{-\frac {1} {10}}$,
$j=k-C,...,k-1$,
\item $(\partial_2 F^{k-j-1}_n)(F^{j+1}_n(x,y)) \leq e^{-\frac {1}
{10} (k-j-1)}$, $j=0,...,k-C-1$,
\end{enumerate}
where $C \geq 1000$ is an upper bound for $2 |\partial_1 F| \geq |\partial_1
F_n|$ and does not depend on $n$.

Let $U=\cup_{i=0}^{N-1} (\frac {i} {N}+\frac {1} {N^2},\frac {i+1} {N}-\frac
{2} {N^2}) \times (-\frac {1} {1000},\frac {1} {1000})$.  Then the first two
conditions are satisfied whenever $F^{k-C}(x,y) \in U$, if $n$ is
sufficiently large.  Thus we only have
to find $(x',y') \in X_n \cap U$ where $X_n$ is the set of points satisfying
\be \label {X}
\sum_{j=1}^l \ln (\partial_2 F_n) F^{-j}_n(x',y') \leq
-\frac {k} {10}, l \geq 1.
\ee

Let $\mu_n$
be a probability measure on $\R/\Z \times \R/\Z$ obtained as a limit of
$\frac {1} {k} \sum_{j=0}^{k-1} (F^j_n)_* ({\mathrm Leb})$.  To conclude the
proof, it is enough to show that
\be
\mu_n(U)+\mu_n(X_n)>1
\ee
for $n$ large enough, since this implies that $X_n \cap U \neq \emptyset$.

The following properties of $\mu_n$ are immediate to verify
\begin{enumerate}
\item The $\mu_n$ are invariant by $F_n$,
\item We have $\lim \int -\ln |\partial_2 F_n| d\mu_n=L(\alpha_n,A^{(n)})$,
\item If $\mu_{n,x}$ is a desintegration of $\mu_n$ ($\int \mu_{n,x}
dx=\mu_n$) then for almost every $x$, $\mu_{n,x}$ is a Dirac mass
supported on the unstable
direction of the Oseledets decomposition corresponding to $x$.
\end{enumerate}

We may assume that $\mu_n$ converges to $\mu$.  Then $\mu$ is invariant
under $F$ and we have
\be
\lim \int -\ln |\partial_2 F| d\mu=L(0,A).
\ee
Thus if $\mu_x$ is a desintegration of $\mu$ then for almost every $x$ such
that $A(x)$ is hyperbolic we must have that $\mu_x$ is a Dirac mass
supported on the unstable eigendirection of $A(x)$.  In particular,
$1-\mu(U) \leq \frac {1} {2000}$.

Let $n$ be large enough so that $1-\mu_n(U) \leq \frac {1} {1000}$.  Then we
have
\be
\int |\ln |\partial_2 F|+\ln 2| d\mu<\frac {1} {100}.
\ee
By the Maximal Ergodic Theorem, the set $Y_n$ of the pairs $(x',y')$ such
that
\be
\sum_{j=1}^l |\ln |(\partial_2 F) (F^{-j}(x',y'))|+\ln 2| \leq \frac {l} {10},
\quad l \geq 1,
\ee
satisfies $\mu_n(Y_n)>\frac {1} {10}$.  Obviously $Y_n \subset X_n$, so
$\mu_n(X_n)+\mu_n(U)>1$ as required.
\end{pf}

Recall that if $A \in C^0(\R/\Z,\SL(2,\R))$ then many cocycles close to
$(0,A)$ have Lyapunov exponent close to
$L(0,A)$: indeed, for every $\delta>0$, there exists $\epsilon>0$,
$0<\rho<\delta$ such that if $|\alpha|<\epsilon$, $\|B-A\|_{C^0}<\epsilon$
then $|\{|\theta|<\rho,\, |L(\alpha,R_\theta B)-L(0,A)|>\delta\}|<\delta
\rho$.  This follows from upper semicontinuity of the Lyapunov exponent and
\cite {AB}.  Moreover, using \cite {BJ} we get the following.

\begin{cor}

The set of $(\alpha,A) \in \R \times
C^\omega(\R/\Z,\SL(2,\R))$ which are non-homotopic to the identity and not
premonotonic has non-empty interior.

\end{cor}

One can obtain examples of non-premonotonic cocycles with any given
frequency by renormalizing the example discussed here.

\comm{
Let $\mu_n$ be probability
measures supported in $\R/\Z \times \D$ which are invariant under
$(\alpha_n,A^{(n)})$ and satisfying
\be
L(\alpha_n,A^{(n)})=\int -\ln |\partial_w A^{(n)}(x) \cdot w| d\mu_n(x,w).
\ee
We may assume that $\mu_n \to \mu$.  Then we have
\be
L(\alpha,A)=\int -\ln |\partial_w A(x) \cdot w| d\mu(x,w).
\ee
Let $\mu_x$ be the desintegration of $\mu$: $\mu=\int d\mu_x dx$.  Then
for almost every $x \in \R/\Z$, if $A(x)$ is hyperbolic then
$\mu_x=\delta_{u(x)}$, where $u(x)$ is the
unstable eigendirection of $A(x)$.  Thus
\be
\int |\ln |\partial_w A(x) \cdot w|+\ln 2| d\mu(x,w) \leq \epsilon_N$,
\ee
where $\epsilon_N \to 0$ with $N$.  Then
\be
\int |\ln |\partial_w A^{(n)}(x) \cdot w+\ln 2| d\mu(x,w) \leq 2 \epsilon_N
\ee
for $n$ large.  By the Maximal Ergodic Theorem, the set of $x$ such that
\be
\sum_{j=1}^k -\ln |\partial_w A^{(n)}(x-j \alpha) \cdot w|
d\mu(x,w) \geq k \frac {\ln 2} {2}, \quad k \geq 1
\ee
has Lebesgue measure close to $1$.  Thus we may select $x$ in this set
satisfying $N x \in [1/3,2/3]$.  Now for $k \geq 1$ we have
\be
\int \partial_x (A^{(n)}_k(x-k \alpha) \cdot w) d\mu_{n,x-k\alpha}=\int
\sum_{j=1}^k \partial_w (A^{(n)}_{j-1}(x-(j-1)\alpha) \cdot w)
\partial_x (A^{(n)}(x-j\alpha) \cdot w) d\mu_{n,x-k\alpha}
}

\begin{problem}

Let $(\alpha,A) \in (\R \setminus \Q) \times C^r(\R/\Z,\SL(2,\R))$, $r \geq
\mathrm {Lip}$.
Does $L(\alpha,A)=0$ imply that $(\alpha,A)$ is premonotonic?

\end{problem}

A positive answer to this problem would show that, for irrational
frequencies, the only obstruction to non-uniform hyperbolicity (in $C^r$,
$r>2$) is to be smoothly conjugated to a cocycle of rotations.
}
}

\comm{
\subsubsection{Codimension}

Let us denote by $M^r_{\deg} \subset C^r(\R/\Z,\SL(2,\R))$ the open set of
monotonic cocycles of degree $\deg$.  Here and in what follows, we shall
assume that $r \geq 6$, so that $A \mapsto L(\alpha,A)$ is a $C^2$ function
of $A \in M^r$.

In this section we are going to show that for each $\alpha \in \R$,
$\{L(\alpha,A)=0\}$ is contained on a smooth submanifold of $M^r$ of
codimension $4 \deg$.  Since $L \geq 0$, the
derivative of $A \mapsto L(\alpha,A)$ must vanish at $\{L(\alpha,A)=0\}$. 
Thus, in order to estimate the codimension of $\{L(\alpha,A)=0\}$, it is
enough to estimate the rank of the second derivative of $A \mapsto
L(\alpha,A)$ at any $A \in M^r$ satisfying $L(\alpha,A)=0$.

We first deal with the case where $\alpha=0$.
It is easy to see that if $A \in M^r$ then $L(0,A)=0$ if and only if each of
the equations $A(x)=\id$ and $A(x)=-\id$ has $\deg$ solutions.  An obvious
estimate shows that the second derivative of $A \mapsto L(\alpha,A)$ has
rank $4 \deg$ at such a cocycle.  Even without using this fact, this
characterization clearly defines $\{L(0,A)=0\}$ as a $C^{r-2}$ smooth
submanifold of $M^r$ of codimension $4 \deg$.

The case of rational frequencies can be studied similar to the case
$\alpha=0$, so we shall concentrate on $\alpha \in \R \setminus \Q$ from now
on.

\begin{thm}

Let $(\alpha,A) \in C^r(\R/\Z,\SL(2,\R))$ be conjugated to a cocycle of
rotations.  Then $D^L$ has rank $4n$.

\end{thm}

\begin{pf}

The result is obvious if $\alpha \in \Q$.

It is enough to consider the case where $A(x)=R_{n x}$.  If $\alpha$ is
close to $0$ then $D^2 L$ has rank $4n$.  The result follows by
renormalization.
\end{pf}
}

\comm{
\subsection{Reduction to the case of premonotonic cocycles}

\begin{thm}

Let $(\alpha,A) \in (\R \setminus \Q) \times
C^2(\R/\Z,\SL(2,\R))$ be $L^2$-conjugate to rotations. 
Then $(\alpha,A)$ is premonotonic.

\end{thm}

\begin{cor}

Let $A_\theta \in C^r(\R/\Z,\SL(2,\R))$, $r=\omega,\infty$,
be $C^{2+\epsilon}$ and monotonic in $\theta$.  For almost every $\theta$,
either $L(\alpha,A)>0$ or $(\alpha,A)$ is $C^r$-conjugate to a
cocycle of rotations.

\end{cor}

\begin{cor}

Let $A_\theta \in C^{r+1+\epsilon}(\R/\Z,\SL(2,\R))$, $1 \leq r<\infty$,
be $C^{2+\epsilon}$ and monotonic in $\theta$.  For almost every $\theta$,
either $L(\alpha,A)>0$ or $(\alpha,A)$ is $C^r$-conjugate to a cocycle of
rotations.

\end{cor}
}

\comm{

In particular, if
$A \in C^3(\R/\Z,\SL(2,\R))$ is monotonic, and $\alpha \in
\R$ then $L(\alpha,A)=0$ implies that
$(\alpha,A)$ is $L^2$-conjugated to a cocycle of rotations.  However, this
result holds in slightly better regularity.

\begin{lemma}

Let us consider a family
$\theta \mapsto A_\theta \in C^0(\R/\Z,\SL(2,\R))$, and let $\alpha \in
\R$.  If $A_{\theta_0}$ is monotonic and
\be
\lim_{\theta \to \theta_0} \frac {\|A_\theta-A_{\theta_0}\|}
{|\theta-\theta_0|} \leq K<\infty,
\ee
then
\be
\lim_{\theta \to \theta_0}
\frac {\rho_{\alpha,A_\theta}-\rho_{\alpha,A_{\theta_0}}} {\theta-\theta_0}
\leq K'<\infty,
\ee
where $K'$ only depends on $K$,
$\|A_{\theta_0}\|_{C^0}$, the monotonicity constant of $A_{\theta_0}$ and
its degree.

\end{lemma}

\begin{cor}

Let $A \in C^0(\R/\Z,\SL(2,\R))$ be monotonic and let $\alpha \in \R$.
Then $\theta \mapsto
\rho_{\alpha,R_\theta A}$ is $K$-Lipschitz where $K$ only depends on
$\|A\|_{C^0}$ and the degree of $A$.

\end{cor}

\begin{rem}

This result obviously generalizes to families of (uniformly) monotonic
cocycles which are monotonic and Lipschitz in $\theta$.

\end{rem}

\begin{cor}

Let $A \in C^0(\R/\Z,\SL(2,\R))$ be monotonic and let $\alpha \in \R$.  If
$L(\alpha,A)=0$ then $(\alpha,A)$ is $L^2$-conjugated to a cocycle of
rotations.

\end{cor}

We can also get some consequences for the regularity of the Lyapunov
exponent.

\begin{cor}

Let $A \in C^0(\R/\Z,\SL(2,\R))$ be monotonic and let $\alpha \in \R$.  Then
$\theta \mapsto L(\alpha,R_\theta A)$ belongs to BMO.

\end{cor}

\begin{rem}

This result generalizes to families of (uniformly) monotonic cocycles which
are $C^3$ and monotonic in $\theta$.

\end{rem}

\begin{cor}

The Lyapunov exponent is a continuous function in the space of
$\epsilon$-monotonic cocycles in $\R \times C^0(\R/\Z,\SL(2,\R))$ 
for every $\epsilon>0$.  In particular, it is a continuous function in the
space of monotonic cocycles in $\R \times C^r(\R/\Z,\SL(2,\R))$, $r \geq
{\mathrm Lip}$.

\end{cor}

We will later get better properties of the Lyapunov exponent of monotonic
cocycles by imposing stronger regularity.
}

\comm{

By Fatou's Lemma, for every $x$ and for almost every $\sigma \in \R/\Z$,
$m(\sigma+it,x)$ converges to some function $m(\sigma,x)$ satisfying
\be
\liminf \frac {1} {t} L(\alpha,A_{\sigma+it}) \geq \frac {\epsilon'}
{2} \int_{\R/\Z} \frac {1} {1-|m(x)|^2} dx,
\ee
provided $\liminf \frac {1} {t} L(\alpha,A_{\sigma+it})$ exists.  But this
happens for almost every $\sigma$ with $L(\alpha,A_\sigma)=0$.
\end{pf}

\begin{prop}

Let $(\alpha_t,A_t) \in \R \times C^0(\R/\Z,\Upsilon)$, $t \in [0,\delta)$
be continuous.  Assume that there exists a sequence $t_k \to 0$ and
measurable functions
$x \mapsto m_k(x) \in \D$ such that $A_{t_k}(x) \cdot
m_k(x)=m_k(x+\alpha_k)$ such that
\be
\int_{\R/\Z} \liminf \frac {1} {1-|m_k(x)|^2} dx<\infty.
\ee
Then there exists a measurable function $x \mapsto m(x) \in \D$ satisfying
$m(x+\alpha)=A(x) \cdot m(x)$ and
\be
\int_{\R/\Z} \frac {1} {1-|m(x)|^2} dx \leq
\int_{\R/\Z} \liminf \frac {1} {1-|m_k(x)|^2} dx.
\ee

\end{prop}

\be
\frac {1} {s_k} \liminf_{k \to \infty} \frac {1} {|m(t_k,x)|^2}<\infty
\ee

Since $A_\theta$ is monotonic, the Cauchy Riemann equations imply that for
every $0<\epsilon'<\epsilon$, there exists $\delta>0$ such that
the analytic extension of $A_\theta$ to $\theta
\in \Omega_\delta$ satisfies
$QA_\theta Q^{-1} \cdot \D \subset \D_{e^{-2 \epsilon' \Im(\theta)}}$.
Using the analiticity in $\theta$ of $A_\theta$, we conclude that
$\theta \mapsto \zeta_{A_\theta}$ is a holomorphic function
$\Omega_\delta \to \H/\Z$, whose imaginary part is continuous up to $\R/\Z$.
In particular, for almost every $\sigma \in
\R/\Z$,
\be
\Im(\zeta_{A_{\sigma+it}})=\lim(\Im(\zeta_{A_\sigma})+t \frac {d} {d\sigma}
\rho_{A_\sigma}+o(t).
\ee
Since the Lyapunov exponent is upper semicontinuous, if we know additionally
that $L(\sigma)=0$, we have
\be
\lim \frac {1} {2 \pi t} L(\sigma+it)=\frac {d} {dt} \rho_{A_\sigma},
\ee
and the result follows (since $L(\sigma+it) \geq \epsilon' t$).
\end{pf}

By Oseledets Theorem, for almost every $x \in \R/\Z$ and for all but at most
onevery $z \in
\overline \D$,\if $\mu$ is a probability measure on $\R/\Z \times
\D$ which projects on Lebesgue measure, then
\be
\frac {-1} {2\pi}=\int_{\R/\Z \times \D} \Im(\chi_n(x,z)) d\mu,
\ee
unless $\mu$ gives positive weight to $\R/\Z \times 

Let $\epsilon>0$.  Define $\Upsilon$ (respectively, $\Upsilon_\epsilon$)
as the space of $\SL(2,\C)$ matrices $A$ such that
$A \cdot \D \subset \D$ (respectively, $A \cdot \D \subset
\D_{e^{-\epsilon}}$).  Let $\Xi=Q\Upsilon Q^{-1}$,
$\Xi_\epsilon=Q\Upsilon_\epsilon Q^{-1}$.

For $\delta>0$, let $\Omega_\delta=\{z \in \C/\Z,\,
0 \leq |\Im(z)|<\delta\}.

Let $U$ be a connected open subset of $\C$ (or $\C/\Z$).
Let us consider a continuous family
$A_\theta \in C^0(\R/\Z, \SL(2,\C))$, $\theta \in U$ such that
$Q A_\theta Q^{-1} \cdot \overline \D \subset \D$.

Given $\alpha \in \R$, there exists a unique continuous function $m:\U
\times \R/\Z \to \D$ such that
\be
m(\lambda,x+\alpha)=Q A_\theta Q^{-1} \cdot m(\lambda,x)
\ee
(this follows from the Schwarz Lemma).

Define a continuous map $\tau:U \times \R/\Z \to \C \setminus \{0\}$ by
\be
Q A_\theta Q^{-1} \cdot
\begin{pmatrix}m(\lambda,x)\\1\end{pmatrix}=\tau(\lambda,x)
\begin{pmatrix}m(\lambda,x+\alpha)\\1\end{pmatrix}.
\ee
By the Schwarz Lemma again, for every $\lambda$,
\be
L(\lambda)=\lim_{n \to \infty} \left
\|\Prod_{k=0}^{n-1} A(x+k \alpha) \right \|=
-2 \int \ln |\tau(\lambda,x)| dx>0,
\ee

Let $\hat \tau:U \times \R \to \C$ be a lift of $\tau$, that is,
$\tau(x)=e^{2 \pi i \hat \tau(x)}$.  Two such lifts differ by a constant
integer.  In particular, $\tau(\lambda,x+1)-\tau(\lambda,x)$ is an integer
$\deg$, which is readily seem to coincide with the degree of the (homotopic)
maps $A_\theta:\R/\Z \to \Xi$.  Define $\zeta:U \to \H/\Z$ by
\be
\zeta(\lambda)=\int_{\R/\Z} \hat (\tau(\lambda,x)-\deg x) dx.
\ee
That $\zeta$ takes values in $\H/\Z$ follows from $\Im(\zeta)=\pi^{-1} L$.

\subsection{Fibered rotation number}
}

\comm{
\subsection{Monotonicity}

Let $I \subset \R$ be an interval.  We say that a continuous function
$f:I \to \R$ is $\epsilon$-monotonic if for every $x \neq x'$ we have
\be
\frac {f(x')-f(x)} {x'-x} \geq \epsilon.
\ee
This definition naturally extends to functions defined on (or taking values
on) $\R/\Z$ (by considering lifts) and on the unit circle $\Sr^1 \subset
\R^2 \equiv \C$ (by considering the identification with $\R/\Z$ given by
$x \mapsto e^{2 \pi i x})$.

We say that a continuous one-parameter family of matrices
$A_\theta \in \SL(2,\R)$ is $\epsilon$-monotonic if,
for every $w \in \R^2 \equiv \C$, the function $\theta \mapsto
\frac {A_\theta \cdot w} {\|A_\theta \cdot w\|}$ is
$\epsilon$-monotonic.

A continuous one-parameter family $A_\theta \in C^0(\R/\Z,\SL(2,\R))$ is
said to be monotonic if for every $x \in \R/\Z$, the family $x \mapsto
A_\theta(x)$ is monotonic.

\subsection{Fibered rotation number}

\begin{thm}

Let $A \in C^0(I,\R/\Z,\SL(2,\R))$ be monotonic, and let $\alpha \in \R$. 
Then $\theta \mapsto \rho(A_\theta)$ is either non-decreasing or
non-incrasing.

\end{thm}

This result can be strenghned if one assumes more regularity.  Although we
will not need the following result in this paper, its simple proof is a good
way to illustrate the use of the complexification machinery.

\begin{thm}

Let $A \in C^3(I,\R/\Z,\SL(2,\R))$ be monotonic, and let
$\alpha \in \R$.  For almost every $\theta \in I$, either
$L(\alpha,A_\theta)>0$ or $\frac {d\rho}
{d\theta}>0$.

\end{thm}

This result is easier to prove when $A \in C^\omega(I,\R/\Z,\SL(2,\R))$,
so we shall start by this case.

\subsection{$L^2$-estimates}

If $I \subset \R$ is an interval, we denote by $C^r(I,\R/\Z,\SL(2,\R))$
the set of continuous families
$\theta \mapsto A_\theta \in C^0(\R/\Z,\SL(2,\R))$ such that
the set $\{\theta \mapsto A_\theta(x)\}_{x \in \R/\Z}$ is a precompact
subset of $C^r(I,\SL(2,\R))$.

\begin{thm}

Let $A \in C^r(I,\R/\Z,\SL(2,\R))$, $r \geq 3$ be monotonic, and let
$\alpha \in \R$.  For almost every $\theta \in I$, either
$L(\alpha,A_\theta)>0$ of $(\alpha,A_\theta)$ is $L^2$-conjugated to a
cocycle of rotations.

\end{thm}
}

\appendix

\section{Conformal barycenter} \label {conformal barycenter}

Let $\MM$ be the set of probability measures on $\D$, and for $\mu \in \MM$,
let $\Phi(\mu)=\int_\D \frac {1} {1-|z|^2} d\mu(z)$.  For $w \in \D$, let
$\Phi_w(\mu)=\Phi(\mu')$ where $\mu'$ is the pushforward of $\mu$ by some
Moebius transformation of $\D$ taking $w$ to $0$.  Notice that
if $\Phi(\mu)<\infty$ then $\Phi_w(\mu)<\infty$ for every $w$.  For every $1
\leq K<\infty$, let $\MM_K=\{\mu \in \MM, \Phi(\mu) \leq K\}$, and let
$\MM_\infty=\cup \MM_K$.  Notice that $\MM_K$ is compact in the weak-*
topology for every $K<\infty$.

The next proposition can be proved using the conformal barycenter of
Douady-Earle \cite {DE}.  The construction is sufficiently simple for us to
give the details here.

\begin{prop}

There exists a Borelian function $\BB:\MM \to \D$, equivariant with
respect to M\"oebius transformations of $\D$ and
such that $\Phi(\delta_{\BB(\mu)}) \leq \Phi(\mu)$.

\end{prop}

\begin{pf}

Following an idea of Yoccoz, let us define a pairing $\D \times \D \to \D$
by setting $z * w$ as the midpoint of the hyperbolic geodesic passing
through $z$ and $w$ if $z \neq w$, and $z * z=z$.  This pairing is
continuous and equivariant, and we have
\be
u_s(z,w) \equiv \Phi_s \left (\frac {\delta_z+\delta_w} {2} \right
)-\Phi_s(\delta_{z*w})\geq 0,
\ee
with equality if and only if $z=w$.
Notice that
\be
u_s(z,s)=(2 \Phi_s(\delta_{z*s})-1) (\Phi_s(\delta_{z*s})-1) \geq
\Phi_s(\delta_{z*s})-1.
\ee
Extend the pairing $*$ to
$\MM \times \MM \to \MM$ linearly.  Thus
\be
\mu * \nu=\int_{\D \times \D} \delta_{z * w} d\mu(z) d\nu(w).
\ee
If $\mu,\nu \in \MM_\infty$ then
\be
u_s(\mu,\nu) \equiv
\Phi_s(\frac {1} {2}(\mu+\nu))-\Phi_s(\mu*\nu)=\int_{\D \times \D}
u_s(z,w) d\mu(z)d\nu(w) \geq 0,
\ee
with equality if and only if $\mu=\nu$ is a Dirac mass.
Notice that $u_s:\MM \times \MM
\to [0,\infty]$ is lower semicontinuous, so if $\mu_k \to \mu$ and
$u_s(\mu_k,\mu_k) \to 0$ then $\mu$ is a Dirac mass.
\comm{If $\mu_k \to \mu$ and $\nu_k \to \nu$ then
\be
\limsup_{k \to \infty}
u_s(\mu_k,\nu_k) \geq \lim_{r \to 1-} \limsup_{k \to \infty}
\int_{\overline \D_r \times \overline \D_r} u_s(z,w) d\mu_k(z) d\nu_k(w)
\geq \lim_{r \to 1-} \int_{\overline \D_r \times
\overline \D_r} u_s(z,w) d\mu(z) d\nu(w)=u_s(\mu,\nu).
\ee
}
If $\mu_k \to \delta_s$ we have
\be
\limsup_{k \to \infty}
u_s(\mu_k,\mu_k) \geq
\limsup_{k \to \infty}
u_s(\mu_k,\delta_s) \geq \limsup_{k \to \infty}
\int_\D \Phi_s(\delta_{z*s})-1 d\mu_k(z),
\ee
\comm{
Lower semicontinuity also implies that if $\mu_k \to \delta_s$ then
\be
\limsup u_s(\mu_k,\mu_k) \leq \limsup u_s(\mu_k,\delta_s)=
\limsup \int_\D \Phi_s(\delta_z)-\Phi_s(\delta_z)^{1/2} d\mu_k(z),
\ee
}
and in particular if additionally
$\lim u_s(\mu_k,\mu_k)=0$ then $\lim \Phi_s(\mu_k*\delta_s)=1$.

Given $\mu \in \MM$, define $\mu^{(k)}$ inductively by $\mu^{(0)}=\mu$ and
$\mu^k=\mu^{(k-1)}*\mu^{(k-1)}$.  If $\mu \in \MM_\infty$ then $\mu^{(k)}
\in \MM_\infty$ and we have
$\Phi(\mu^{(k+1)})=\Phi(\mu^{(k)})-u(\mu^{(k)},\mu^{(k)})$.  Thus
$u_s(\mu^{(k)},\mu^{(k)}) \to 0$, and any limit of $\mu^{(k)}$
(which exists by compactness) must be a Dirac mass.  Moreover, if
$\mu^{(n_k)} \to \delta_s$ then $\Phi_s(\mu^{(n_k)}*\delta_s) \to 1$, so
$\Phi_s(\mu^{(n)}*\delta_s) \to 1$ as well and $\delta_s$ must be the unique limit of
$\mu^{(n)}$.  Now we can set $\BB(\mu)=s$,
which is clearly Borelian.\footnote{Although we do not need this fact,
it is easy to see that
$\BB$ is continuous in each $\MM_K$, $1 \leq K<\infty$.}
\end{pf}

The estimates above allow us to obtain compactness result for invariant
sections of cocycles.  For instance, we have the following.

\begin{prop}

Let $f_k:X \to X$ be a sequence of homeomorphisms of $X$ preserving a
probability measure $\mu$ and converging uniformly to a homeomorphism $f:X
\to X$.  Let $A_k \in C^0(X,\Upsilon)$ be a sequence
converging to $A \in C^0(X,\Upsilon)$.  Assume there exists
measurable $m_k:X \to \D$ satisfying $A_k(x) \cdot
m_k(x)=m_k(f_k(x))$, such that
\be \label {H}
H \equiv \liminf_{K,k \to \infty} \int_X \min \{K,\frac {1}
{1-|m_k(x)|^2} \}d\mu(x)<\infty.
\ee
Then there exists a measurable
$m:X \to \D$ such that $A(x) \cdot m(x)=m(f(x))$
and $\int_X \frac {1} {1-|m(x)|^2} d\mu(x) \leq H$.

\end{prop}

\begin{pf}

Let $X_{K,k}=\{x \in X,\, \frac {1} {1-|m_k(x)|^2}<K\}$, and let
$\nu_{K,k}=\int_{X_{K,k}} \delta_{m_k(x)} d\mu(x)$.  Let $\nu$
be any limit of
$\nu_{K,k}$ along a sequence $K_i \to \infty$, $k_i \to \infty$ attaining
the $\liminf$ in (\ref{H}).  Then $\nu$
is a probability measure which projects onto $\mu$ and
satisfies $\int_{X \times \D} \frac {1} {1-|z|^2} d\nu(x,z) \leq H$. 
Let $\nu_x$, $x \in \R/\Z$ be a desintegration of $\nu$: $\int_{X \times
\D} \phi(x,z) d\nu(x,z)=\int_X (\int_\D \phi(x,z) d\nu_x(z)) d\mu(x)$. 
Then $\nu_{f(x)}$ is the pushforward of $\nu_x$ by $w \mapsto \mA(x) \cdot
w$, and $\nu_x \in
\MM_\infty$ for $\mu$-almost every $x$.  Let $m(x)=\BB(\nu_x)$.  Then
$m(f(x))=\mA(x) \cdot m(x)$ and we have
$\int \frac {1} {1-|m(x)|^2} d\mu(x) \leq \int \int
\frac {1} {1-|z|^2} d\nu_x(z) d\mu(x) \leq H$.
\end{pf}

\section{Transitivity of the projective action} \label {transitivity}

We follow the notation of section \ref {minim}.  Our goal is to show the transitivity of the
projective action of multidimensional quasiperiodic cocycles which are not
homotopic to the identity.  The one-dimensional case was
considered in \cite {KKHO}, and in fact the topological ideas that make the
one-dimensional argument work are easily implemented in the multidimensional
case as well.  Let $f_\alpha:x \mapsto x+\alpha$ be an ergodic translation
in $\R^d/\Z^d$.

\begin{prop} \label {tra}

Let $A \in C^0(\R^d/\Z^d,\SL(2,\R))$ be non-homotopic to the identity.
If $f_\alpha:x \mapsto x+\alpha$ be an ergodic translation
in $\R^d/\Z^d$ then $(f_\alpha,A)$ is transitive on $\R^d/\Z^d \times
\partial \D$.

\end{prop}

\begin{pf}

Up to change of coordinate, we may
assume that $x_1 \mapsto A(x_1,...,x_d)$ has positive degree $\deg \geq 1$.
To prove transitivity, it is enough to show that for any open set $U \subset
\R^d/\Z^d \times \partial \D$, the set $\cup_{k \geq 0}
(f_\alpha,A)^k(U)$ is dense in $\R^d/\Z^d \times \partial \D$.

We will actually show a stronger statement.  Let
$\Pi_1:\R^d/\Z^d \times \partial \D \to \R/\Z$,
$\Pi_2:\R^d/\Z^d \times \partial \D \to \R^{d-1}/\Z^{d-1}$ and
$\Pi_3:\R^d/\Z^d \times \partial \D \to \partial \D$ be given by
$\Pi_1(x_1,...,x_d,z)=x_1$, $\Pi_2(x_1,...,x_d,z)=(x_2,...,x_d)$ and
$\Pi_3(x_1,...,x_d,z)=z$.

Let $0<\epsilon<1/10$.  Let us say that a point $x \in \R^d/\Z^d$ is
$\epsilon$-short if there exists a sequence of paths
$\gamma_n:[0,1] \to \R^d/\Z^d \times \partial \D$
such that $\Pi_1 \circ \gamma_n(t)$ converges uniformly to
$\Pi_1(x)+\epsilon t$, $\Pi_2 \circ \gamma_n(t)$ converges uniformly to
$\Pi_2(x)$, and $\Pi_3((f_\alpha,A)^k(\gamma_n(t)))$ has (algebraic) length at most $2
\pi-1/10$ for every $k \geq 0$.  It is clear that the set of
$\epsilon$-short $x$ is forward invariant and closed, so for each
$\epsilon$, either every point is $\epsilon$-short or no point is
$\epsilon$-short.

If for every $\epsilon>0$ there is no point which
is $\epsilon$-short, then for any
$(x,z)=(x_1,...,x_d,z)$ and for
every $\delta>0$, there exists $k \geq 0$ such that, letting
$J_\delta(x,z)=[x_1,x_1+\delta] \times \{(x_2,...,x_d,z)\}$, we have
$|\Pi_3(f_\alpha,A)^k(J_\delta(x,z))|>2 \pi-\delta$.  It  follows
that for every $\delta_0$, the closure of
$\cup_{k \geq 0} (f_\alpha,A)^k(J_{\delta_0}(x,z))$ contains some circle
$\{y\} \times \partial \D$.  Since this set is also forward invariant, it
must contain also the circles $\{y+l \alpha\} \times \partial \D$ for every
$l \geq 0$, and hence, by minimality of $x \mapsto x+\alpha$,
the whole $\R^d/\Z^d \times \partial \D$.  Since $(x,z)$ and $\delta_0>0$
are arbitrary, transitivity follows.

Assume now that there exists $\epsilon>0$ such that every point is
$\epsilon$-short.  We may assume that $\epsilon=\frac {1} {k}$ for some $k
\geq 2$.  Let $\gamma_{n,i}$, $1 \leq i \leq k$, $n \geq 1$, be the
sequences of paths
associated to $(\frac {i-1} {k},0,...,0)$.  We define a sequence
of paths $\tilde \gamma_n:\R/\Z \to \R^d/\Z^d \times \partial \D$ so
that
\begin{enumerate}
\item $\tilde \gamma_n|[(i-1)/k,(3i-2)/3k]$ is given by
$\tilde \gamma_n(t)=\gamma_{n,i}(3kt-3i+3)$,
\item $\tilde \gamma_n|[(3i-2)/3k,(3i-1)/3k]$ is
such that $\Pi_1 \circ \tilde
\gamma_n$ and $\Pi_2 \circ \tilde \gamma_n$ are constant and
$\Pi_3 \circ \tilde \gamma_n$ is a homeomorphism,
\item the diameter of the image of
$\tilde \gamma_n|[(3i-1)/3k,i/k]$ converges to $0$.
\end{enumerate}

One readily checks that these properties imply that for every $l \geq 0$, if
$n$ is sufficiently large, then $\Pi_1 \circ
(f_\alpha,A)^l \circ \tilde \gamma_n$ has topological degree $1$, $\Pi_2 \circ
(f_\alpha,A)^l \circ \tilde \gamma_n$ is homotopic to a constant and $\Pi_3
\circ (f_\alpha,A)^l \circ \tilde \gamma_n$ has topological degree
$\deg_{l,n}$ satisfying $|\deg_{l,n}| \leq 2k-1$.  But
$\deg_{l+1,n}=\deg_{l,n}+\deg \geq \deg_{l,n}+1$ for every $l$ and $n$,
since $A$ is not homotopic to the identity.  Thus for large $n$ we have both
$\deg_{4k,n}-\deg_{0,n} \geq
4k$ and $|\deg_{4k,n}|,|\deg_{0,n}| \leq 2k-1$, a contradiction.
\comm{
is a homeomorphism onto a
$\delta$-neighborhood of $\Pi_1(x)$,
$\Pi_2 \circ \gamma$ is constant equal to $\Pi_2(x)$ and
$\Pi_3((\alpha,A)^k(\gamma([0,1])))$ has length at most $2 \pi-\epsilon$.

if there exists a compact
connected set $K$ with $\Pi_2(K)=\{\Pi_2(x)\}$ and $\Pi_1(x) \in \inter
\Pi_1(K)$ such that for every $k \geq 0$, $\partial \D
\setminus \Pi_3((\alpha,A)^k(K))$ contains an open
interval of length $\epsilon$.  It is easy to see that the set of
$\epsilon$-short points is forward invariant and closed, so either it is
empty or the whole $\R^d/\Z^d$.

If it is empty for every $\epsilon$,
then for any non-trivial
compact interval $I \subset \R/\Z$ and any $w \in \R^{d-1}/\Z^{d-1} \times
\partial \D$, the set $K=I \times \{w\}$
is such that the closure
$\cup_{k \geq 0} (\alpha,A)^k(K)$ contains a circle $\{x_0\}
\times \partial \D$ for some $x_0 \in \R^d/\Z^d$, and since $x \mapsto
x+\alpha$ is minimal we conclude that
$\cup_{k \geq 0} (\alpha,A)^k(K)$ is dense in $\R^d/\Z^d$.

Assume now that there exists $\epsilon>0$ such that every $x$ is
$\epsilon$-short.  Fix $x_0 \in \R^{d-1}/\Z^{d-1}$.
We can find connected
compact sets $K_1,...,K_n$ with $\Pi_2(K_i)=\{x_0\}$, $1 \leq i \leq n$
and $\Pi_1(\cup_{i=1}^n K_i)=\R/\Z$ such that for every $k \geq 0$ and $1
\leq i \leq n$,
$\R/\Z \setminus \Pi_3((\alpha,A)^k(K_i))$ contains an open
interval of length $\epsilon$.  We may assume that
$\Pi_1(K_i)=[a_i,b_i] \neq \R/\Z$ and
$b_i \in J_{i+1}$, $1 \leq i \leq n-1$, $b_n
\in [a_1,b_1]$, and any point in $\R/\Z$ belongs to at most two $J_i$'s.
Let us consider topological closed
intervals $L_1,...,L_n$ with
$\Pi_1(L_i)=b_i$, $\Pi_2(L_i)=\{x_0\}$, $\Pi_3(L_i)$ of length at most half
of the length of $\partial \D$,
and $K_i \cap L_i \neq \emptyset$, $1 \leq i \leq n$,
$L_i \cap K_{i+1} \neq \emptyset$, $1 \leq i \leq n-1$ and
$L_n \cap K_1 \neq \emptyset$.  Consider a sequence of continuous maps
$\gamma_l:\R/\Z \to \Pi_2^{-1}(x_0)$ such that for $1 \leq j \leq n$,
$\gamma_l([(j-1)/n,(2j-1)/2n])$ is contained in a $1/l$-neighborhood of
$K_j$ and $\gamma_k([(2j-1)/2n,j/n])=L_j$.  Then for $k$ large, all
$\gamma_k$ are homotopic and $\Pi_1 \circ \gamma_k:\R/\Z \to \partial \D$
has degree $1$.

Then $K=\bigcup_{i=1}^n (K_i \cup L_i)$ is a
compact connected set.  It is easy to check that $\Pi_1^*$ is an injective
homomorphism between the Cech cohomology groups (with integer coefficients)
$H^1(\R/\Z)$ and $H^1(K)$, and

and
for any $w \in \partial \D$, we have $K(U \times \{w\})=\R^d/\Z^d \times
\partial \D$.

Let $\Pi_1:\R/\Z \times \partial \D \to \R/\Z$ and $\Pi_2:\R/\Z \times
\partial \D \to \partial \D$ be the coordinate projections.
We claim that for any open interval $J \subset J'$, there exists
$k \geq 0$ such that $\Pi_2((\alpha,A)^k(J \times \{w\}))=\partial
\D$.  To see this, fix
some $N>0$ such that the translates of $J$ by $k\alpha$, $0 \leq k \leq N-1$
cover $\R/\Z$.  Then there are $N$ points $x_0<x_1<...<x_{N-1}<x_0+1$
such that for every $0 \leq i \leq N-1$, $[x_i,x_{i+1}]$ is contained in
$J+k(i)\alpha$ for some $0 \leq k(i) \leq N-1$, where $x_N=x_{N-1}+1$.
Let $\gamma:\R/\Z \to \R/\Z \times \partial \D$ be such that for every
$0 \leq i \leq N-1$, $\gamma|[\frac {2i} {2N},\frac {2i+1} {2N}]$ is an
injective map onto $(\alpha,A)^{k(i)}([x_i-k(i)\alpha,x_{i+1}-k(i)\alpha]
\times \{w\}) \subset (\alpha,A)^k(J \times \{w\})$, and $\gamma|[\frac
{2i+1} {2N},\frac {2i+2} {2N}]$ is an injective map onto an interval
contained in $\{x_{i+1}\} \times \partial \D$.  Since $A$ is not homotopic
to the identity, there exists $k>0$ such that the degree of
$\Pi_2 \circ (\alpha,A)^k \circ \gamma$ is at least $2N+1$.  Since $\Pi_2
\circ (\alpha,A)^k \circ \gamma|[\frac {2i+1} {2N},\frac {2i} {2N}]$ is not
surjective, there exists $0 \leq i \leq N-1$ such that $\Pi_2 \circ
(\alpha,A)^k \circ \gamma|[\frac {2i+1} {2N},\frac {2i} {2N}]$ is
surjective.  Then $\Pi_2((\alpha,A)^{k+k(i)}(J \times \{w\}))=\partial \D$
as required.  This concludes the proof of the claim.

Let $x \in J'$ and for $n>0$ let $x \in J_n \subset J'$ be an interval of
length at most $1/n$.  By the claim, there exists
$k_n \geq 0$ such that $\Pi_2 \circ (\alpha,A)^{k_n} (J_n \times
\{w\})=\partial \D$.  Let $y$ be an accumulation point of $x+k_n \alpha$. 
Then $\{y\} \times \partial \D \subset K(J' \times \{w\})$.
The set $K(J' \times \{w\})$ is forward invariant,
so it must contain $\R/\Z \times \partial \D$.
}
\end{pf}


\end{document}